\documentclass[11pt,a4paper]{article}
\usepackage{bbm}
\usepackage{mathrsfs}
\usepackage{amssymb}
\usepackage{amsfonts}
\usepackage{theorem}
\usepackage{graphicx}
\usepackage{enumerate,color}
\usepackage[mathscr]{euscript}
\usepackage{amsmath,amscd}

\makeatletter
\def\tank#1{\protected@xdef\@thanks{\@thanks
		\protect\footnotetext[0]{#1}}}
\def\bigfoot{
	
	\@footnotetext}
\makeatother

\makeatletter 
\@addtoreset{equation}{section}
\makeatother  

\def\Ex {{\mathbb E}}

\def\R {{\mathbb R}}

\def\P {{\mathbb P}}
\def\Q {{\mathbb Q}}
\def\D {{\mathbf D}}
\def\S {{\mathcal S}}

\def\v{{\mathbf{v}}}
\def\bu{{\mathbf{u}}}

\newcommand{\norm}[1]{\lVert#1\rVert}
\newcommand{\qv}[1]{\langle#1\rangle}
\newcommand{\wb}[1]{\overline{#1}}
\def\wt{\widetilde}
\def\wh{\widehat}
\def\D{\mathbf D}

\def\1{{\mathbf 1}}

\def\<{{\langle}}
\def\>{{\rangle}}
\def\eps{\varepsilon}
\def\wh{\widehat}

\def\qed{{\hfill $\square$ \bigskip}}

\newcommand{\ea}{\end{array}}
\newtheorem{theorem}{Theorem}[section]

\newtheorem{lemma}{Lemma}[section]
\newtheorem{definition}{Definition}[section]
\newtheorem{remark}{Remark}[section]

{\theorembodyfont{\rmfamily}
}

\newenvironment{proof}{Proof. }

\begin{document}
		\title{Asymptotic behavior of 3-D stochastic primitive equations of large-scale moist atmosphere with additive noise}
		\author{{Lidan Wang} \quad \hbox{and} \quad {Guoli Zhou}}
		\date{}
		\maketitle

		\begin{abstract}
			
Using a new and general method, we prove the existence of random attractor for the three dimensional stochastic primitive equations defined on a  manifold $\D\subset\R^3$ improving the existence of weak attractor for the deterministic model. Furthermore, we show the existence of the invariant measure.
		\end{abstract}
		
		

		\bigskip
		\section{Introduction}
		Given the space domain as a manifold $\D=\mathbf S^2\times(0,1)$, where $\mathbf S^2$ is a two-dimensional unit sphere, we consider the following three-dimensional viscous stochastic primitive equations in the pressure coordinate system $(\theta,\varphi,\xi)$:
		\begin{subequations}
			\label{eqn:moist1}
			\begin{align}
				& \partial_t\v+\nabla_\v\v+ w \partial_\xi\v +\frac{f}{R_0} \vec k\times\v +\nabla \Phi+L_{1}\v=\dot{W}_{1},\\
				&\label{eqn:Phi}
				\partial_\xi\Phi+\frac{br_s}{r}(1+aq)T=0,\\
				& \label{eqn:div}
				\text{div } \v+\partial_\xi w=0,\\
				& \partial_t T+\nabla_\v T+w\partial_\xi T-\frac{br_s}{r}(1+aq)w+L_2T=Q_T+\dot{W}_2,\\
				& \partial_tq+\nabla_\v q+w\partial_\xi q+L_3q=Q_q+\dot{W}_3.
			\end{align}
		\end{subequations}
		In this geophysical system, unknown functions are $\v, w, \Phi, T, q$ and the physical meanings are as follows: $(\v,w)=(v_\theta,v_\varphi,w)$ is the 3-D fluid velocity field, with $\v=(v_\theta,v_\varphi)$ being the horizontal velocity and $w$ being the vertical velocity in the pressure coordinate system; $\Phi$ is the geopotential, $T$ is the temperature, and $q$ is the mixing ratio of water vapor in the air.
		\par
		$f=2\cos\theta$ is the Coriolis parameter, $R_0$ is the Rossby number, $\vec k$ is the vertical unit vector. $r$ is the pressure function depending on the variable $\xi$:  $r=(r_s-r_0)\xi+r_0$, where $0<r_0\leq r\leq r_s$. $Q_T,Q_q$ are given functions on $\mathbf S^2\times(0,1)$, $a, b$ are positive constants. The viscosity and the heat diffusion operators $L_{1}, L_{2}, L_3$ are given by
		\begin{equation*}
			L_{i}=-\nu_{i}\Delta-\mu_{i}\partial_{\xi\xi} , \hspace{1cm} i=1,2,3,
		\end{equation*}
		where the positive constants $\nu_{i}, \mu_{i}, i=1,2,3$ are the horizontal and vertical Reynolds numbers. To simplify the notations, we assume
		$ \nu_{i}= \mu_{i}=1,\ \ i=1,2,3.$ The results in this paper are still valid when we consider the general cases.  Note that in this paper we are discussing operators defined on manifolds, and in the pressure coordinate system. The formal definitions of $\nabla_\v\v, \Delta\v, \Delta T, \Delta q, \nabla_\v T, \nabla_\v q, \text{div }\v, \nabla\Phi$ will be given in Section {\ref{sec:preliminary}}. $\dot{W}_{j},j=1,2,3$ are  independent Gaussian white noise processes which
are formally delta correlated in time.
		\par
		The boundary conditions of the system (\ref{eqn:moist1}) are given by
		\begin{subequations}
			\label{eqn:BC}
			\begin{align}
				& \text{On } \xi=1:\ \partial_\xi\v=0, w=0, \partial_\xi T=-\alpha(T-T^*), \partial_\xi q=-\beta(q-q^*),\\
				& \text{On } \xi=0:\ \partial_\xi\v=0, w=0, \partial_\xi T=0, \partial_\xi q=0,
			\end{align}
		\end{subequations}
		where $\alpha, \beta$ are positive constants, $T^*, q^*$ are the given temperature and mixing ratio of water vapor on the surface of the earth, respectively. For simplicity, we assume that $T^*=q^*=0$. One can always homogenize boundary conditions for nonzero $T^*, q^*$ (see \cite{GH}).
		\par
		Now integrating (\ref{eqn:Phi}) and (\ref{eqn:div}) and applying boundary conditions (\ref{eqn:BC}), with $\Phi_s(t;\theta,\varphi)$ being a certain unknown function on the isobaric surface $\xi=1$, one can get
		\begin{equation}
			\Phi(t;\theta,\varphi,\xi)=\Phi_s(t;\theta,\varphi)+\int_{\xi}^1\frac{br_s}{r}(1+aq)Td\xi',
		\end{equation}
		\begin{equation}
			w(t;\theta,\varphi,\xi)=\int_{\xi}^{1}\text{div }\v(t;\theta,\varphi,\xi')d\xi',
		\end{equation}
		\begin{equation}
			\int_0^1\text{div }\v d\xi=0.
		\end{equation}
		In addition, we supply the system with the initial conditions:
		\begin{subequations}
			\label{eqn:IC}
			\begin{align}
				&\v(t_0;\theta,\varphi,\xi)=\v_0(\theta,\varphi,\xi),\\
				&T(t_0;\theta,\varphi,\xi)=T_0(\theta,\varphi,\xi),\\
				&q(t_0;\theta,\varphi,\xi)=q_0(\theta,\varphi,\xi).
			\end{align}
		\end{subequations}
		With all the above discussion, we have the following equivalent formulation for the 3-D stochastic PEs:
		\begin{subequations}
			\label{eqn:moist2}
			\begin{align}
				\label{eqn:v}
				&\partial_t\v+\nabla_\v\v+\bigg(\int_{\xi}^1\text{div }\v(t;\theta,\varphi,\xi')d\xi'\bigg)\partial_\xi\v+\frac{f}{R_0}\vec k\times\v+\nabla\Phi_s\notag\\
				&\hspace{1cm}+\int_{\xi}^1\frac{br_s}{r}\nabla[(1+aq)T]d\xi'-\Delta\v-\partial_{\xi\xi}\v     =\dot{W}_1,\\
				\label{eqn:T}
				& \partial_t T+\nabla_\v T+\bigg(\int_{\xi}^1\text{div }\v(t;\theta,\varphi,\xi')d\xi'\bigg)\partial _\xi T-\frac{br_s}{r}(1+aq)\bigg(\int_{\xi}^1\text{div }\v(t;\theta,\varphi,\xi')d\xi'\bigg)\notag\\
				&\hspace{1cm}-\Delta T-\partial_{\xi\xi}T=Q_T+\dot{W}_2,\\
				\label{eqn:q}
				&\partial_t q+\nabla_\v q+\bigg(\int_{\xi}^1\text{div }\v(t;\theta,\varphi,\xi')d\xi'\bigg)\partial_\xi q-\Delta q-\partial_{\xi\xi}q=Q_q+\dot{W}_3,\\
				&\int_0^1\text{div }\v d\xi=0,\\
				&\text{On } \xi=1\ (r=r_s):\ \partial_\xi\v=0, w=0, \partial_\xi T=-\alpha T, \partial_\xi q=-\beta q,\\
				& \text{On } \xi=0\ (r=r_0):\ \partial_\xi\v=0, w=0, \partial_\xi T=0, \partial_\xi q=0.
			\end{align}
		\end{subequations}
		The primitive equations (PEs) are the basic models to study the mechanism of long-term weather prediction and climate changes, whose mathematical study was initiated by Lions, Teman and Wang (\cite{LTW1}-\cite{LTW4}). This research field has received a wide attention from mathematical community over the last two decades. Taking advantage of the fact that the pressure is essentially two-dimensional in the PEs, Cao and Titi \cite{CT} proved the global results for the existence of strong solutions of
		three-dimensional PEs. Subsequently, Kukavica and Ziane \cite{KZ} developed a different proof which allows one to treat non-rectangular domains as well as different and physically realistic boundary conditions. We refer the reader to the survey
		papers \cite{PTZ, RTT} for further references and background about the deterministic mathematical theory of the PEs.

\par
		Although the PEs express very fundamental laws of physics, the deterministic models are numerically intractable.
		Studies have shown that resolved states are associated with many possible unresolved states.
		This calls for stochastic methods for numerical weather and climate prediction which potentially allow a proper
		representation of the uncertainties, a reduction of systematic biases and improved representation of long-term
		climate variability (see \cite{BSL, MT, ZF}).

		\par
		Despite the developments in the deterministic case, the theory for the
		stochastic PEs remains underdeveloped. Ewald, Petcu, Teman \cite{EPT} and Glatt-Holtz, Ziane \cite{GHZ} considered a two-dimensional stochastic PEs. Then Glatt-Holtz and Temam \cite{GHT1, GT} extended the case to the greater generality of physically relevant boundary conditions and nonlinear multiplicative noise. Following the methods similar to \cite{CT}, Guo and Huang \cite{GH} studied the global well-posedness of the three-dimensional system with an additive noise in the horizontal momentum equations and obtained some kind of weak type compactness properties of the solutions to the stochastic system. Using methods different from \cite{GH},  Debussche, Glatt-Holtz,  Temam and  Ziane considered a three-dimensional system with multiplicative noise. Dong,  Zhai, and   Zhang obtained the large deviation principle, Markov selection and ergodicity for the three dimensional PEs with non-degenerate noises (see \cite{DZZ1,DZZ2,DZZ3}).

        \par
In this paper, we mainly study the existence of random attractor and invariant measure for the stochastic PEs. Note that the definition of  attractors in our paper is different from that in \cite{GH}. The random attractor obtained in our work is $\mathbb{P}$-a.e. $\omega$ compact in $(H^{1}(\D))^{4}$ and attracts any orbit starting from $-\infty$ in the strong topology of $(H^{1}(\D))^{4}$. While the attractor studied in \cite{GH} is not necessarily  a compact subset in  $(H^{1}(\D))^{4}$, and the attractor attracts any orbit in the weak topology of $(H^{1}(\D))^{4}$.

		\par
		
	     Since the uniqueness of the weak solution to the 3D stochastic PEs is still open, we have to choose $(H^{1}(\D))^{4}$ as the phase space to work with. After following the method in \cite{GH11} to prove the global existence of the strong solutions, we will show the continuity of the strong solutions to the 3D stochastic PEs in the space $(H^{1}(\D))^{4}$ with respect to time $t$ as well as with respect to the initial condition $(\v_{0}, T_{0}, q_0).$ Notice that \cite{GH11} only proved the strong solution is Lipschitz continuous in the space $( L^{2}(\D))^{4}$ with respect to the initial data but this is not enough to study the asymptotic behavior in $(H^{1}(\D))^{4}$ considered  here. The first new difficulty arises here in obtaining the regularities of the strong solution about time $t$ and initial condition is that we have no valid boundedness  for the derivatives of the vertical velocity. The second difficulty is that the geometric structure of the manifold is more complicated than the case of $\mathbb{R}^{n}.$ For example, in order to obtain  \textit{a priori} estimates in $(L^{4}(\mathbf{D}))^{4},$ there is no estimate like the following:
   $$|\nabla_{e_\varphi} \v^{3}|\leq C|\nabla_{e_\varphi} \v| |\v|^{2},$$
   where $\v=(v_{1}, v_{2})\in \mathbb{R}^{2}, \v^{3}=(v_{1}^{3}, v_{2}^{3})\in \mathbb{R}^{2}$ and $C$ is a constant.
To overcome the difficulties, we should delicately and carefully use the geometric structure of the manifold involved with the velocity  to obtain the $a\ priori$ estimates.
		Finally,  with the above mentioned $a\ priori$ estimates and continuity properties, we try to find a compact absorbing set in $(H^{1}(\D))^{4}$ to establish the existence of the random attractor which is in fact the most difficult problem for this stochastic moist PEs.  For a bounded domain, the common method is to find an absorbing set for the solutions in the function space with higher regularity than $(H^{1}(\D))^{4}.$  However, this seems to be very complicated for the 3D stochastic PEs.
To overcome the difficulties, we will adopt a new method in this paper. The main idea is that we
		firstly prove that $\mathbb{P}$-a.e. $\omega$ the solution operator $\mathcal S(t,s;\omega)_{t\geq s;\omega\in \Omega}$ of stochastic PEs is compact in the function space $ (H^{1}(\D))^{4}$ for any fixed time $t,s\in \mathbb{R}$. Then by virtue of the regularity of strong solution, we use the solution operator to act on an absorbing ball to construct a compact absorbing ball, which implies the existence of random attractor and invariant measure.
We would like to mention that our method provides a general way for proving the existence of random attractor for common classes of dissipative stochastic partial differential equations and has some advantages over the common method of using compact Sobolev embedding theorem, i.e., if an absorbing ball for the solutions in space $(H^{2}(\D))^{4}$ does exist, our method will then further imply the existence of global random attractor in $(H^{2}(\D))^{4}$.
		\par
		The remaining of this paper is organized as follows. In Section 2, we state some preliminaries and then give the main results, including local and global  existence of solutions, as well as the existence of random attractor and invariant measure,  the proof of global existence is given in Section 3. Since in Section 2 we point out that the existence of the random attractor implies the existence of the invariant measure, in Section $4$ we only study the existence of  random attractor. The $a \ priori$ estimates for the global existence of strong solutions are shown in Appendix the Section $5$. As usual, constant
		$C$ may change from one line to the next, unless we give a special declaration; we denote by $C(a)$ a constant which depends on some parameter $a.$
		
		\section{Preliminaries and main results}
		\label{sec:preliminary}
		In this section, we will give the formal definitions for differential operators in the pressure coordinate system, and the stochastic terms, then reformulate
		this geophysical system into an abstract setting. We will present the main results at the end of this section.
		\par
		Let $e_\theta, e_\varphi, e_\xi$ be the unit vectors in $\theta, \varphi, \xi$ directions of the space domain $\D=\mathbf S^2\times(0,1)$, respectively,
		$$e_\theta=\partial_\theta,\ e_\varphi=\frac{1}{\sin\theta}\partial_\varphi,\ e_\xi=\partial_\xi.$$
		Correspondingly, define the following spaces
		$$L^p(\D):=\{h:\D\to\R, \int_{\D}|h|^p d\D<\infty\},\text{ for }1\leq p<\infty,$$
		$$L^2(T\D|T\mathbf S^2):=\{\v=(v_\theta,v_\varphi):\mathbf D\to T\mathbf S^2, \int_{\D}(|v_\theta|^2+|v_\varphi|^2)d\D<\infty\},$$
		$$C^\infty(\mathbf S^2),\  C^\infty(\D)\text{ are smooth function spaces defined on }\mathbf S^2, \D,\text{ respectively},$$
		$$C^\infty(T\D|T\mathbf S^2):=\{\v=(v_\theta,v_\varphi): v_\theta, v_\varphi\in C^\infty(\D)\}.$$
		$H^m(\D)$ is the Sobolev space of functions which are in $L^2$, with all covariant derivatives with respect to $e_\theta, e_\varphi, e_\xi$ of order $\leq m$, then analogously, we define
		$$H^m(T\D|T\mathbf S^2):=\{\v=(v_\theta,v_\varphi): v_\theta,v_\varphi\in H^m(\D)\}.$$
		In the pressure coordinate system, given $\v=v_\theta e_\theta+v_\varphi e_\varphi, \bu=u_\theta e_\theta+u_\varphi e_\varphi\in C^\infty(T\D|T\mathbf S^2)$, $T,q\in C^\infty(\D)$ and $\Phi_s\in C^\infty(\mathbf S^2)$, we first define the horizontal gradient $\nabla$ for $T, \Phi_s$ on $\mathbf S^2$ as follows:
		\begin{subequations}
			\begin{align}
				&\nabla T=(\partial_\theta T)e_\theta+\Big(\frac{1}{\sin\theta}\partial_\varphi T\Big) e_\varphi,\\
				&\nabla\Phi_s=(\partial_\theta\Phi_s) e_\theta+\Big(\frac{1}{\sin\theta}\partial_\varphi\Phi_s \Big)e_\varphi.
			\end{align}
		\end{subequations}
		We then define the covariant derivative of $\bu, T, q$ with respect to $\v$ as follows:
		\begin{subequations}
			\begin{align}
			\label{eqn:2.2a}
				&\nabla_\v\bu=\bigg(v_\theta\partial_\theta u_\theta+\frac{v_\varphi}{\sin\theta}\partial_\varphi u_\theta-v_\varphi u_\varphi\cot\theta\bigg)e_\theta+\bigg(v_\theta\partial_\theta u_\varphi+\frac{v_\varphi}{\sin\theta}\partial_\varphi u_\varphi+v_\varphi u_\theta\cot\theta\bigg)e_\varphi,\\
				&\nabla_\v T=v_\theta\partial_\theta T+\frac{v_\varphi}{\sin\theta}\partial_\varphi T,\\
				&\nabla_\v q=v_\theta\partial_\theta q+\frac{v_\varphi}{\sin\theta}\partial_\varphi q.
			\end{align}
		\end{subequations}
		Finally, the divergence form of $\v$ and the Laplace-Beltrami operator $\Delta$ for scalar and vector functions are defined as
		\begin{subequations}
			\begin{align}
			\label{eqn:2.3a}
				&\text{div }\v=\text{div}(v_\theta e_\theta+v_\varphi e_\varphi)=\frac{1}{\sin\theta}(\partial_\theta(v_\theta\sin\theta)+\partial_\varphi v_\varphi),\\
				&\Delta T=\text{div}(\nabla T)=\frac{1}{\sin\theta}[\partial_\theta(\sin\theta\partial_\theta T)+\frac{1}{\sin\theta}\partial_{\varphi\varphi}T],\\
				&\Delta q=\text{div}(\nabla q)=\frac{1}{\sin\theta}[\partial_\theta(\sin\theta\partial_\theta q)+\frac{1}{\sin\theta}\partial_{\varphi\varphi}q],\\
				&\Delta\v=\Big(\Delta v_\theta-\frac{2\cos\theta}{\sin^2\theta}\partial_\varphi v_\varphi-\frac{v_\theta}{\sin^2\theta}\Big)e_\theta+\Big(\Delta v_\varphi+\frac{2\cos\theta}{\sin^2\theta}\partial_\varphi v_\theta-\frac{v_\varphi}{\sin^2\theta}\Big)e_\varphi.
			\end{align}
		\end{subequations}
		We then define our working spaces for the stochastic PE system as follows:
		\begin{align*}
			&\mathcal V_1^0:=\left\{\v\in C^\infty(T\D|T\mathbf S^2): \partial_\xi\v\big|_{\xi=0}=\partial_\xi\v\big|_{\xi=1}=0,\int_0^1\text{div }\v d\xi=0\right\},\\
			& \mathcal V_2^0:=\left\{T\in C^\infty(\D): \partial_\xi T\big|_{\xi=0}=0,\ \partial_\xi T\big|_{\xi=1}=-\alpha T\right\},\\
			& \mathcal V_3^0:=\left\{q\in C^\infty(\D): \partial_\xi q\big|_{\xi=0}=0,\ \partial_\xi q\big|_{\xi=1}=-\beta q\right\}.
		\end{align*}
		We denote by $\mathcal V_1, \mathcal V_2$ and $\mathcal V_3$ the closure spaces of $\mathcal V_1^0$ in $H^1(T\D|T\mathbf S^2)$, $\mathcal V_2^0$  and $\mathcal V_3^0$ in $H^1(\D)$ with respect to $H^1$ norm. Also define $\mathcal H_1$ as the closure space of $\mathcal V_1^0$ with respect to $L^2$ norm in $L^2(T\D|T\mathbf S^2)$. Now set
		$$\mathcal V=\mathcal V_1\times\mathcal V_2\times\mathcal V_3,\ \mathcal H=\mathcal H_1\times(L^2(\D))^2.$$
		Let $U:=(\v,T,q), \wt U:=(\wt\v,\wt T,\wt q)\in \mathcal V$ and we equip $\mathcal V$ with the inner product
		\begin{align*}
			&\qv{U,\wt U}_{\mathcal V}:=\qv{\v,\wt\v}_{\mathcal V_1}+\qv{T,\wt T}_{\mathcal V_2}+\qv{q,\wt q}_{\mathcal V_3},\\
			&\qv{\v,\wt\v}_{\mathcal V_1}:=\int_\D (\nabla_{e_\theta}\v\cdot\nabla_{e_\theta}\wt v+\nabla_{e_\varphi}\v\cdot\nabla_{e_\varphi}\wt\v+\partial_\xi\v\partial_\xi\wt\v+\v\cdot\wt\v)d\D,\\
			&\qv{T,\wt T}_{\mathcal V_2}:=\int_\D(\nabla T\cdot\nabla\wt T+\partial_\xi T\partial_\xi\wt T)d\D+\alpha\int_{\mathbf S^2}T|_{\xi=1}\wt T|_{\xi=1}d\mathbf S^2,\\
			&\qv{q,\wt q}_{\mathcal V_3}:=\int_\D(\nabla q\cdot\nabla\wt q+\partial_\xi q\partial_\xi\wt q)d\D+\beta\int_{\mathbf S^2}q|_{\xi=1}\wt q|_{\xi=1}d\mathbf S^2.
		\end{align*}
		Similarly, we equip the Hilbert space $\mathcal H$ with the inner product
		\begin{align*}
			&\qv{U,\wt U}_{\mathcal H}:=\qv{\v,\wt\v}+\qv{T,\wt T}+\qv{q,\wt q},\\
			&\qv{\v,\wt\v}:=\int_{\D}\v\cdot\wt\v d\D,\\
			&\qv{T,\wt T}:=\int_{\D}T\wt Td\D,\\
			&\qv{q,\wt q}:=\int_{\D}q\wt qd\D.
		\end{align*}
		We also denote different norms for $U,\v,T,q$ by
		$$\norm{U}_{1}=\qv{U,U}_{\mathcal V}^{\frac{1}{2}},\ \norm{\v}_{1}=\qv{\v,\v}_{\mathcal V_1}^{\frac{1}{2}},\ \norm{T}_{1}=\qv{T,T}_{\mathcal V_2}^{\frac{1}{2}},\ \norm{q}_{1}=\qv{q,q}_{\mathcal V_3}^{\frac{1}{2}}.$$
		$$|U|_{2}=\qv{U,U}_{\mathcal H}^{\frac{1}{2}},\ |\v|_{2}=\qv{\v,\v}^{\frac{1}{2}},\ |T|_{2}=\qv{T,T}^{\frac{1}{2}},\ |q|_{2}=\qv{q,q}^{\frac{1}{2}}.$$

		We here directly quote Lemma 4.1 in \cite{GH11}, without a proof, to obtain some integral equalities:
		\begin{lemma}
			\label{lemma:ibp}
			Let $\bu=(u_\theta,u_\varphi), \wt\bu=(\wt u_\theta,\wt u_\varphi)\in C^\infty(T\D|T\mathbf S^2)$, and $h\in C^\infty(\mathbf S^2)$, we have
			\begin{enumerate}[\rm i)]
				\item $$\int_{\mathbf S^2}h\text{div }\bu d\mathbf S^2=-\int_{\mathbf S^2}\nabla h\cdot\bu d\mathbf S^2,\hspace{1cm}\text{in particular, }\int_{\mathbf \D}\nabla h\cdot\v d\mathbf S^{2} =0,\text{ for any }\v\in \mathcal V_{1}.$$
				\item $$\int_{\D}(-\Delta\bu)\cdot\wt\bu d\D=\int_{\D}(\nabla_{e_\theta}\bu\cdot\nabla_{e_\theta}\wt\bu+\nabla_{e_\varphi}\bu\cdot\nabla_{e_\varphi}\wt\bu+\bu\cdot\wt\bu)d\D.$$
			\end{enumerate}
		\end{lemma}
In our paper, we will frequently use the following inequalities, so we state them in the following lemmas. For their proof, one can refer to \cite{GH}.
\begin{lemma}
	\label{lemma:2.2}
Let $v\in H^{2}(T\D|T\mathbf S^2), \mu\in H^{1}(T\D|T\mathbf S^2)  \Big{(}\text{ or }\mu\in H^{1}(\D)\Big{)}$ and $\nu\in  L^{2}(T\D|T\mathbf S^2)  \Big{(}\text{ or } \nu \in L^{2}(\D)\Big{)}.$ Then, there exists a positive constant $c$ independent of $v,\mu$ and $\nu$ such that
\begin{align*}
&\Big{|}\langle \Big{(} \int_{\xi}^{1} \mathrm{div} v(t; \theta, \phi, \xi' )d\xi'  \Big{)} \mu,  \nu \rangle\Big{|}\\
\leq& c|\mathrm{div} v|_{2}^{\frac{1}{2}}(|\mathrm{div} v|_{2}^{\frac{1}{2}}+ |\Delta v|_{2}^{\frac{1}{2}} )|\mu|_{2}^{\frac{1}{2}}(|\nabla_{e_{\theta}} \mu|_{2}^{\frac{1}{2}}+|\nabla_{e_{\varphi}} \mu|_{2}^{\frac{1}{2}}   + |\Delta \mu|_{2}^{\frac{1}{2}})|\nu|_{2},\\
\Big{(}\text{ or }\leq& c|\mathrm{div} v|_{2}^{\frac{1}{2}}(|\mathrm{div} v|_{2}^{\frac{1}{2}}+ |\Delta v|_{2}^{\frac{1}{2}} )|\mu|_{2}^{\frac{1}{2}}(|\nabla \mu|_{2}^{\frac{1}{2}} + |\Delta \mu|_{2}^{\frac{1}{2}})|\nu|_{2}\Big{)}.
\end{align*}
\end{lemma}
 \begin{lemma}
 	\label{lemma:2.3}
Let $v\in H^{1}(T\D|T\mathbf S^2), \mu\in H^{1}(T\D|T\mathbf S^2) \Big{(}\text{ or }\mu\in H^{1}(\D)\Big{)}$ and   $\nu\in H^{1}(T\D|T\mathbf S^2)     \Big{(}\text{ or }  \nu\in H^{1}(\D)\Big{)}.$ Then, there exists a positive constant $c$ independent of $v,\mu$ and $\nu$ such that
\begin{align*}
&\Big{|}\langle \Big{(} \int_{\xi}^{1} \mathrm{div} v(t; \theta, \phi, \xi' )d\xi'  \Big{)} \mu,  \nu \rangle\Big{|}\\
\leq& c| \mathrm{div} v|_{2}|\mu|_{2}^{\frac{1}{2}}( |\mu|_{2}^{\frac{1}{2}}+ |\nabla_{e_{\theta}}\mu|_{2}^{\frac{1}{2}}+ |\nabla_{e_{\varphi}}\mu|_{2}^{\frac{1}{2}} )
|\nu|_{2}^{\frac{1}{2}}( |\nu|_{2}^{\frac{1}{2}}+ |\nabla_{e_{\theta}}\nu|_{2}^{\frac{1}{2}}+ |\nabla_{e_{\varphi}}\nu|_{2}^{\frac{1}{2}} ),\\
\Big{(} \text{ or }\leq& c| \mathrm{div} v|_{2}|\mu|_{2}^{\frac{1}{2}}( |\mu|_{2}^{\frac{1}{2}}+ |\nabla\mu|_{2}^{\frac{1}{2}})
|\nu|_{2}^{\frac{1}{2}}( |\nu|_{2}^{\frac{1}{2}}+ |\nabla\nu|_{2}^{\frac{1}{2}})   \Big{)}.
\end{align*}
\end{lemma}
Now define linear operators $A_{i}: \mathcal V_{i}\mapsto \mathcal V_{i}', i=1,2,3$ :
\begin{eqnarray*}
	&\langle A_{1}\v, \wt\v    \rangle= \langle \v,  \wt\v   \rangle_{\mathcal V_{1}},&\text{ for any }\v,\wt\v\in\mathcal V_1,\\
	&\langle A_{2}T ,  \wt T   \rangle= \langle T,  \wt T   \rangle_{\mathcal V_{2}},&\text{ for any } T, \wt T\in\mathcal V_2,\\
	&\langle A_{3} q ,  \wt q   \rangle= \langle q,  \wt q   \rangle_{\mathcal V_{3}},&\text{ for any }q,\wt q\in\mathcal V_3.
\end{eqnarray*}
Denote $D(A_{i})=\{\eta\in \mathcal V_{i},    A_{i}\eta\in \mathcal H_{i} \}.$ 	From the second part in Lemma \ref{lemma:ibp}, we see that $A_i$'s are  positive definite, self-adjoint operators, according to the classic spectral theory we can define the power $A_{i}^{s}$ for any $s\in \mathbb{R}.$ Then we have $D(A_{i}^{\frac{1}{2}})=\mathcal V_{i}$ and $D(A_{i}^{-\frac{1}{2}})=\mathcal V_{i}'.$ Moreover,
\begin{eqnarray*}
	D(A_{i})\subset \mathcal V_{i}\subset \mathcal H_{i}\subset \mathcal V_{i}'\subset D(A_{i})',
\end{eqnarray*}
where $D(A_{i})' $ is the dual space of $ D(A_{i})$,  and the embeddings above are all compact.
Now we denote by $P_H$ the Leray type projection operator from $L^2(T\D|T\mathbf S^2)\times(L^2(\D))^2$ onto $\mathcal H$, and take $\mathcal V^{(2)}$ as the closure of $\mathcal V$ in the $H^2(T\D|T\mathbf S^2)\times(H^2(\D))^2$ norm,  then we define the principal linear portion of the system:
$$AU=P_H\begin{pmatrix}
-\Delta\v-\partial_{\xi\xi}\v\\
-\Delta T-\partial_{\xi\xi}T\\
-\Delta q-\partial_{\xi\xi}q
\end{pmatrix},\text{ for any }U\in D(A),$$
where
\begin{align*}
	D(A):=\{U=(\v,T,q)\in \mathcal V^{(2)};\  &\text{On } \xi=1:\ \partial_\xi\v=0, w=0, \partial_\xi T=-\alpha T, \partial_\xi q=-\beta q;\\
	&\text{On } \xi=0:\ \partial_\xi\v=0, w=0, \partial_\xi T=0, \partial_\xi q=0\},
\end{align*}
and $\qv{AU,\wt U}=\qv{U,\wt U}_{\mathcal V}$ for all $U,\wt U\in D(A)$.
		Let $\bu=(u_\theta,u_\varphi) \in C^\infty(T\D|T\mathbf S^2), \wt\bu= (-\Delta-\partial_{\xi\xi} )\bu.$ Applying Lemma $2.1$ and boundary conditions , we have
		\begin{align*}
			&\langle (-\Delta-\partial_{\xi\xi} )\bu ,  (-\Delta-\partial_{\xi\xi}) \bu  \rangle \\
			&=\langle (-\Delta-\partial_{\xi\xi} )\bu ,  \wt\bu \rangle\\
			&=\langle -\Delta \bu , \wt\bu \rangle+ \langle -\partial_{\xi\xi} \bu ,  \wt\bu \rangle\\
			 &=\int_{\D}(\nabla_{e_\theta}\bu\cdot\nabla_{e_\theta}\wt\bu+\nabla_{e_\varphi}\bu\cdot\nabla_{e_\varphi}\wt\bu+\bu\cdot\wt\bu)d\D+\int_{\D}\partial_{\xi}\bu
			\cdot \partial_{\xi}\wt\bu d\D.
		\end{align*}

		\par
		Next, we define the diagnostic function:
		\begin{equation}
			\label{eqn:diagfunc}
			w(\v):=\int_{\xi}^1\text{div }\v(t;\theta,\varphi,\xi')d\xi',\hspace{1cm}\v\in \mathcal V_{1}.
		\end{equation}
		Now take $U,\wt U\in D(A)$ and define the nonlinear operator as
		\begin{equation}
			\label{eqn:nonlinear}
			B(U,\wt U):=P_H\begin{pmatrix}
				\nabla_\v\wt\v\\\nabla_\v\wt T\\\nabla_\v\wt q
			\end{pmatrix}+P_H\begin{pmatrix}
				w(\v)\partial_\xi\wt\v\\w(\v)\partial_\xi \wt T\\w(\v)\partial_\xi\wt q
			\end{pmatrix}.
		\end{equation}
		Also, we define the pressure operator, Coriolis operator and external operator as
		\begin{equation}
			\label{eqn:pressureop}
			A_pU:=P_H\begin{pmatrix}
				\int_{\xi}^1\frac{br_s}{r}\nabla[(1+aq)T]d\xi'\\
				-\frac{br_s}{r}(1+aq)\bigg(\int_{\xi}^1\text{div }\v(t;\theta,\varphi,\xi')d\xi'\bigg)\\0
			\end{pmatrix},\ U\in \mathcal V,
		\end{equation}
		\begin{equation}
			\label{eqn:corop}
			EU:=P_H\begin{pmatrix}
				\frac{f}{R_0}\vec k\times\v\\0\\0
			\end{pmatrix},\ U\in \mathcal H,
		\end{equation}
		\begin{equation}
			F:=P_H\begin{pmatrix}
				0\\Q_T\\Q_q
			\end{pmatrix}.
		\end{equation}
\par
		Finally, Let $(B_{i}(t))_{i\in \mathbb{N}^{+}}$ be a sequence of one-dimensional, independent, identically distributed, two-sided Brownian motions, defined on the complete probability space $(\Omega, \mathcal{F}, \mathbb{P} )$. For $j=1,2,3$, we write $(e_{i,j})_{i\in \mathbb{N}^{+}},$ for an orthonormal basis of $H_{j},$ consisting of eigenfunctions of the operator $A_{j}$ and $(\gamma_{i,j})_{i\in \mathbb{N}^{+},}$ for the sequence of the corresponding eigenvalues.
We introduce the $H_{j}$-valued Wiener process $(W_{j}(.,t))_{t\in \mathbb{R}^{+}}$ with $j=1,2,3$ by setting
\begin{eqnarray}
W_{j}(.,t):=\sum^{\infty}_{i=1}\lambda_{i,j}^{\frac{1}{2}}e_{i,j}(.)B_{i}(t),
\end{eqnarray}
where $(\lambda_{i,j})_{i\in \mathbb{N}^{+}} $ is a sequence of positive numbers such that the series converge a.s. in the strong topology of $H_{j}$ .

	\par
		With all the above operator notations, we could reformulate (\ref{eqn:moist2}) into the following abstract evolution system,
		\begin{equation}
			\label{eqn:abstractmoist}
			dU+(AU+B(U)+A_pU+EU)dt=Fdt+dW,\ U(t_0)=U_0,
		\end{equation}
		where $U_0=(\v_0,T_0,q_0)(\theta,\varphi,\xi)$.
\begin{definition}
For $\mathbb{P}$-a.e. $\omega\in \Omega,$  we say a continuous $\mathcal{V}$-valued $(\mathcal{F}_{t_{0},t})=\sigma(W_{i}(s)-W_{i}(t_{0}), s\in [t_{0}, t] , i=1,2,3 )$adapted random field $(U(.,t))_{t\in [t_{0}, \tau]}:=(\mathbf{v}(.,t),T(.,t), q(.,t))_{t\in [t_{0}, \tau]}$ defined on $(\Omega, \mathcal{F}, \mathbb{P})$ is a \textbf{strong solution} to problem $(1.7a)$-$(1.7f)$ with $\v_{0}\in \mathcal{V}_{1}, T_{0}\in \mathcal{V}_{2}, q_{0} \in \mathcal{V}_{3}$ and $t_{0},\tau\in \mathbb{R}, \tau\geq t_{0},$  if $(U(.,t))_{t\in [t_{0}, \tau]}$ satisfies $(1.7a)$-$(1.7f)$ in the weak  sense such that
\begin{eqnarray*}
&&\mathbf{v}\in C([t_0,\mathcal{\tau}];\mathcal{V}_{1})\cap L^{2}([t_0,\mathcal{\tau}]; (H^{2}(\D))^{2}),\\
&&T\in C([t_0,\mathcal{\tau}];\mathcal{V}_{2})\cap L^{2}([t_0,\mathcal{\tau}]; H^{2}(\D)),\\
&&q\in C([t_0,\mathcal{\tau}];\mathcal{V}_{3})\cap L^{2}([t_0,\mathcal{\tau}]; H^{2}(\D)).
\end{eqnarray*}
\end{definition}
Similarly, we can define the strong solution to $(5.3a)-(5.3h).$
	In this part, we state our results about the local well-posedness, the global well-posedness, the existence of random attractor and invariant measure.
\begin{theorem}[Existence of local solutions]
	\label{thm:2.1}
	 If $Q_T, \partial_{\xi}Q_{T}$, $Q_q$, $\partial_{\xi}Q_{q}\in L^{2}(\D),$ $(\mathbf{v}_{0}, T_{0}, q_{0})\in \mathcal{V}$ then, for $\mathbb P$-a.e. $\omega \in\Omega,$ there exists a stopping time $T^{*}>0$
	such that $( \mathbf{v},T,q)$ is a strong solution of the system $(1.7a)-(1.7f)$ on the interval $[0, T^{*}].$
\end{theorem}
To consider the local well-posedness,  we separate the equation $(1.7a)-(1.7f)$ into a deterministic linear equation corresponding to $(1.7a)$ and the stochastic nonlinear part with zero initial condition, i.e., $(1.7a)-(1.7f)$ with $\bu(t_0)=\bu_{0}$ replaced by $\bu(t_0)=0.$ The global well-posedness of linear part is well known. The proof of the local well-posedness of the nonlinear part with zero initial condition is also classic. We first obtain that the sequence of solution to the approximation of $(1.7a)-(1.7f)$ is bounded in $L^{2}([t_{0}, \tau]; (H^{2}(\D))^{4})$, then we use Aubin-Lions Lemma to obtain an strongly convergent subsequence of solution in $L^{2}([t_{0}, \tau]; (H^{1}(\D))^{4})$. Reasoning on weakly and strongly convergent subsequences one gets the existence of a solution with the regularity specified by Theorem \ref{thm:2.1}. The proof is similar to \cite{GH11}, so we omit it here.
\begin{theorem}[Existence of global solutions]
	\label{thm:2.2}
	Let $Q_T$, $\partial_{\xi}Q_{T}$, $Q_q$, $\partial_{\xi}Q_{q}\in L^{2}(\D),$ $(\mathbf{v}_{0}, T_{0}, q_{0})\in \mathcal{V}$, and $\sum_{i=1}^{\infty}\lambda_{i,j}^{2}\gamma_{i,j}^{2+\sigma}<\infty,j=1,2,3,$ for small positive constant $\sigma$.  Then, for arbitrary $\tau>t_0$, there exists a unique strong
	solution $( \mathbf{v},T,q)$ to the system $(1.7a)-(1.7f)$ on the interval $[t_0, \tau],$ which is Lipschitz continuous with respect to the initial data in $\mathcal{V}.$
\end{theorem}
Now we give preliminary knowledge about random attractors. 	Let $(X, d)$ be a polish space and $(\wt{\Omega}, \wt{\mathcal{F}}, \wt{\mathbb{P}}  )$ be a probability space, where $ \wt{\Omega}$ is the two -sided Wiener space $C_{0}(\mathbb{R}; X )$ of continuous functions with values in $X$, equal to $0$ at $t=0$. We consider a family of mappings
$\S(t,s;\omega):X\rightarrow X,\ -\infty<s \leq t< \infty,$ parametrized by $\omega\in \wt{\Omega},$ satisfying for $ \wt{\mathbb{P}}$-$a.s.\ \omega$,  the following properties (i)-(iv):
\begin{enumerate}[(i)]
	\item $\S(t,r;\omega)\S(r,s;\omega)x= \S(t,s;\omega)x$ for all $s\leq r\leq t$ and $x\in X$,
	\item $\S(t,s;\omega)$ is continuous in $X,$ for all $s\leq t$,
	\item for all $s<t$ and $x\in X$, the mapping
	$$\omega\mapsto \S(t,s;\omega)x$$
	is measurable from $(\wt{\Omega},\wt{\mathcal{F}})$ to $(X, \mathcal{B}(X )  )$ where $\mathcal{B}(X ) $ is the Borel-$\sigma$- algebra of $ X$,
	\item  for all $t, x\in X,$ the mapping $s\mapsto \S(t,s;\omega)$ is right continuous at any point.
\end{enumerate}
We define for $A,B\in 2^X$ with $A,B\neq\emptyset$, $d(A,B)=\sup\{\inf\{d(x,y):y\in B \}:x\in A \}$, and it follows that $d(x,B)=d(\{x\},B)$. We now give the following definitions.
\begin{definition}
	A set-valued map $K: \wt{\Omega}\rightarrow 2^{X}$ taking values in the closed subsets of $X$ is said to be \textbf{measurable}, if for each $x\in X$, the map
	$\omega\mapsto d(x, K(\omega))$ is measurable. A closed set-valued measurable map $K:\wt{\Omega}\rightarrow 2^{X}$ is called a \textbf{random closed set}.  	 
\end{definition}
\begin{definition}
	Given $t\in \mathbb{R}$ and $\omega\in \wt{\Omega}, K(t,\omega)\subset X$ is called an \textbf{attracting set} at time $t$ if , for all bounded sets $B\subset X,$
	$$
	d(\S(t,s;\omega)B, K(t,\omega) )\rightarrow 0,\ \ \text{provided }\ s\rightarrow -\infty.
	$$
	Moreover, if for all bounded sets $B\subset X,$ there exists $t_{B}(\omega)$ such that for all $s\leq t_{B}(\omega)$,
	$$
	\S(t,s;\omega)B\subset K(t,\omega),
	$$
	we say $K(t,\omega) $ is an \textbf{absorbing set} at time $t.$
\end{definition}
Let $\{\vartheta_{t}:\wt{\Omega}\rightarrow \wt{\Omega}   \}_{t\in\R}$ be a family of measure preserving transformations of the probability space $(\wt{\Omega}, \wt{\mathcal{F}},\wt{ \mathbb{P}} )$ such that for all $s< t$ and $\omega\in \wt{\Omega}$,
\begin{enumerate}[(a)]
	\item $(t,\omega)\rightarrow \vartheta_{t}\omega$ is measurable,
	\item $\vartheta_{t}(\omega)(s)=\omega(t+s)-\omega(t)$,
	\item $\S(t,s;\omega)x=\S(t-s,0;\vartheta_{s}\omega)x.$
\end{enumerate}
We defined $(\vartheta_{t} )_{t\in T}$ as a flow, and
$((\wt{\Omega}, \wt{\mathcal{F}},\wt{ \mathbb{P}} ), (\vartheta_{t} )_{t\in \R} )$ is a measurable dynamical system.
\begin{definition}
	Given a bounded set $B\subset X$, the set
	\begin{eqnarray*}
		\mathcal{A}(B,t,\omega)=\bigcap\limits_{T\leq t}\overline{\bigcup\limits_{s\leq T}\S(t, s,\omega)B}
	\end{eqnarray*}
	is said to be the \textbf{$\Omega$-limit set of $B$} at time $t$. Obviously, if we denote $\mathcal{A}(B,0,\omega)=\mathcal{A}(B,\omega),$ we have
	$\mathcal{A}(B,t,\omega)=\mathcal{A}(B,\vartheta_{t}\omega).$
\end{definition}
We may identify
\begin{equation}
\mathcal A(B,t,\omega)=\{x\in X:\text{ there exist }s_n\to-\infty,\ x_n\in B\text{ such that }\lim_{n\to\infty}\S(t,s_n,\omega)x_n=x\}.
\end{equation}
Furthermore, if there exists a compact attracting set $K(t,\omega)$ at time $t,$ it is not difficult to check that $\mathcal{A}(B,t,\omega)$ is a nonempty compact subset of $X$ and $\mathcal{A}(B,t,\omega)\subset K(t,\omega)$. Now we are ready to give the definition of random attractors as follows:
\begin{definition}
	If, for any $t\in \mathbb{R}$ and $\omega\in \wt{\Omega},$ the random closed set $\omega\rightarrow \mathcal{A}(t,\omega)$ satisfying the following properties:
	\begin{enumerate}[(1)]
		\item  $\mathcal{A}(t,\omega)$ is a nonempty compact subset of $X,$
		\item  $\mathcal{A}(t,\omega)$ is the minimal closed attracting set,
		i.e., if $\wt{\mathcal{A}}(t,\omega)$ is another closed attracting set, then $\mathcal{A}(t,\omega)\subset \wt{\mathcal{A}}(t,\omega),$
		\item it is invariant,  in the sense that, for all $s\leq t,$
		$$
		\S (t,s;\omega)\mathcal{A}(s,\omega)=\mathcal{A}(t,\omega).
		$$	
	\end{enumerate}
	$\mathcal{A}(t,\omega)$  is called the \textbf{random attractor}.
\end{definition}
We finish this section with our main result, the existence of random attractor and invariant measure for $(1.7a)$-$(1.7f)$.
\begin{theorem}[Existence of random attractor ]
	\label{thm:2.3}
	In addition to the conditions in Theorem \ref{thm:2.2}, we assume $|\frac{br_{s}}{r}|\leq \min\{\frac{1}{2},\alpha,\beta\}$. Then the solution operator $(\mathcal S(t,s;\omega))_{t\geq s,\omega\in \tilde{\Omega}} $  of 3D stochastic PEs $(1.7a)-(1.7f): \mathcal S(t,s;\omega)(\v_{s}, T_{s},q_s)=(\v(t), T(t),q(t) ) $ has properties $\mathrm{(i)}-\mathrm{(iv)}$ and possesses a compact absorbing ball $\mathcal{B}(0,\omega)$ in $\mathcal{V}$ at time $0.$  Furthermore, for $\tilde{\mathbb{P}}$-a.e. $\omega,$ the set
	$$\mathcal{A}(\omega)=\overline{\bigcup_{B\subset \mathcal{V}}\mathcal{A}(B,\omega) } $$
 is the random attractor of stochastic PEs, where the union is taken over all the bounded subsets of $\mathcal V$.
\end{theorem}

With the above conclusions, we can prove the existence of invariant measures for the system (\ref{eqn:moist2}).
\par
Let $U_{0}:=(\v_{0},T_{0},q_{0})\in  \mathcal{V}$. In the following, we denote by $$U(t,\omega; U_{0}):=(\v(t,\omega;t_0,\v_{0}), T(t,\omega;t_0,T_{0}), q(t,\omega;t_0,q_{0})) $$ the solution to $(\ref{eqn:moist2})$ with $(\v(t_0)=\v_{0}, T(t_0)=T_{0}, q(t_{0})=q_{0}).$ Following the standard argument, we can show that $U(t,\omega; U_{0}), t\in [t_0, \mathcal{T}],t_{0}\leq \mathcal{T}$ is Markov in the following sense:\\
for every bounded, $\mathcal{B}(\mathcal{V})$-measurable $F:\mathcal{V}\rightarrow \mathbb{R},$ and all $s,t\in [t_{0}, \mathcal{T}]$, $t_{0}\leq s\leq t\leq  \mathcal{T}$
\begin{eqnarray*}
\mathbb{E}(F(U(t,\omega; U_{0}))|\mathcal{F}_{s} )(\omega)=\mathbb{E}(F(U(t,s, U(s))) )\ \ \mathrm{for}\ \tilde{\mathbb{P}}-a.e.\ \omega\in \Omega.
\end{eqnarray*}
where $\mathcal{F}_{s}=\mathcal{F}_{t_{0},s}$ (see Definition 2.1) and $ U(t,s, U(s))$ is the solution to (\ref{eqn:moist2}) at time $t$ with initial data $U(s).$
\par
For $B\in \mathcal{B}(\mathcal{V})$ the collection of Borel measurable subset on $\mathcal{V},$ we define
\begin{eqnarray*}
\tilde{\mathbb{P}}_{t}(U_{0}, B)=\tilde{\mathbb{P}}((U(t,\omega; U_{0})\in B ).
\end{eqnarray*}
\par
For any probability  measure $\nu$ defined on $\mathcal{B}(\mathcal{V}),$ we denote by $(\nu \tilde{\mathbb{P}}_{t})(\cdot)= \int_{\mathcal{V}}\tilde{\mathbb{P}}_{t}(x,\cdot)\nu(dx)$ the distribution at time $t$ of the solution to  (\ref{eqn:moist2}) with initial condition having the distribution $\nu.$
\par
For $t\geq t_0$ and any function $f\in C_{b}(\mathcal{V};\mathbb{R})$ the set of continuous and bounded functions from $\mathcal{V}$ into $\mathbb{R},$ denote
\begin{eqnarray*}
\tilde{\mathbb{P}}_{t}f(U_{0})=\mathbb{E}[f(U(t,\omega; U_{0})]=\int_{V}f(x)\tilde{\mathbb{P}}_{t}(U_{0}, dx).
\end{eqnarray*}
\begin{definition}
Let $\rho$ be a probability measure on $\mathcal{B}(\mathcal{V})$. We say that $\rho$ is $\mathbf{an}$ $\mathbf{invariant}$ $\mathbf{measure}$ for $\tilde{\mathbb{P}}_{t}$ if we have
\begin{eqnarray*}
\int_{\mathcal{V}}f(x)\rho(dx)=\int_{\mathcal{V}}\tilde{\mathbb{P}}_{t}f(x)\rho(dx)
\end{eqnarray*}
for all $f\in C_{b}(\mathcal{V};\mathbb{R})$ and $t\geq 0.$
\end{definition}
\par
Let $\mu_{\cdot}$ be a transition probability from $\tilde{\Omega}$ to $\mathcal{V}$, i.e., $\mu_{\cdot}$ is a Borel probability measure on $\mathcal{V}$ and $\omega\rightarrow \mu_{\cdot}(B)$ is measurable for every Borel set $B\subset \mathcal{V}.$ Denote by $\mathcal{P}_{\tilde{\Omega}}(\mathcal{V}) $ the set of transition probabilities with $\mu_{\cdot}$ and $\nu_{\cdot}$ identified if $\tilde{\mathbb{P}}\{ \omega:\mu_{\omega}\neq \nu_{\omega}\}=0. $
\par
In view of Proposition $4.5$ in $\cite{CF}$, the existence of random attractor obtained in Theorem $2.3$ implies the existence of  an invariant Markov measure $\mu_{\cdot}\in \mathcal{P}_{\tilde{\Omega}}(V)$ for $\mathcal{S}$ such that $ \mu_{\omega}(\mathcal{A}(\omega))=1\ \tilde{\mathbb{P}}$-a.e.. Therefore, by $\cite{C}$ there exists an invariant measure for the markov semigroup $\tilde{\mathbb{P}}_{t}$ and it is given by
$$\rho(B)=\int_{\tilde{\Omega}}\mu_{\omega}(B)\tilde{\mathbb{P}}(d\omega),$$
where $B\subseteq \mathcal{V}$ is a Borel set and $f\in C_{b}(\mathcal{V};\mathbb{R})$. If the invariant measure $\rho$ for $\tilde {\mathbb{P}}$ is unique, the invariant Markov measure $\mu_{\cdot}$ for $\mathcal{S}$ is unique and given by
$$\mu_{\omega}=\lim\limits_{t\rightarrow \infty}\mathcal{S}(0,-t, \omega)\rho. $$
Summarizing the above argument, we arrive at
\begin{theorem}
The Markov semigroup $ (\tilde{\mathbb{P}}_{t})_{t\geq 0}$ induced by the solution $(U(t,\omega; U_{0}))_{t\geq 0}$  to (\ref{eqn:moist2}) has an invariant measure $\rho$ with $ \rho(\mathcal{A}(\omega))=1\ \tilde{\mathbb{P}}$-a.e..
\end{theorem}

\section{Global well-posedness of strong solutions}

We need the regularity of the strong solution and\textit{ a priori} estimates to prove the compact property of solution operator, which is the key to prove Theorem \ref{thm:2.3}.
\par
The following Lemma, a special case of a general result of Lions and Magenes $\cite{LM}$,  will help us to show the continuity of the solution to stochastic PEs with respect to time in $(H^{1}(\D))^{4}.$ For the proof of the Lemma we can refer to $\cite{T}$ for details.
\begin{lemma}
Let $V , H, V'$ be three Hilbert spaces such that $V \subset H = H¡ä\subset V' $
, where $ H'$ and $V'$ are the dual spaces of $H$ and $V$ respectively. Suppose $u \in
L^{2}(0, T; V )$ and $u'\in L^{2}(0, T; V')$. Then $u$ is almost everywhere equal to a function
continuous from $[0, T]$ into $H$.
\end{lemma}

\par
Before giving our proof, we should notify that the global well-posedness of $(1.7a)-(1.7f)$ with initial condition $(1.6a)-(1.6c)$ is equivalent to the system  $(5.3a)-(5.3h).$
\par
\textbf{Proof of Theorem \ref{thm:2.2}.}
In the following, we will complete our proof of the global well-posedness of stochastic PEs by three steps. Firstly, we will prove the global existence of strong solution. Then, we will show that the solution is continuous in the space $\mathcal{V}$ with respect to $t$. Finally,  we will obtain the continuity in $\mathcal{V}$ with respect to the initial data.  \\
$\mathbf{Step}\ 1$: We prove the global existence of the strong solutions.\\
 We denote by $[t_0, \tau_{*})$ the maximal interval of existence of the solution of $(5.3a)-(5.3h),$ we infer that $\tau_{*}=\infty,$ a.s..
Otherwise, if there exists $A\in \mathcal{F}$ such that $\mathbb{P}(A)>0$ and for fixed $\omega\in A, \tau_{*}(\omega)<\infty,$ it is clear that
\begin{eqnarray*}
\limsup\limits_{t\rightarrow \tau_{*}^{-}(\omega)}(\| \bu(t)\|_{1}+\|S(t)\|_{1} +\|p\|_{1})=\infty,\ \ \mathrm{for}\ \mathrm{any}\ \omega\in A,
\end{eqnarray*}
which contradicts the estimates  \eqref{eqn:5.56}, \eqref{eqn:5.78}, \eqref{eqn:5.84}
and \eqref{eqn:5.93} given in the Appendix. Therefore $\tau_{*}=\infty,$ a.s.,  and the strong solution $(\bu, S,p)$ exists globally in time a.s..\\
$\mathbf{Step}\ 2$: We show the continuity of strong solutions with respect to $t$.\\
Taking inner product between $\partial_tA_1^{\frac{1}{2}}\bu$ and $\eta$, by $(5.3a)$, one can get
\begin{align*}
&\langle \partial_{t}A_{1}^{\frac{1}{2}} \bu, \eta \rangle= \langle \partial_{t}\bu, A_{1}^{\frac{1}{2}}  \eta \rangle=- \langle A_{1}\bu, A_{1}^{\frac{1}{2}}\eta\rangle - \langle \nabla_{\bu+Z_{1}} (\bu+Z_{1}),A_{1}^{\frac{1}{2}}\eta\rangle\nonumber\\
& -\langle w(\bu+Z_{1}) \partial_{\xi}(\bu+Z_{1}), A_{1}^{\frac{1}{2}}  \eta \rangle-\frac{f}{R_{0}}\langle  ( \bu+Z_{1})^{\bot}, A_{1}^{\frac{1}{2}}  \eta \rangle  \nonumber\\
& - \langle\int_{\xi}^{1}\frac{br_{s}}{r}\nabla[(1+a(Z_{3}+p))(Z_{2}+S) ] d\xi',   A_{1}^{\frac{1}{2}}  \eta \rangle+ \gamma\langle Z_{1}, A_{1}^{\frac{1}{2}}  \eta \rangle ,
\end{align*}
where we have used $ \langle \nabla \Phi_{s},   A_{1}^{\frac{1}{2}}  \eta   \rangle=0$ which follows by integration by parts formula.
Taking a similar argument in  \eqref{eqn:5.87}, we get
$$\qv{w(\bu+Z_1)\partial_{\xi}(\bu+Z_1),A_1^{\frac{1}{2}}\eta}\leq C\norm{\bu+Z_1}_1\norm{\bu+Z_1}_2|A_1^{\frac{1}{2}}\eta|_2.$$
By the H\"older inequality and the Sobolev embedding theorem, we have
\begin{align*}
&- \langle\int_{\xi}^{1}\frac{br_{s}}{r}\nabla[(1+a(Z_{3}+p))(Z_{2}+S) ] d\xi',   A_{1}^{\frac{1}{2}}  \eta \rangle\notag\\
&\leq C|A_{1}^{\frac{1}{2}}\eta|_{2}(\|Z_{3}+p\|_{2}\|Z_{2}+S\|_{1}+ \|Z_{3}+p\|_{1}\|Z_{2}+S\|_{2} ).
\end{align*}
Similarly,
$$
- \langle \nabla_{\bu+Z_{1}} (\bu+Z_{1}),A_{1}^{\frac{1}{2}}\eta\rangle\leq C \|\bu+Z_{1}\|_{1}\|\bu+Z_{1}\|_{2}|A_{1}^{\frac{1}{2}} \eta|_{2}.
$$
Therefore, combining the above estimates yields
\begin{align*}
\|\partial_{t}(A_{1}^{\frac{1}{2}}\bu)\|_{\mathcal{V}_{1}'}\leq C(&\|\bu\|_{2}+\|\bu+Z_{1}\|_{1}\|\bu+Z_{1}\|_{2}+|\bu|_{2}+|Z_{1}|_{2}\notag\\
&+\|Z_{3}+p\|_{2}\|Z_{2}+S\|_{1}+ \|Z_{3}+p\|_{1}\|Z_{2}+S\|_{2}).
\end{align*}
By \textit{a priori} estimates of $\bu$ in Appendix,
\begin{eqnarray*}
\bu\in L^{\infty}([t_0, \tau]; \mathcal{V}_{1})\cap L^{2}([t_0, \tau]; H^{2}(T\D|T\mathbf{S}^{2})),  Z_{1}\in C([t_0,\tau]; H^{3}(T\D|T\mathbf{S}^{2})),
\end{eqnarray*}
and $Z_{2}, Z_{3}\in C([t_0,\tau]; H^{3}(\D))$ for all $\tau>t_0,$
we obtain
\begin{eqnarray*}
A_{1}^{\frac{1}{2}}\bu\in L^{2}([t_0,\tau]; \mathcal{V}_{1} ),\ \ \ \partial_{t}(A_{1}^{\frac{1}{2}}\bu )\in L^{2}([t_0,\tau]; \mathcal{V}_{1}' ),
\end{eqnarray*}
which together with Lemma 3.1 implies $$A_{1}^{\frac{1}{2}}\bu\in C([t_0,\tau]; H_{1})\ \mathrm{or}\  \bu\in C([t_0,\tau]; \mathcal{V}_{1})\text{ a.s. }. $$
Similarly, we can prove
\begin{eqnarray*}
S\in C([t_0,\tau]; \mathcal{V}_{2})\ \mathrm{and}\   p\in C([t_0,\tau]; \mathcal{V}_{3}).
\end{eqnarray*}

$\mathbf{Step}\ 3$: We obtain the continuity of strong solutions in $\mathcal{V}$ with respect to the initial data.\\ Let $(\v_1, T_1, q_1 )$ and $(\v_2, T_{2}, q_{2}) $ be two solutions of the system
$(1.7a)-(1.7f)$ with corresponding pressure $\Phi_{s}{'}$ and $\Phi_{s}{''},$ and initial data $(\v_{t_0}^{1}, T_{t_0}^{1}, q_{t_0}^{1})$ and
$(\v_{t_0}^{2}, T_{t_0}^{2}, q_{t_0}^{2})$ respectively. Denote by $\v=\v_{1}-\v_{2},T=T_{1}-T_{2},  q=q_{1}-q_{2} $  and $\Phi_{s}=\Phi_{s}{'}-\Phi_{s}{''}.$
Then we derive from $(1.7a)-(1.7f)$ that
\begin{align}
\label{eqn:3.1}
&\partial_{t}\v+L_{1}\v+\nabla_{\v_{1}}\v+\nabla_{\v}\v_{2}+\Big{(}\int_{\xi}^{1}\mathrm{div} \v_{1}(x,y,\xi',t )d\xi'   \Big{)}\partial_{\xi}\v\nonumber\\
&+\Big{(}\int_{\xi}^{1}\mathrm{div} \v(x,y,\xi',t )d\xi'   \Big{)}\partial_{\xi}\v_{2}+\frac{f}{R_{0}}\v^{\bot}+ \mathrm{grad} \Phi_{s}\nonumber\\
&+\int_{\xi}^{1}\frac{bP}{p} \mathrm{grad} T d\xi'+\int_{\xi}^{1}\frac{abP}{p} \mathrm{grad} (q_{1} T) d\xi'+\int_{\xi}^{1}\frac{abP}{p} \mathrm{grad} (q T_{2}) d\xi'=0,\\
\label{eqn:3.2}
&\partial_{t}T+L_{2}T+\nabla_{\v_{1}}T+\nabla_{\v}T_{2}+\Big{(}\int_{\xi}^{1}\mathrm{div} \v_{1}(x,y,\xi',t )d\xi'   \Big{)}\partial_{\xi}T\nonumber\\
&+\Big{(}\int_{\xi}^{1}\mathrm{div} \v(x,y,\xi',t )d\xi'   \Big{)}\partial_{\xi}T_{2}-\frac{bP}{p}\Big{(}\int_{\xi}^{1}\mathrm{div} \v(x,y,\xi',t )d\xi'   \Big{)}\nonumber\\
&-\frac{abP}{p}q_{1}\Big{(}\int_{\xi}^{1}\mathrm{div} \v(x,y,\xi',t )d\xi'   \Big{)}
-\frac{abP}{p}q\Big{(}\int_{\xi}^{1}\mathrm{div} \v_{2}(x,y,\xi',t )d\xi'   \Big{)}=0,\\
\label{eqn:3.3}
&\partial_{t}q+L_{3}q+\nabla_{\v_{1}}q+\nabla_{\v}q_{2}+\Big{(}\int_{\xi}^{1}\mathrm{div} \v_{1}(x,y,\xi',t )d\xi'   \Big{)}\partial_{\xi}q \nonumber\\ &+\Big{(}\int_{\xi}^{1}\mathrm{div} \v(x,y,\xi',t )d\xi'   \Big{)}\partial_{\xi}q_{2}=0,\\
&\v|_{t_0}= \v_{t_0}^{1}-\v_{t_0}^{2},\ \  T|_{t_0}=T_{t_0}^{1}-T_{t_0}^{2},\ \ q|_{t_0}=q_{t_0}^{1}-q_{t_0}^{2},\\
&\xi=1:\ \ \partial_{\xi}\v=0,\ \ \partial_{\xi}T=-\alpha T,\ \ \partial_{\xi} q=-\beta q,\\
&\xi=0:\ \ \partial_{\xi}u=0,\ \ \partial_{\xi}T=0,\ \ \partial_{\xi} q=0.
\end{align}
Taking inner product of \eqref{eqn:3.1} with $A_{1}\v$ in $L^{2}(T\D|T\mathbf{S}^{2})$ we obtain
\begin{align}
\label{eqn:3.7}
\frac{1}{2}\frac{d}{dt}\|\v\|_{1}^{2}+|A_{1}\v|_{2}^{2}=& -\langle \nabla_{\v_{1}}\v,     A_{1}\v \rangle-\langle \nabla_{\v}\v_{2},    A_{1}\v \rangle-\langle\Big{(}\int_{\xi}^{1}\mathrm{div} \v(x,y,\xi',t )d\xi'   \Big{)}\partial_{\xi}\v_{2},    A_{1} \v  \rangle\nonumber\\
&-\langle \Big{(}\int_{\xi}^{1}\mathrm{div} \v_{1}(x,y,\xi',t )d\xi'   \Big{)}\partial_{\xi} \v,   A_{1}\v \rangle-\langle \int_{\xi}^{1}\frac{bP}{p} \mathrm{grad }\mathbb{T} d\xi',    A_{1} \v  \rangle\notag\\
&-\langle (\frac{f}{R_{0}} \v^{\bot}+\mathrm{grad} \Phi_{s}),   A_{1}\v   \rangle-\langle \int_{\xi}^{1}\frac{abP}{p} \mathrm{grad}(q_{1}T)d\xi',     A_{1} \v    \rangle\notag\\
&-\langle \int_{\xi}^{1}\frac{abP}{p} \mathrm{grad}(qT_{2})d\xi'   ,     A_{1} \v\rangle\nonumber\\
=:&\sum_{i=1}^{8}k_{i}.
\end{align}
By the H$\mathrm{\ddot{o}}$lder inequality, the Agmon inequality and Young's inequality , we have
\begin{align*}
k_{1}\leq& |\v_{1}|_{\infty}(|\nabla_{e_{\theta}}\v|_{2}+  |\nabla_{e_{\varphi}}\v|_{2})|A_{1}\v|_{2}\\
\leq &c\|\v_{1}\|_{1}^{\frac{1}{2}}|A_{1}\v_{1}|_{2}^{\frac{1}{2}}\|\v\|_{1}|A_{1}\v|_{2}\\
\leq & \varepsilon |A_{1}\v|_{2}^{2}+c\|\v_{1}\|_{1}|A_{1}\v_{1}|_{2}\|\v\|_{1}^{2}.
\end{align*}
Similarly, we obtain
$$
k_{2}\leq |\v|_{\infty}\|\v_{2}\|_{1}|A_{1}\v|_{2}\leq \varepsilon |A_{1}\v|_{2}^{2}+c\|\v\|_{1}^{2}\|\v_{2}\|_{1}^{4}.
$$
Taking an analogous argument as \eqref{eqn:5.87}, we have
$$
k_{3}+k_{4}\leq \varepsilon |A_{1}\v|_{2}^{2}+c\|\v\|_{1}^{2}\|\v_{2}\|_{1}^{2}\| \v_{2}\|_{2}^{2}+c\|\v\|_{1}^{2}\|\v_{1}\|_{1}^{2}\|\v_{1}\|_{2}^{2} .
$$
In view of Lemma \ref{lemma:ibp} and the H\"older inequality, we obtain
$$
k_{5}+k_{6}\leq \varepsilon |A_{1}\v|_{2}^{2}+c\|T\|_{1}^{2}+c|\v|_{2}^{2}.
$$
To estimate $k_{7}$,  by the H\"older inequality, the Minkowski inequality and the interpolation inequality  we have
\begin{align*}
k_{7}\leq& c\int_{\mathbf S^2}(\int_{0}^{1}|\mathrm{grad} q_{1}| |T|d\xi \int_{0}^{1}|A_{1}\v|d\xi  )d\mathbf S^2\\
&+c\int_{\mathbf S^2}(\int_{0}^{1}| q_{1}| |\mathrm{grad} T|d\xi \int_{0}^{1}|A_{1}\v|d\xi  )d\mathbf S^2\\
\leq &c |A_{1}\v|_{2} \Big{(}\int_{\mathbf S^2}(\int_{0}^{1}|\mathrm{grad} q_{1}|^{2}d\xi)^{2}d\mathbf S^2\Big{)}^{\frac{1}{4}}\Big{(}\int_{\mathbf S^2}(\int_{0}^{1}|T|^{2}d\xi)^{2}d\mathbf S^2\Big{)}^{\frac{1}{4}}\\
&+c |A_{1}\v|_{2} \Big{(}\int_{\mathbf S^2}(\int_{0}^{1}| q_{1}|^{2}d\xi)^{2}d\mathbf S^2\Big{)}^{\frac{1}{4}}\Big{(}\int_{\mathbf S^2}(\int_{0}^{1}|\mathrm{grad} T|^{2}d\xi)^{2}d\mathbf S^2\Big{)}^{\frac{1}{4}}\\
\leq & \varepsilon |A_{1}\v|_{2}^{2}+c|\mathrm{grad }q_{1}|_{2}( |\mathrm{grad }q_{1}|_{2}+|\Delta q_{1}|_{2}   )(|T|_{2}^{2}+|T|_{2}|\nabla T|_{2})\\
&+c|\mathrm{grad }T|_{2}( |\mathrm{grad }T|_{2}+|\Delta T|_{2}   )|q_{1}|_{4}^{2}\\
\leq & \varepsilon |A_{1}\v|_{2}^{2}+\varepsilon |\Delta T|_{2}^{2}+c\|T\|_{1}^{2}(\|q_{1}\|_{1}^{4}+\|q_{1}\|_{2}^{2} ).
\end{align*}
Similarly, we have
\begin{align*}
k_{8}
\leq& \varepsilon |A_{1}\v|_{2}^{2}  +c|\mathrm{grad }q|_{2}( |\mathrm{grad }q|_{2}+|\Delta q|_{2}   )|T_{2}|_{4}^{2}\\
&+c|\mathrm{grad }T_{2}|_{2}( |\mathrm{grad }T_{2}|_{2}+|\Delta T_{2}|_{2}   )(|q|_{2}^{2}+|q|_{2}|\nabla q|_{2})\\
\leq & \varepsilon |A_{1}\v|_{2}^{2}+\varepsilon |\Delta q|_{2}^{2}+c\|q\|_{1}^{2}(\|T_{2}\|_{1}^{4}+\|T_{2}\|_{2}^{2} ).
\end{align*}
By \eqref{eqn:3.7} and estimates of $k_{i}$, $i=1,\dots,8$,  we get
\begin{align}
\frac{1}{2}\frac{d\|\v\|_{1}^{2}}{dt}+|A_{1}\v|_{2}^{2}\leq& \varepsilon|A_{1}\v|_{2}^{2}+\varepsilon|A_{2}T|_{2}^{2}+\varepsilon|A_{3}q|_{2}^{2}\nonumber\\
&+c\|\v\|_{1}^{2}(\|\v_{1}\|_{1}^{2}\|\v_{1}\|_{2}^{2}+ \|\v_{2}\|_{1}^{2}\|\v_{2}\|_{2}^{2}+\|\v_{2}\|_{1}^{4}+1)\nonumber\\
&+c\|T\|_{1}^{2}(\|q_{1}\|_{1}^{4}+\|q_{1}\|_{2}^{2})\nonumber\\
&+c\|q\|_{1}^{2}(\|T_{2}\|_{1}^{4}+\|T_{2}\|_{2}^{2}).
\end{align}
Taking an analogous argument as above, from \eqref{eqn:3.2} and \eqref{eqn:3.3} we have
\begin{align}
\frac{1}{2}\frac{d\| T\|_{1}^{2}}{dt}+|A_{2}T|_{2}^{2}
\leq&\varepsilon |A_{2}T|_{2}^{2}+\varepsilon |A_{1}\v|_{2}^{2}+c\|q\|_{1}^{2}\|\v_{2}\|_{1}\|\v_{2}\|_{2}+c\|T\|_{1}^{2}(1+\|\v_{1}\|_{1}^{2}\|\v_{1}\|_{2}^{2} )\nonumber\\&+c \|\v\|_{1}^{2}(1+\|q_{1}\|_{1}^{4}+\|T_{2}\|_{2}^{2}+\|T_{2}\|_{1}^{2}\|T_{2}\|_{2}^{2}  ),
\end{align}
and
\begin{equation}
\frac{1}{2}\frac{d\|q\|_{1}^{2}}{dt}+|A_{3}q|_{2}^{2}\leq\varepsilon |A_{3}q|_{2}^{2}+c\|q\|_{1}^{2}(1+ \|\v_{1}\|_{1}^{2}\|\v_{1}\|_{2}^{2})+c\|\v\|_{1}^{2}(  \|q_{2}\|_{2}^{2}+  \|q_{2}\|_{1}^{2}\|q_{2}\|_{2}^{2}).
\end{equation}
Let
\begin{align*}
&g_{1}:=1+\|\v_{1}\|_{1}^{2}\|\v_{1}\|_{2}^{2}+\|\v_{2}\|_{1}^{4}+\|\v_{2}\|_{1}^{2}\|\v_{2}\|_{2}^{2}+\|T_{2}\|_{2}^{2}\nonumber\\
&\hspace{1cm}+\|T_{2}\|_{1}^{2}\|T_{2}\|_{2}^{2}+\|q_{1}\|_{1}^{4}+\|q_{2}\|_{2}^{2}+\|q_{2}\|_{1}^{2}\|q_{2}\|_{2}^{2},\\
&g_{2}:=1+\|q_{1}\|_{2}^{2}+\|q_{1}\|_{1}^{4}+\|\v_{1}\|_{1}^{2}\|\v_{1}\|_{2}^{2},
\end{align*}
and
$$
g_{3}:=1+\|T_{2}\|_{2}^{2}+\|T_{2}\|_{1}^{4}+\|\v_{2}\|_{1}\|\v_{2}\|_{2} +\|\v_{1}\|_{1}^{2}\|\v_{1}\|_{2}^{2}.
$$
It is obvious that for arbitrary $0\leq a< b<\infty$,
$$
\int_{a}^{b}(g_{1}(t)+g_{2}(t)+g_{3}(t))dt< \infty.
$$
Therefore, we get
$$
\frac{d(\|\v\|_{1}^{2}+\|T\|_{1}^{2}+\|q\|_{1}^{2}) }{dt}\leq c( g_{1}(t)+g_{2}(t)+g_{3}(t) )(\|\v\|_{1}^{2}+\|T\|_{1}^{2}+\|q\|_{1}^{2}),
$$
 applying Gronwall inequality implies
\begin{align*}
&\|\v(t)\|_{1}^{2}+\|T(t)\|_{1}^{2}+\|q(t)\|_{1}^{2}\\
&\leq c(\|\v_{t_0}^1-\v_{t_0}^2\|_{1}^{2} + \|T_{t_0}^1-T_{t_0}^2\|_{1}^{2} +\|q_{t_0}^1-q_{t_0}^2\|_{1}^{2}  )
e^{\int_{0}^{t}(g_{1}(s)+g_{2}(s)+g_{3}(s))ds }.
\end{align*}
So far, we have shown that for $t>t_0$, the strong solution $\ (\v(t), T(t), q(t) )$ to ({\ref{eqn:moist1}}) is Lipschitz continuous in $\mathcal{V}$ with respect to the initial data $(\v_0, T_0, q_0 ) .$
\hspace{\fill}$\square$

	\section{Random Attractors}
	\label{sec:random}

 We denote by $W=(W_{1},W_{2}, W_{3})$ the $\mathcal{V}$-valued Wiener process, which has a version $\omega$ in $ C_{0}(\R,\mathcal V):=\wt{\Omega}$, the space of continuous functions which are $\mathrm{zero}$ at $\mathrm{zero}.$ In what follows we consider a canonical version of $W$ given by the probability space $(C_{0}(\R,\mathcal V), B(C_{0}(\R,\mathcal V)), \wt \P )$ where $\wt \P$ is the Wiener-measure generated by $W.$
On this probability space we can also introduce the shift
	$$\vartheta_s\omega(t)=\omega(t+s)-\omega(s),\hspace{1cm}s,t\in\R.$$
	\par
Going back to the abstract evolution system defined in (\ref{eqn:abstractmoist}),
\begin{equation}
\label{eqn:abstractmoist2}
dU+(AU+B(U)+A_pU+EU)dt=Fdt+dW,\ U(t_0)=U_0,
\end{equation}
	we define an Ornstein-Uhlenbeck process by
	$$Z(t)=\int_{-\infty}^te^{-(A+\gamma)(t-s)}dW(s),\hspace{1cm}\wh U=U-Z.$$
	$Z$ is a stationary process and its trajectories are $\wt \P$-a.s. continuous. $\wh U$ satisfies another evolution system:
	\begin{equation}
	\label{eqn:Z}
	\frac{d\wh U}{dt}+A\wh U+B(\wh U+Z)+A_p(\wh U+Z)+E(\wh U+Z)=F+\gamma Z-AZ.
	\end{equation}
	Again, using Galerkin approximating method, we have for any $\omega\in\wt\Omega$, and any fixed $s\in\R$, and $\wh U_s\in \mathcal V$,  a.s., there exists a unique solution, $\wh U(t,\omega)$,  defined on $[s,\infty)$, satisfying the above equation and
	\begin{equation}
	\label{eqn:Zs}
	\wh U(s,\omega)=\wh U_s(\omega),\ \wt\P\text{-a.s.}
	\end{equation}
	We now define the stochastic dynamical system $(\S(t,s;\omega))_{t\geq s,\omega\in\wt\Omega}$ by
	$$\S(t,s;\omega)U_s=\wh U(t,\omega)+Z(t,\omega),$$
	with $\wh U(s,\omega)=\wh U_s(\omega)=U_s-Z_s(\omega)$. It's obvious that $\S(t,s;\omega)$ satisfies (i)-(iv)(see Section 2), and also satisfies for any $s<t$ and $h\in \mathcal V$,
	$$\S(t,s;\omega)h=\S(t-s,0;\vartheta_s\omega)h,\ \wt\P\text{-a.s.}$$

To prove the compact property of solution operator, we need Aubin's Lemma which is cited below.
	\begin{lemma}[Aubin's Lemma]
		Let $B_0, B, B_1$ be Banach spaces such that $B_0, B_1$ are reflexive and $B_0\overset{c}{\subset}B\subset B_1$. Define for $0<K<\infty$,
		$$X:=\{h|h\in L^2([0,K],B_0),h'(t)\in L^2([0,K];B_1)\}.$$
		Then $X$ is a Banach space equipped with the norm $|h|_{L^2([0,K];B_0)}+|h|_{L^2([0,K];B_1)}$. Moreover,
		$$X\overset{c}{\subset}L^2([0,K];B_1).$$
	\end{lemma}
	Finally, we restate our main result of the existence of the random attractor for stochastic PEs.
	\begin{theorem}
	Let $Q_T$, $\partial_{\xi}Q_{T}$, $Q_q$, $\partial_{\xi}Q_{q}\in L^{2}(\D),$ $(\mathbf{v}_{0}, T_{0}, q_{0})\in \mathcal{V}$, and $\sum_{i=1}^{\infty}\lambda_{i,j}^{2}\gamma_{i,j}^{2+\sigma}<\infty,j=1,2,3,$ for small positive constant $\sigma$. Furthermore we assume $|\frac{br_{s}}{r}|\leq \min\{\frac{1}{2},\alpha,\beta\}$. Then	the solution operator $\left(\S(t,s;\omega)\right)_{t\geq s,\omega\in\wt\Omega}$ of the 3D stochastic PEs (\ref{eqn:moist1}), defined as $\S(t,s;\omega)(\v_s,T_s,q_s)=(\v_t,T_t,q_t)$ has properties $\mathrm{(i)}-\mathrm{(iv)}$ and possesses a compact absorbing set $\mathcal B(0,\omega)$ at time $t$.
		Moreover, for $\wt\P$-a.s. $\omega$, the set $\mathcal A(\omega)=\wb{\underset{B\subset\mathcal V}{\bigcup}\mathcal A(B,\omega)}$ is the random attractor of stochastic PEs, where the union is taken over all the bounded subsets of $\mathcal{V}.$
	\end{theorem}
\proof
Following the classic arguments (see \cite{DZ}), we can prove that for arbitrary small  $\varepsilon>0$, we can choose $\gamma$ big enough such that $E\|Z_{i}(0)\|_{3}^{2}\leq \varepsilon, i=1,2,3 $  and $\|Z_{i}(t)\|_{3}$ has polynomial growth when $t\rightarrow -\infty$.  Furthermore, the process $Z(t)$ is stationary and ergodic, thus, we know from the ergodic theorem that
\begin{equation}
-\frac{1}{s}\int_s^0(\norm{Z_1}_3^2+\norm{Z_2}_3^2+\norm{Z_3}_3^2)dr\to\Ex[\norm{Z_1(0)}_3^2+\norm{Z_2(0)}_3^2+\norm{Z_3(0)}_3^2]\text{ as }s\to-\infty.
\end{equation}
Since we can choose $\gamma$ big enough such that
	$$\Ex[\norm{Z_1(0)}_3^2+\norm{Z_2(0)}_3^2+\norm{Z_3(0)}_3^2]\leq \frac{\gamma_1}{4},$$
 there exists $s_0(\omega)$ such that for any $s<s_0(\omega)$.
$$-\frac{1}{s}\int_s^0(\norm{Z_1}_3^2+\norm{Z_2}_3^2+\norm{Z_3}_3^2)dr\leq \frac{\gamma_1}{4}.$$
	Using similar discussion with respect to negative time $t$, and referring to the result (\ref{eqn:gronwall}), we have
	\begin{align}
	\label{eqn:L2negative}
	&|\bu(-4)|_2^2+|S(-4)|_2^2+|p(-4)|_2^2\notag\\
	 &\leq(|\v(t_0)|_2^2+|T(t_0)|_2^2+|q(t_0)|_2^2)\exp\bigg[C\int_{t_0}^{-4}(-\gamma_1+\norm{Z_1(s)}_3^2+\norm{Z_2(s)}_3^2+\norm{Z_3(s)}_3^2)ds\bigg]\notag\\
	 &\hspace{5mm}+\int_{t_0}^{-4}(|Q_T|_2^2+|Q_q|_2^2+|Z|_2^2)\exp\bigg[C\int_{s}^{-4}(-\gamma_1+\norm{Z_1(r)}_3^2+\norm{Z_2(r)}_3^2+\norm{Z_3(r)}_3^2)dr\bigg]ds.
	\end{align}
	
	\par
	Applying (\ref{eqn:gronwall}) again, we have  for $t\in [-4,0]$
	\begin{align}
	\label{eqn:L2integral}
	&|\bu(t)|_2^2+|S(t)|_2^2+|p(t)|_2^2\notag\\
	&\leq(|\bu(-4)|_2^2+|S(-4)|_2^2+|p(-4)|_2^2)\exp\bigg[C\int_{-4}^{t}(-\gamma_1+\norm{Z_1(s)}_3^2+\norm{Z_2(s)}_3^2+\norm{Z_3(s)}_3^2)ds\bigg]\notag\\
	 &\hspace{5mm}+\int_{-4}^{t}(|Q_T|_2^2+|Q_q|_2^2+|Z|_2^2)\exp\bigg[C\int_{s}^{-4}(-\gamma_1+\norm{Z_1(r)}_3^2+\norm{Z_2(r)}_3^2+\norm{Z_3(r)}_3^2)dr\bigg]ds.
	\end{align}
	Now we denote by $(\bu(t,\omega;t_0,\bu_*), S(t,\omega;t_0,S_*), p(t,\omega;t_0,p_*)) $ the solution to the system (\ref{eqn:system}) with $(\bu(t_0),S(t_0),p(t_0))=(\bu_*,S_*,p_*)$.
	Then, by (\ref{eqn:L2negative}) and (\ref{eqn:L2integral}), there exists $r_1(\omega)$ depending on $\gamma_1, Z_1, Z_2$ and $Z_3$, such that for arbitrarily fixed $\rho>0$ there exists $t(\omega)\leq -4$, $\P$-as for all $t_0\leq t(\omega)$ and $(\bu_*,S_*,p_*)\in\mathcal V$ with $\|\bu_*\|_1+\|S_*\|_1+\|p_*\|_1<\rho$, the solution $(\bu(t,\omega;t_0,\bu_*), S(t,\omega;t_0,S_*), p(t,\omega;t_0,p_*)) $ on $[t_0,\infty)$ satisfies
	\begin{equation}
		\label{eqn:L2ofuSp}
	|\bu(t,\omega;t_0,\bu_*)|_2^2+|S(t,\omega;t_0,S_*)|_2^2+|p(t,\omega;t_0,p_*)|_2^2\leq r_1(\omega)\text{ for }t\in[-4,0].
	\end{equation}	
	Moreover, integrating (\ref{eqn:diff}), we have
	\begin{align}
		&\int_{-4}^{0}(\norm{\bu}_1^2+\norm{S}_1^2+\norm{p}_1^2)ds\notag\\
		&\leq |\bu(-4)|_2^2+|S(-4)|_2^2+|p(-4)|_2^2+ C\int_{-4}^0(|\bu|_2^2+|S|_2^2+|p|_2^2)(\norm{Z_1(s)}_3^2+\norm{Z_2(s)}_3^2+\norm{Z_3(s)}_3^2))ds\notag\\
		&\hspace{5mm}+C\int_{-4}^{0}(|Q_T|_2^2+|Q_q|_2^2+|Z|_2^2)ds,
	\end{align}
	thus, there exists $c_1(\omega)$ depending on $\gamma_1, Z_1, Z_2$ and $Z_3$,
	\begin{equation}
	\label{eqn:normfromneg}
		\int_{-4}^{0}(\norm{\bu(t,\omega;t_0,\bu_*)}_1^2+\norm{S(t,\omega;t_0,S_*)}_1^2+\norm{p(t,\omega;t_0,p_*)}_1^2)ds\leq c_1(\omega).
	\end{equation}
	We continue the discussion of $L^4$ norms. For $t<-3,$ by (\ref{eqn:gronwallofL4p}) we have
	\begin{align}
		|p(-3,\omega;t_0,p_*)|_4^2\leq& |p(t,\omega;t_0,p_*)|_4^2e^{-C(-3-t)}\notag\\
		&+C\int_t^{-3}e^{-C(-3-s)}(|Q_q|_2^{2}+\norm{Z_3}_3^{2}+\norm{Z_1}_1^{2}+\norm{Z_3}_3^{2}\norm{\bu(s,\omega;t_0,\bu_*)}_1^{2})ds.
	\end{align}
	We now integrate both sides over $[-4,-3]$,
	\begin{align}
		&|p(-3,\omega;t_0,p_*)|_4^2\notag\\
		&\leq\int_{-4}^{-3} |p(t,\omega;t_0,p_*)|_4^2e^{-C(-3-t)}dt\notag\\
		 &\hspace{5mm}+C\int_{-4}^{-3}\int_t^{-3}e^{-C(-3-s)}(|Q_q|_2^{2}+\norm{Z_3}_3^{2}+\norm{Z_1}_1^{2}+\norm{Z_3}_3^{2}\norm{\bu(s,\omega;t_0,\bu_*)}_1^{2})dsdt\notag\\
		&\leq C\int_{-4}^{-3}(\norm{p(t,\omega;t_0,p_*)}_1^2+\norm{\bu(t,\omega;t_0,\bu_*)}_1^2)dt+C\int_{-4}^{-3}e^{-C(-3-s)}(|Q_q|_2^{2}+\norm{Z_3}_3^{2}+\norm{Z_1}_1^{2})dt.
	\end{align}
	By (\ref{eqn:normfromneg}), there exists  $c_{2}(\omega)$, depending only on
	$\gamma_{1}, Z_{1}, Z_{2}, Z_3$  such that for arbitrarily fixed $\rho>0$ there exists $t(\omega)\leq -3$, $\wt\P$-a.s. for all $t_0\leq t(\omega)$ and $(\bu_*,S_*,p_*)\in\mathcal V$ with $\|\bu_*\|_1+\|S_*\|_1+\|p_*\|_1<\rho$, and the solution $p(t,\omega;t_0,p_*)$ satisfies
	\begin{equation}
		|p(-3,\omega;t_0,p_*|_{4}\leq c_2(\omega).
	\end{equation}
	Similar as in (\ref{eqn:L2ofuSp}), there exists $c_3(\omega)$depending only on
	$\gamma_{1}, Z_{1}, Z_{2}, Z_3$  such that for arbitrarily fixed $\rho>0$ there exists $t(\omega)\leq -3$, $\wt\P$-a.s. for all $t_0\leq t(\omega)$ and $(\bu_*,S_*,p_*)\in\mathcal V$ with $\|\bu_*\|_1+\|S_*\|_1+\|p_*\|_1<\rho$, and the solution $p(t,\omega;t_0,p_*)$ satisfies
	\begin{equation}
		|p(t,\omega;t_0,p_*)|_4^2+ \int_{-3}^{0}|p|_{\xi=1}|_{4}^{4}ds\leq c_3(\omega),\text{ for any }t\in[-3,0].
	\end{equation}
	Analogously, by \eqref{eqn:5.23} and \eqref{eqn:5.24}, we conclude that there exists  $c_4(\omega)$ depending on $\gamma_{1}, Z_{1}, Z_{2}, Z_3$  such that for arbitrarily fixed $\rho>0$ there exists $t(\omega)\leq -3$, $\wt\P$-a.s. for all $t_0\leq t(\omega)$ and $(\bu_*,S_*,p_*)\in\mathcal V$ with $\|\bu_*\|_1+\|S_*\|_1+\|p_*\|_1<\rho$, and the solution $S(t,\omega;t_0,S_*)$ satisfies
	\begin{equation}
		|S(t,\omega;t_0,S_*)|_4^2+ \int_{-3}^{0}|S|_{\xi=1}|_{4}^{4}\leq c_4(\omega), \text{ for any }t\in[-3,0].
	\end{equation}
	By  \eqref{eqn:5.39}, \eqref{eqn:5.40}, and with similar discussion as the above, there exist constants $c_5(\omega)$, $c_6(\omega)$ depending on $\gamma_{1}, Z_{1}, Z_{2}, Z_3$  such that for arbitrarily fixed $\rho>0$ there exists $t(\omega)\leq-3$, $\wt\P$-a.s. for all $t_0\leq t(\omega)$ and $(\bu_*,S_*,p_*)\in\mathcal V$ with $\|\bu_*\|_1+\|S_*\|_1+\|p_*\|_1<\rho$, and the solution $\wt\bu(t,\omega;t_0,\wt\bu_*)$ satisfies
	\begin{align}
	\label{eqn:L4ofwtuneg}
		&|\wt\bu(t,\omega;t_0,\wt\bu_*)|_4^2\leq c_5(\omega),\text{ for any }t\in[-3,0],\\
		\label{eqn:L2ofwtuandder}
		&\int_{-3}^{0}\int_{\mathbf D}|\wt\bu|^2\Big(|\nabla_{e_\theta}\wt\bu|^2+|\nabla_{e_\varphi}\wt\bu|^2+|\partial_\xi\wt\bu|^2\Big)d\mathbf Dds\leq c_6(\omega).
	\end{align}
	By \eqref{eqn:5.43}, we proceed to have the existence of a constant $c_7(\omega)$ depending on $\gamma_{1}, Z_{1}, Z_{2}, Z_3$  such that for arbitrarily fixed $\rho>0$ there exists $t(\omega)\leq -2$, $\wt\P$-a.s. for all $t_0\leq t(\omega)$ and $(\bu_*,S_*,p_*)\in\mathcal V$ with $\|\bu_*\|_1+\|S_*\|_1+\|p_*\|_1<\rho$, and the solution $\bar\bu(t,\omega;t_0,\bar\bu_*)$ satisfies
	\begin{equation}
	\label{eqn:L2ofbaruder}
		\int_{\mathbf S^2}|\nabla_{e_\theta}\bar\bu|^2+|\nabla_{e_\varphi}\bar{\bu}|^2d\mathbf S^2\leq c_7(\omega),\text{ for any }t\in[-2,0].
	\end{equation}
	By (\ref{eqn:H1ofv}), and together with results in (\ref{eqn:normfromneg}) and (\ref{eqn:L4ofwtuneg}), (\ref{eqn:L2ofwtuandder}), (\ref{eqn:L2ofbaruder}), we conclude that there exist constants $r_2(\omega)$ and $c_{8}(\omega)$ depending on $\gamma_{1}, Z_{1}, Z_{2}, Z_3$  such that for arbitrarily fixed $\rho>0$ there exists $t(\omega)\leq-1$, $\wt\P$-a.s. for all $t_0\leq t(\omega)$ and $(\bu_*,S_*,p_*)\in\mathcal V$ with $\|\bu_*\|_1+\|S_*\|_1+\|p_*\|_1<\rho$, $\bu_{\xi}(t,\omega;t_0,\bu_*)$ satisfies
	\begin{align}
		&|\bu_{\xi}(t,\omega;t_0,\bu_*)|_2^2\leq r_2(\omega),\text{ for any }t\in[-1,0],\\
		&\int_{-1}^0\int_{\mathbf S^2}|\nabla_{e_\theta}\bu_{\xi}(s,\omega;t_0,\bu_*)|^2+|\nabla_{e_\varphi}\bu_{\xi}(s,\omega;t_0,\bu_*)|^2+ |\bu_{\xi\xi}(s,\omega;t_0,\bu_*)|^2 d\mathbf S^2ds\leq c_{8}(\omega).
	\end{align}

By \eqref{eqn:5.77}, and together with results in (\ref{eqn:normfromneg}) and (\ref{eqn:L4ofwtuneg}), (\ref{eqn:L2ofwtuandder}), (\ref{eqn:L2ofbaruder}), $(4.18)$ and $(4.19)$ we conclude that there exist constants $r_3(\omega)$ and $c_{9}(\omega)$ depending on $\gamma_{1}, Z_{1}, Z_{2}, Z_3$  such that for arbitrarily fixed $\rho>0$ there exists $t(\omega)\leq-1$, $\wt\P$-as for all $t_0\leq t(\omega)$ and $(\bu_*,S_*,p_*)\in\mathcal V$ with $\|\bu_*\|_1+\|S_*\|_1+\|p_*\|_1<\rho$, $p_{\xi}(t,\omega;t_0,p_*)$ satisfies
	\begin{align}
		&|p_{\xi}(t,\omega;t_0,p_*)|_2^2\leq r_3(\omega),\text{ for any }t\in[-1,0],\\
		&\int_{-1}^0(|\nabla p_{\xi}|_{2}^{2}+|p_{\xi\xi}|_{2}^{2})ds\leq c_{9}(\omega).
	\end{align}

By \eqref{eqn:5.83}, and together with results in (\ref{eqn:normfromneg}) and (\ref{eqn:L4ofwtuneg}), (\ref{eqn:L2ofwtuandder}), (\ref{eqn:L2ofbaruder}), $(4.18)-(4.21)$ we conclude that there exist constants $r_4(\omega)$ and $c_{10}(\omega)$ depending on $\gamma_{1}, Z_{1}, Z_{2}, Z_3$  such that for arbitrarily fixed $\rho>0$ there exists $t(\omega)\leq-1$, $\wt\P$-a.s. for all $t_0\leq t(\omega)$ and $(\bu_*,S_*,p_*)\in\mathcal V$ with $\|\bu_*\|_1+\|S_*\|_1+\|p_*\|_1<\rho$, $ S_{\xi}(t,\omega;t_0,p_*)$ satisfies
	\begin{align}
		&| S_{\xi}(t,\omega;t_0,S_*)|_2^2\leq r_4(\omega),\text{ for any }t\in[-1,0],\\
		&\int_{-1}^0(|\nabla S_{\xi}|_{2}^{2}+|S_{\xi\xi}|_{2}^{2})ds\leq c_{10}(\omega).
	\end{align}

In view of \eqref{eqn:5.92}, there exist constant $r_5(\omega)$ such that for arbitrarily fixed $\rho>0$, there exists $t(\omega)\leq-1$, $\wt\P$-a.s. for all $t_0\leq t(\omega)$ and $(\bu_*,S_*,p_*)\in\mathcal V$ with $\|\bu_*\|_1+\|S_*\|_1+\|p_*\|_1<\rho$, $\bu(t,\omega;t_0,\bu_*), S$ and $ p$ satisfy
	\begin{align}	
&|\nabla_{e_\theta}\bu(t,\omega;t_0,\bu_*)|_2^2+|\nabla_{e_\varphi}\bu(t,\omega;t_0,\bu_*)|_2^2\notag\\
&\hspace{5mm}+|\nabla S(t,\omega;t_0,S_*)|_2^2+|\nabla p(t,\omega;t_0,p_*)|_2^2\leq r_5(\omega), \text{ for any }t\in[-1,0],
	\end{align}
which together with $(4.7), (4.18), (4.20)$ and $(4.22)$ imply that there exists a constant $r_6(\omega)$ such that for arbitrarily fixed $\rho>0$, there exists $t(\omega)\leq-1$, $\wt\P$-a.s. for all $t_0\leq t(\omega)$ and $(\bu_*,S_*,p_*)\in\mathcal V$ with $\|\bu_*\|_1+\|S_*\|_1+\|p_*\|_1<\rho$, $(\bu, S, p)$ satisfy
	\begin{equation}
\norm{\bu(t,\omega;t_0,\bu_*)}_1^2+		\norm{S(t,\omega;t_0,S_*)}_1^2+\norm{p(t,\omega;t_0,p_*)}_1^2\leq r_6(\omega),\text{ for any }t\in[-1,0].
	\end{equation}
	Now we are ready to prove the desired compact result.
	Let $r(\omega)= r_{6}(\omega)+\norm{Z(-1)}_1^2$, then  $B(-1,r(\omega))$, the ball of center $0\in \mathcal{V}$ and radius $r(\omega),$ is a absorbing set
	at time $-1$ for $(\S(t,s;\omega))_{t\geq s,\omega\in \tilde{\Omega}}$. Therefore, in order to prove the existence of the random attractor of  the stochastic dynamical system in space $\mathcal{V}$, we need to to construct a compact absorbing set at time $0$ in $\mathcal{V}$. Denote by $\mathcal{B}$  a bounded subset of $\mathcal{V}$ and set $\mathcal{C}_{T}$ as a subset of the function space:
	\begin{align}
		\mathcal{C}_{T,q}:=\Big\{(A_1^{1/2}\v,A_2^{1/2}T,A_3^{1/2}q)\Big|&(\v(-1),T(-1),q(-1))\in\mathcal B,\notag\\
		&(\v(t),T(t),q(t))=\S(t,-1;\omega)(\v(-1),T(-1),q(-1)),t\in[-1,0].
		\Big\}
	\end{align}

	Obviously the embedding $\mathcal{V} \subset \mathcal{H}$ is compact.  Let $(\v(-1), T(-1),q(-1))\in \mathcal{B} ,$ by the continuity of strong solutions with respect to time $t$, we know
	\begin{align*}
		&(A_1^{1/2}\bu,A_2^{1/2}S,A_3^{1/2}p)\in L^2([-1,0];\mathcal V_1\times\mathcal V_2\times\mathcal V_3),\\
		&(\partial_tA_1^{1/2}\bu,\partial_tA_2^{1/2}S,\partial_tA_3^{1/2}p)\in L^2([-1,0];\mathcal V'_1\times\mathcal V'_2\times\mathcal V'_3).
	\end{align*}
	Now we apply Aubin's Lemma with
	$$
	B_{0}=\mathcal{V}_{1}\times \mathcal{V}_{2}\times\mathcal V_3,\ \ B=\mathcal H_{1}\times (L^2(\D))^2,\ \ B_{1}=\mathcal{V}_{1}'\times \mathcal{V}_{2}'\times\mathcal V_3',
	$$
	$\mathcal{C}_{T,q}$ is compact in $L^{2}([-1,0];\mathcal{H} ).$\\
	In order to show that for any fixed $t\in (-1,0],\omega\in \tilde{\Omega}, \S(t,-1;\omega) $ is a compact operator in $\mathcal{V} ,$ we take any bounded sequences
	$\{(\v_{0,n}, T_{0,n} ,q_{0,n}) \}_{n\in \mathbb{N}}\subset \mathcal{B}$ and we want to extract, for any fixed $t\in (-1,0]$ and $\omega\in \tilde{\Omega},$ a convergent subsequence from
	$\{\S(t,-1;\omega)( \v_{0,n}, T_{0,n},q_{0,n}) \}$.
	Since $\{(A_{1}^{\frac{1}{2}}\v, A_{2}^{\frac{1}{2}}T, A_3^{\frac{1}{2}}q  ) \}\subset \mathcal{C}_{T,q} ,$ by Aubin's Lemma, there is a function $(\v_{*},T_{*}, q_*)$:
	\[
	(\v_{*},T_{*},q_*)\in  L^{2}([-1,0];\mathcal{V}),
	\]
	and there exists a subsequence of $\{\S(t,-1;\omega)(\v_{0,n},T_{0,n},q_{0,n} ) \}_{n\in \mathbb{N}},$ for simplicity, we still denote it by  $\{\S(t,-1;\omega)(\v_{0,n},T_{0,n},q_{0,n} ) \}_{n\in \mathbb{N}}$, and it satisfies
	\begin{equation}
		\lim\limits_{n\rightarrow \infty}\int_{-1}^{0}\|\S(t,-1;\omega)(\v_{0,n},T_{0,n},q_{0,n} )-(\v_{*}(t),T_{*}(t), q_*(t))  \|_{1}^{2}dt=0.
	\end{equation}
	By elementary measure theory, there exists a subsequence of $\{\S(t,-1;\omega)(\v_{0,n},T_{0,n},q_{0,n} ) \}_{n\in \mathbb{N}}, $
	still denoted by $\{\S(t,-1;\omega)(\v_{0,n},T_{0,n},q_{0,n} ) \}_{n\in \mathbb{N}}$ for simplicity , such that
	\begin{equation}
		\lim\limits_{n\rightarrow \infty}\|\S(t,-1;\omega)(\v_{0,n},T_{0,n},q_{0,n} )-(\v_{*}(t),T_{*}(t), q_*(t))  \|_{ 1}=0,\ \ a.e.\ t\in (-1,0].
	\end{equation}
	Fix any $t\in (-1,0]$, we can select a $t_{0}\in (-1,t)$ such that
	$$
	\lim\limits_{n\rightarrow \infty}\|\S(t_{0},-1,\omega)(\v_{0,n},T_{0,n},q_{0,n} )-(\v_{*}(t_{0}),T_{*}(t_{0}),q_*(t))  \|_{ 1}=0.
	$$
	Then by the continuity of the map $\S(t-t_{0},t_{0};\omega) $ in $ \mathcal{V}$ with respect to initial value, we have
	\begin{align*}
		\S(t,-1;\omega)(\v_{0,n},T_{0,n},q_{0,n})=&\S(t-t_{0},t_{0};\omega)\S(t_{0},-1;\omega)(\v_{0,n},T_{0,n},q_{0,n})\\
		\to & \S(t-t_{0},t_{0};\omega)(\v_{*}(t_{0}),T_{*}(t_{0}),q_*(t_0))\ \ \ \text{in}\ \mathcal{V}.
	\end{align*}
	Hence for any $t\in(-1,0], \{\S(t,-1;\omega)(\v_{0,n},T_{0,n},q_{0,n} ) \}_{n\in \mathbb{N}}$   contains a subsequence which is convergent $\mathrm{in}\ \mathcal{V},$ which implies that for any fixed $t\in (-1,0],\omega\in \tilde{\Omega}, \S(t,-1;\omega) $ is a compact operator in $\mathcal{V}.$ Let $\mathcal{B}(0,\omega)= \overline{\S(0,-1;\omega)B(-1,r(\omega))}$ be the closed set of $ \S(0,-1;\omega)B(-1,r(\omega))$ in $\mathcal{V}.$ Then, by the above arguments, we know $\mathcal{B}(0,\omega)$ is a random compact set in $\mathcal{V}.$  More precisely,  $\mathcal{B}(0,\omega)$ is a compact absorbing set in $\mathcal{V}$ at time $0.$ Indeed, for $(\v_{0,n},T_{0,n},q_{0,n})\in \mathcal{B}, $ there exists $s(\mathcal{B})\in \mathbb{R}_{-}$ such that if $s\leq s(\mathcal{B}),$ we have
	\begin{align*}
			\S(0,s;\omega)(\v_{0,n},T_{0,n},q_{0,n})=&\S(0,-1;\omega)\S(-1,s;\omega)(\v_{0,n},T_{0,n},q_{0,n})\\
			\subset& \S(0,-1;\omega) B(-1,r(\omega))\subset \mathcal{B}(0,\omega) .
	\end{align*}
	Therefore, the existence of the random attractor follows.
\hspace{\fill}$\square$

		\section{Appendix: A Priori Estimates}
		
		\subsection{Decomposition}
		We first introduce a modified stochastic convolution, which is an Ornstein-Uhlenbeck process satisfying
		\begin{equation}
			\label{eqn:OU}
			dZ+(AZ+\gamma Z)dt=dW,\hspace{1cm}\text{with }\gamma>0.
		\end{equation}
		Denote $Z=(Z_1,Z_2,Z_3)$, and define $\wh{U}=U-Z=(\bu,S,p)$ which satisfies
		\begin{align}
			&\frac{\partial \wh U}{\partial t}+A\wh U+B(Z+\wh U)+E(Z+\wh  U)+A_p(Z+\wh U)=F+\gamma W,\notag\\
			&\wh U(t_0)=U_0-Z(t_0)=(\v_0,T_0,q_0).
		\end{align}
Using Kolmogorov test theorem we can get the regularity of $Z_{i},i=1,2,3$. We can also see the standard arguments in \cite{DZ}, so we omit the proof of the following lemma.
\begin{lemma}
Assume $\tau\geq t_{0}$ and $\sum_{i=1}^{\infty}\lambda_{i,j}\gamma_{i,j}^{2+\sigma}<\infty,j=1,2,3,$ for small positive constant $\sigma.$ Then $$Z_{1}\in C([t_0,\tau]; H^{3}(T\D|T\mathbf{S}^{2}))\ \mathrm{and}\
Z_{2}, Z_{3}\in C([t_0,\tau]; H^{3}(\D)).$$

\end{lemma}
		In details, the components of $\wh U$, $\bu, S, p$, satisfy
		\begin{subequations}
			\label{eqn:system}
			\begin{align}
				&\partial_t\bu+\nabla_{Z_1+\bu}(Z_1+\bu)+w(Z_1+\bu)\partial_\xi(Z_1+\bu)+\frac{f}{R_0}\vec{k}\times(Z_1+\bu)+\nabla\Phi_s\notag\\
				\label{eqn:u}
				&=\Delta\bu+\partial_{\xi\xi}\bu-\int_{\xi}^{1}\frac{br_s}{r}\nabla[(1+a(Z_3+p))(Z_2+S)]d\xi'+\gamma Z_1,\\
				&\partial_{t} S+\nabla_{Z_1+\bu}(Z_2+S)+w(Z_1+\bu)\partial_\xi(Z_2+S)\notag\\
				\label{eqn:S}
				&=\Delta S+\partial_{\xi\xi}S+\frac{br_s}{r}(1+a(Z_3+p))w(Z_1+\bu)+Q_T+\gamma Z_2,\\
				\label{eqn:p}
				&\partial_{t} p+\nabla_{Z_1+\bu}(Z_3+p)+w(Z_1+\bu)\partial_{\xi}(Z_3+p)=\Delta p+\partial_{\xi\xi}p+Q_q+\gamma Z_3,\\
				&w(Z_1+\bu)=\int_{\xi}^{1}\text{div }(Z_1+\bu)d\xi',\\
				&\int_0^1\text{div }\bu d\xi=0,\\
				&\text{On }\xi=1\ (r=r_s):\ \partial_\xi\bu=0, w=0, \partial_\xi S=-\alpha S, \partial_\xi p=-\beta p,\\
				&\text{On }\xi=0\ (r=r_0):\ \partial_\xi\bu=0, w=0, \partial_\xi S=0, \partial_\xi p=0,\\
				&(\bu(t_0),S(t_0),p(t_0))=(\bu_0,S_0,p_0).
			\end{align}
		\end{subequations}
		\subsection{$L^2$ estimates of $\bu, S, p$}
		We first take the inner product of equation (\ref{eqn:p}) with $p$, in $L^2(\mathbf D)$, we get
		\begin{align*}
			&\frac{1}{2}\frac{d|p|_2^2}{dt}+|\nabla p|_2^2+|\partial_\xi p|_2^2+\beta|p(\xi=1)|_2^2\\
			&=\int_{\mathbf D}(Q_q+\gamma Z_3)pd\mathbf D-\int_{\mathbf D}\Big(\nabla_{Z_1+\bu}(Z_3+p)+w(Z_1+\bu)\partial_{\xi}(Z_3+p)\Big)pd\mathbf D.
		\end{align*}
		By Lemma \ref{lemma:ibp} and integration by parts, we have
		$$\int_{\mathbf D}\Big(\nabla_{Z_1+\bu} p+w(Z_1+\bu)\partial_\xi p\Big)pd\mathbf D=0.$$
		Applying the H\"older's inequality and the Sobolev embedding theorem, we get
		\begin{align*}
			&\int_{\mathbf D}\left(w(Z_1+\bu)\partial_\xi Z_3+ \nabla_{Z_1+\bu} Z_3\right)pd\mathbf D\\
			&\leq C|\partial_\xi Z_3|_\infty|\text{div } Z_1+\text{div }\bu|_2|p|_2+C|\nabla Z_3|_\infty|Z_1+\bu|_2|p|_2\\
			&\leq\eps(|\nabla_{e_\theta}\bu|_2^2+|\nabla_{e_\varphi}\bu|^2+|\bu|_2^2)+C\norm{Z_3}_3^2|p|_2^2+C\norm{Z_1}_1^2.
		\end{align*}
		We also have
			$$|p|_2^2-|p(\xi=1)|_2^2\leq2\int_{\D}\int_0^1|p\partial_\xi p|d\D
			\leq\frac{1}{2}|p|_2^2+2|\partial_\xi p|_2^2.$$
		This gives $$|p|_2^2\leq 2|p(\xi=1)|_2^2+4|\partial_\xi p|_2^2,$$ similarly we also have
$$|S|_2^2\leq 2|S(\xi=1)|_2^2+4|\partial_\xi S|_2^2.$$
Therefore, we obtain that
		\begin{align*}
			\int_{\mathbf D}(Q_q+\gamma Z_3)pd\mathbf D\leq& C(|Q_q|_2^2+|Z_3|_2^2)+\frac{\eps}{16}|p|_2^2\\
			\leq& C(|Q_q|_2^2+|Z_3|_2^2)+\eps|\partial_\xi p|_2^2+\eps|p(\xi=1)|_2^2.
		\end{align*}
		Altogether, we have the estimate for $p$ as
		\begin{align}
			\label{eqn:diffofp}
			&\frac{1}{2}\frac{d|p|_2^2}{dt}+|\nabla p|_2^2+(1-\eps)|\partial_\xi p|_2^2+(\beta-\eps)|p(\xi=1)|_2^2\notag\\
			&\leq C(|Q_q|_2^2+|Z_3|_2^2)+\eps(|\nabla_{e_\theta}\bu|^2+|\nabla_{e_\varphi}\bu|^2+|\bu|_2^2)+C\norm{Z_1}_1^2+C\norm{Z_3}_3^2|p|_2^2.
		\end{align}
		Similarly, we get the estimate for $S$ from \eqref{eqn:S} as
\begin{align*}
			&\frac{1}{2}\frac{d|S|_2^2}{dt}+|\nabla S|_2^2+|\partial_\xi S|_2^2+\alpha|S(\xi=1)|_2^2\\
			&=\int_{\mathbf D}(Q_T+\gamma Z_2) S d\mathbf D-\int_{\mathbf D}\Big(\nabla_{Z_1+\bu} Z_2+w(Z_1+\bu)\partial_{\xi}Z_2\Big)Sd\mathbf D\\
            &+\int_{\mathbf D}\frac{br_{s}}{r}(1+a(Z_{3}+p))w(Z_1+\bu)Sd\mathbf D\\
            &=:I_{1}+I_{2}+I_{3}.
		\end{align*}
With same discussion, we have
$$I_1\leq C(|Q_T|_2^2+|Z_{2}|_{2}^{2})+\eps|\partial_\xi S|_2^2+\eps|S(\xi=1)|_2^2.$$
In view of the H\"older inequality, the Sobolev imbedding theorem and Young's inequality, we reach
\begin{align*}
            I_{2}&\leq \varepsilon (|\bu|_{2}^{2}+ |\nabla_{e_{\theta}}\bu |_{2}^{2}+  |\nabla_{e_{\phi}}\bu |_{2}^{2}+|S|_{2}^{2}  ) \\
            &+C\|Z_{2}\|_{3}^{2}|S|_{2}^{2}+C\|Z_{1}\|_{1}^{2}.
	\end{align*}
Similarly, we deduce
\begin{align*}
            I_{3}=&\int_{\mathbf{D}}\frac{b r_{s}}{r}(1+a(Z_{3}+p))w(Z_{1})S d\mathbf{D}+\int_{\mathbf{D}} \frac{b r_{s}}{r}Z_{3}w(u)S d\mathbf{D}\\
            &+\int_{\mathbf{D}} \frac{br_{s}}{r}(1+ap)w(\bu)S d\mathbf{D}\\
            \leq &\varepsilon(|S|_{2}^{2}+  |\nabla_{e_{\theta}}\bu |_{2}^{2}+ |\nabla_{e_{\phi}}\bu |_{2}^{2} +|p|_{2}^{2})\\
            &+C\|Z_{1}\|_{1}^{2}+C\|Z_{3}\|_{2}^{2}|S|_{2}^{2}+C\|Z_{1}\|_{3}^{2}|S|_{2}^{2}\\
             &+\int_{\mathbf{D}} \frac{br_{s}}{r}(1+ap)w(\bu)S d\mathbf{D}.
	\end{align*}
Therefore, combining the above estimates about $I_{1}-I_{3}$, we obtain
		\begin{align}
			\label{eqn:diffofS}
			&\frac{1}{2}\frac{d|S|_2^2}{dt}+(1-\eps)|\nabla S|_2^2+(1-\eps)|\partial_\xi S|_2^2+(\alpha-\eps)|S(\xi=1)|_2^2\notag\\
			&\leq C(|Q_T|_2^2+|Z_2|_2^{2})+\eps(|\nabla_{e_\theta}\bu|^2+|\nabla_{e_\varphi}\bu|^2+|\bu|_2^2)+C\norm{Z_1}_1^2\notag\\
			&\hspace{5mm}+C(\norm{Z_2}_3^2+\norm{Z_1}_3^2)|S|_2^2+\varepsilon|p|_2^{2}\nonumber\\
&\hspace{5mm}+\int_{\mathbf D}\frac{br_s}{r}(1+ap)w(\bu)Sd\mathbf D.
		\end{align}
		For $\bu$ in \eqref{eqn:u}, since $(\frac{f}{R_0}\times\bu)\cdot\bu=0$, and by Lemma \ref{lemma:ibp} and integration by parts,
		\begin{align}
			\label{eqn:diffofu}
			 &\frac{1}{2}\frac{d|\bu|_2^2}{dt}+(1-2\eps)\left(|\nabla_{e_\theta}\bu|_2^2+|\nabla_{e_\varphi}\bu|_2^2\right)+|\partial_\xi\bu|_2^2+|\bu|_2^2\notag\\
			&\leq C\norm{Z_1}_1^2+C\norm{Z_1}_3^2|\bu|_2^2+\varepsilon|\bu|_{2}^{2}-\int_{\mathbf D}\bigg(\int_{\xi}^{1}\frac{br_s}{r}\nabla[(1+aq)T]d\xi'\bigg)\cdot\bu d\mathbf D.
		\end{align}
By integration by parts,
    \begin{align}
    \label{eqn:5.7}
\int_{\mathbf D}\frac{br_s}{r}(1+aq)Tw(\bu)d\mathbf D&= -\int_{\mathbf D}\nabla\left[\frac{br_s}{r}(1+aq)T\right]\int_{\xi}^{1}ud\xi'd\mathbf D\nonumber\\
&=\int_{\mathbf D}\Big{(}\partial_{\xi}\int_{\xi}^{1} \nabla\left[\frac{br_s}{r}(1+aq)T\right]d\xi' \Big{)}\int_{\xi}^{1}ud\xi'  d\mathbf D \nonumber\\
&=\int_{\mathbf D}\Big{(}\int_{\xi}^{1} \nabla\left[\frac{br_s}{r}(1+aq)T\right]d\xi' \Big{)}\bu d\mathbf D\nonumber.
     \end{align}
Therefore, we estimate the sum of the last term on the right hand side of \eqref{eqn:diffofS} and \eqref{eqn:diffofu} that
   \begin{align}
&\int_{\mathbf D}\frac{br_s}{r}(1+ap)w(\bu)Sd\mathbf D -\int_{\mathbf D}\Big{(}\int_{\xi}^{1} \frac{br_s}{r}\nabla[(1+aq)T]d\xi' \Big{)}\cdot\bu d\mathbf D\nonumber\\
&=\int_{\mathbf D}\frac{br_s}{r}(1+aq)w(\bu)Z_{2}d\mathbf D+ \int_{\mathbf D}\frac{br_s}{r}(1+aZ_{3})w(\bu)Td\mathbf D\nonumber\\
&\leq \varepsilon (|\nabla_{e_{\theta}}\bu|_{2}^{2}+|\nabla_{e_{\varphi}}\bu|_{2}^{2}  )+\frac{br_{s}}{2r}(|\nabla_{e_{\theta}}\bu|_{2}^{2}+|\nabla_{e_{\varphi}}\bu|_{2}^{2}  )\nonumber\\
&+C|Z_{2}|_{\infty}^{2}(1+|q|_{2}^{2})+|T|_{2}^{2}(\frac{br_{s}}{2r}+C|Z_{3}|_{\infty}^{2})\nonumber\\
&\leq \varepsilon (|\nabla_{e_{\theta}}\bu|_{2}^{2}+|\nabla_{e_{\varphi}}\bu|_{2}^{2}  )+\frac{br_{s}}{2r}(|\nabla_{e_{\theta}}\bu|_{2}^{2}+|\nabla_{e_{\varphi}}\bu|_{2}^{2}  +|S|_{2}^{2})\nonumber\\
&+C\|Z_{2}\|_{2}^{2}(1+|p|_{2}^{2})+C\|Z_{2}\|_{2}^{2}\|Z_{3}\|_{2}^{2}.
     \end{align}
\begin{remark}
    For the deterministic moist primitive equations, the authors use the uniform Gronwall inequality intensively to obtain the uniform estimates (see \cite{GH11}). But, for the stochastic moist primitive equations, the global attractor is a pullback attractor. In addition, Ornstein-Uhlenbeck	 process has the properties of polynomial growth,  the uniform Gronwall inequality is not valid. Therefore, in order to overcome the difficulties,  the first way is to use Young's inequality to estimate the nonlinear term; the second way is to take advantage of the ergodicity of Ornstein-Uhlenbeck process to make the coefficient of the term be small enough that the term can be absorbed by the dissipative terms of the equations.  Therefore, in our paper, to study the existence of random attractor we assume $\frac{br_{s}}{r}$ to be smaller than the minimum of $\{\frac{1}{2}, \alpha , \beta\}$. If we only want to study the global well-posedness of the strong solution to stochastic PEs, we do not need the assumption.
\end{remark}
		Combining \eqref{eqn:diffofp}-\eqref{eqn:5.7}, we get
		\begin{align}
			\label{eqn:diff}
			 \frac{1}{2}\frac{d(|\bu|_2^2+|S|_2^2+|p|_2^2)}{dt}+&(1-\varepsilon)(|\nabla_{e_\theta}\bu|_2^2+|\nabla_{e_\varphi}\bu|_2^2+|\bu|_2^2+|\nabla S|_2^2+|\nabla p|_2^2)\nonumber\\
&+(1-\varepsilon)(|\partial_\xi\bu|_{2}^{2}+|\partial_\xi S|_{2}^{2}+|\partial_\xi p|_2^2)\notag\\
			&+(\alpha-\eps)|S(\xi=1)|_2^2+(\beta-\eps)|p(\xi=1)|_2^2\notag\\
			\leq
&\frac{br_{s}}{2r}(|\nabla_{e_{\theta}}\bu|_{2}^{2}+|\nabla_{e_{\varphi}}\bu|_{2}^{2}  +|S|_{2}^{2})\nonumber\\
&+C(|\bu|_2^2+|S|_2^2+|p|_2^2)(\norm{Z_1}_3^2+\norm{Z_2}_3^2+\norm{Z_3}_3^2)\nonumber\\
&+C(|Q_T|_2^2+|Q_q|_2^2+\|Z_1\|_1^2+|Z_2|_2^2+|Z_3|_2^2+\|Z_{2}\|_{2}^{2}\|Z_{3}\|_{2}^{2}).
		\end{align}
	With the assumption that $\frac{br_s}{2r}$ is small enough,	there exists $\gamma_1>0$ such that
		\begin{align*}
			&(2-2\eps-\frac{br_s}{r})(|\nabla_{e_\theta}\bu|^2+|\nabla_{e_\varphi}\bu|^2)+(2-2\eps)(|\nabla S|_2^2+|\nabla p|_2^2)+(2-2\eps)|\partial_\xi\bu+\partial_\xi S+\partial_\xi p|_2^2\\
			&+(2\alpha-2\eps)|S(\xi=1)|_2^2+(2\beta-2\eps)|p(\xi=1)|_2^2> (\gamma_1-2+2\eps)(|\bu|_2^2+|S|_2^2+|p|_2^2).
		\end{align*}
		Therefore, by the Gronwall inequality we have for any $t$, there exists $\gamma_2$
		\begin{align}
			\label{eqn:gronwall}
			&|\bu(t)|_2^2+|S(t)|_2^2+|p(t)|_2^2\notag\\
			 &\leq(|\bu_0|_2^2+|S_0|_2^2+|p_0|_2^2)\exp\bigg[-\gamma_1t+\gamma_2\int_{t_0}^t(\norm{Z_1(s)}_3^2+\norm{Z_2(s)}_3^2+\norm{Z_3(s)}_3^2)ds\bigg]\notag\\
			 &\hspace{5mm}+\int_{t_0}^t(|Q_T|_2^2+|Q_q|_2^2+\|Z_1\|_1^2+|Z_2|_2^2+|Z_3|_2^2+\|Z_{2}\|_{2}^{2}\|Z_{3}\|_{2}^{2})\notag\\
			 &\hspace{15mm}\times\exp\bigg[-\gamma_1(t-s)+\gamma_2\int_s^t(\norm{Z_1(r)}_3^2+\norm{Z_2(r)}_3^2+\norm{Z_3(r)}_3^2)dr\bigg]ds.
		\end{align}
		\subsection{$L^4$ estimates of $p, S$.}
		Taking the inner product of Equation (\ref{eqn:p}) with $p^3$ in $L^2(\mathbf D)$, we obtain that
		\begin{align}
			\label{eqn:L4ofp}
			&\frac{1}{4}\frac{d|p|_4^4}{dt}+\frac{3}{4}|\nabla p^2|_2^2+\frac{3}{4}|\partial_\xi p^2|_2^2+\beta\int_{\mathbf S^2}|p(\xi=1)|^4d\mathbf S^2\notag\\
			&=\int_{\mathbf D}(Q_q+\gamma Z_3)p^3d\mathbf D-\int_{\mathbf D}\Big(\nabla_{Z_1+\bu}(Z_3+p)+w(Z_1+\bu)\partial_{\xi}(Z_3+p)\Big)p^3d\mathbf D.
		\end{align}
		Again, applying integration by parts, we have
		$$\int_{\mathbf D}\Big(\nabla_{Z_1+\bu} p+w(Z_1+\bu)\partial_\xi p\Big)p^3d\mathbf D=0.$$
		Applying the H\"older inequality and the Sobolev imbedding theorem, we get
		\begin{align}
			\label{eqn:L4ofp1}
			&\bigg|\int_{\mathbf D}(w(Z_1+\bu)\partial_\xi Z_3+\nabla_{Z_1+\bu} Z_3)p^3d\mathbf D\bigg|\notag\\
			&\leq C|\partial_\xi Z_3|_\infty|\text{div }\bu+\text{div } Z_1|_{2}|p^3|_{2}+C|\nabla Z_3|_\infty|Z_1+\bu|_2|p^3|_2\notag\\
            &\leq C\|Z_{3}\|_{3}(|\nabla_{e_\theta}\bu|_{2}+|\nabla_{e_\varphi}\bu|_{2}+|u|_{2} )|p^{2}|_{3}^{3/2}\notag\\
            &\hspace{5mm}+C\|Z_{3}\|_{3}\|Z_{1}\|_{3}|p^{2}|_{3}^{3/2}.
		\end{align}
Applying the interpolation inequality to $|p^{2}|_{3},$ we obtain
\[
|p^{2}|_{3}\leq C|p^{2}|_{2}^{\frac{1}{2}}(|\nabla p^{2} |_{2}^{\frac{1}{2}}+|  \partial_{\xi}p^{2}|_{2}^{\frac{1}{2}}+\beta|p^{2}(\xi=1)|_{2}^{\frac{1}{2}} ).
\]
Therefore, by Young's inequality we have
	\begin{align}
&\|Z_{3}\|_{3}(|\nabla_{e_\theta}\bu|_{2}+|\nabla_{e_\varphi}\bu|_{2}+|u|_{2} )|p^{2}|_{3}^{3/2}\notag\\
&\leq C\|Z_{3}\|_{3}(|\nabla_{e_\theta}\bu|_{2}+|\nabla_{e_\varphi}\bu|_{2}+|u|_{2} )|p^{2}|_{2}^{\frac{3}{4}}(|\nabla p^{2} |_{2}^{\frac{3}{4}}+|  \partial_{\xi}p^{2}|_{2}^{\frac{3}{4}}+\beta|p^{2}(\xi=1)|_{2}^{\frac{3}{4}} )\notag\\
&\leq \varepsilon (|\nabla p^{2} |_{2}^{2}+|  \partial_{\xi}p^{2}|_{2}^{2}+\beta|p^{2}(\xi=1)|_{2}^{2} )\notag\\
&\hspace{5mm}+C\|Z_{3}\|_{3}^{\frac{8}{5}}(|\nabla_{e_\theta}\bu|_{2}^{\frac{8}{5}}+|\nabla_{e_\varphi}\bu|_{2}^{\frac{8}{5}}+|u|_{2}^{\frac{8}{5}} )|p|_{4}^{\frac{12}{5}}.
    \end{align}
Similarly, we obtain
\begin{equation}
\label{eqn:5.13}
\|Z_{3}\|_{3}\|Z_{1}\|_{3}|p^{2}|_{3}^{3/2}
\leq \varepsilon (|\nabla p^{2} |_{2}^{2}+|  \partial_{\xi}p^{2}|_{2}^{2}+\beta|p^{2}(\xi=1)|_{2}^{2} )+C\|Z_{3}\|_{3}^{\frac{8}{5}}\|Z_{1}\|_{3}^{\frac{8}{5}}|p|_{4}^{\frac{12}{5}},
\end{equation}
and
\begin{align}
\label{eqn:5.14}
&\int_{\mathbf D}(Q_{q}+\gamma Z_{3})p^{3}d\mathbf{D}\notag\\
&\leq (| Q_{q}|_{2}+ \gamma|Z_{3}|_{2}  )|p^{2}|_{3}^{3/2}\notag\\
&\leq \varepsilon (|\nabla p^{2} |_{2}^{2}+|  \partial_{\xi}p^{2}|_{2}^{2}+\beta|p^{2}(\xi=1)|_{2}^{2})
+C(| Q_{q}|_{2}^{8/5}+ \gamma|Z_{3}|_{2}^{8/5}  )|p|_{4}^{12/5}.
\end{align}
From \eqref{eqn:L4ofp1}-\eqref{eqn:5.13}, we conclude that
\begin{align}
\label{eqn:5.15}
&\bigg|\int_{\mathbf D}(w(Z_1+\bu)\partial_\xi Z_3+\nabla_{Z_1+\bu} Z_3)p^3d\mathbf D\bigg|\notag\\
&\leq \varepsilon (|\nabla p^{2} |_{2}^{2}+|  \partial_{\xi}p^{2}|_{2}^{2}+\beta|p^{2}(\xi=1)|_{2}^{2} )\notag\\
&\hspace{5mm}+C\|Z_{3}\|_{3}^{\frac{8}{5}}(|\nabla_{e_\theta}\bu|_{2}^{\frac{8}{5}}+|\nabla_{e_\varphi}\bu|_{2}^{\frac{8}{5}}+|u|_{2}^{\frac{8}{5}} )|p|_{4}^{\frac{12}{5}}\notag\\
&\hspace{5mm}+ C\|Z_{3}\|_{3}^{\frac{8}{5}}\|Z_{1}\|_{3}^{\frac{8}{5}}|p|_{4}^{\frac{12}{5}}.
\end{align}
		
		Now back to the equation (\ref{eqn:L4ofp}), together with estimates \eqref{eqn:5.14}, \eqref{eqn:5.15}, we have
		\begin{align*}
			&\frac{d|p|_4^4}{dt}+(3-2\eps)(|\nabla p^2|_2^2+|\partial_\xi p^2|_2^2)+(4\beta-2\eps)\int_{\mathbf S^2} |p(\xi=1)|^4d\mathbf S^2\\
			&\leq C\|Z_{3}\|_{3}^{8/5}(|\nabla_{e_\theta}\bu|_2^{8/5}+|\nabla_{e_\varphi}\bu|_{2}^{8/5}+|\bu|_2^{8/5})|p|_{4}^{12/5}\\
&\hspace{5mm}+C(|Q_q|_2^{8/5}+|Z_3|_2^{8/5})|p|_4^{12/5}+C\norm{Z_1}_3^{8/5}\norm{Z_3}_3^{8/5}|p|_4^{12/5}.
 		\end{align*}
		By Young's inequality, we have
		\begin{align*}
			|p|_4^4=&\int_{\mathbf D}p^4d\mathbf D=-\int_{\mathbf S^2}\int_0^1\int_{\xi}^{1}\partial_\xi p^4+\int_{\mathbf S^2}\int_0^1p^4(\xi=1)\\
			\leq & 2|\partial_\xi p^2|_2^2+\frac{1}{2}|p|_4^4+\int_{\mathbf S^2}p^4(\xi=1)d\mathbf S^2,
		\end{align*}
		which implies that
		$$\frac{|p|_4^2}{2}\frac{d|p|_4^2}{dt}+|p|_4^4\leq C(|Q_q|_2^{8/5}+|Z_3|_2^{8/5}+\norm{Z_1}_3^{8/5}\norm{Z_3}_3^{8/5}+\norm{Z_3}_3^{8/5}\norm{\bu}_1^{8/5})|p|_4^{12/5}.$$
		Thus, we can apply Gronwall's inequality to
		$$\frac{d|p|_4^2}{dt}+|p|_4^2\leq C(|Q_q|_2^{8/5}+|Z_3|_2^{8/5}+\norm{Z_1}_3^{8/5}\norm{Z_3}_3^{8/5}+\norm{Z_3}_3^{8/5}\norm{\bu}_1^{8/5})|p|_4^{2/5},$$
		and get
		\begin{equation}
			\label{eqn:gronwallofL4p}
			|p|_4^2\leq |q_0|_4^2e^{-Ct}+C\int_{t_0}^te^{-C(t-s)}(|Q_q|_2^{8/5}+|Z_3|_2^{8/5}+\norm{Z_1}_3^{8/5}\norm{Z_3}_3^{8/5}+\norm{Z_3}_3^{8/5}\norm{\bu}_1^{8/5})ds.
		\end{equation}
		Taking inner product of Eq. (\ref{eqn:S}) with $S^3$ in $L^{2}(\mathbf D)$, we have
		\begin{align}
			\label{eqn:L4ofS}
			&\frac{1}{4}\frac{d|S|_4^4}{dt}+\frac{3}{4}|\nabla  S^2|_2^2+\frac{3}{4}|\partial_\xi \mathbf  S^2|_2^2+\alpha\int_{ \mathbf S^2}|S(\xi=1)|_4d\mathbf S^2\notag\\
			&=\int_{\mathbf D}(Q_T+\gamma Z_2)S^3d\mathbf D-\int_{\mathbf D}\Big(\nabla_{Z_1+\bu}(Z_2+S)+w(Z_1+\bu)\partial_\xi(Z_2+S)\Big)S^3d\mathbf D\notag\\
			&\hspace{5mm}+\int_{\mathbf D}\frac{br_s}{r}(1+a(Z_3+p))w(Z_1+\bu)S^3d\mathbf D.
		\end{align}
		Similarly as the evaluations in the previous subsection, by integration by parts, the H\"older inequality and the Sobolev imbedding theorem, one have that
		$$\int_{\mathbf D}\Big(\nabla_{Z_1+\bu}S+w(Z_1+\bu)\partial_\xi S\Big)S^3d\mathbf D=0.$$
Taking an  analogous argument of \eqref{eqn:5.14} and \eqref{eqn:5.15}, we have
\begin{align}
\label{eqn:5.18}
&\int_{\mathbf D}(Q_{T}+\gamma Z_{2})S^{3}d\mathbf{D}\notag\\
&\leq \varepsilon (|\nabla S^{2} |_{2}^{2}+|  \partial_{\xi}S^{2}|_{2}^{2}+\beta|S^{2}(\xi=1)|_{2}^{2})
+C(| Q_{T}|_{2}^{8/5}+ |Z_{2}|_{2}^{8/5}  )|S|_{4}^{12/5},
\end{align}
and
\begin{align}
&\bigg|\int_{\mathbf D}(w(Z_1+\bu)\partial_\xi Z_2+\nabla_{Z_1+\bu} Z_2)S^3d\mathbf D\bigg|\notag\\
&\leq \varepsilon (|\nabla S^{2} |_{2}^{2}+|  \partial_{\xi}S^{2}|_{2}^{2}+\beta|S^{2}(\xi=1)|_{2}^{2} )\notag\\
&\hspace{5mm}+C\|Z_{2}\|_{3}^{\frac{8}{5}}(|\nabla_{e_\theta}\bu|_{2}^{\frac{8}{5}}+|\nabla_{e_\varphi}\bu|_{2}^{\frac{8}{5}}+|u|_{2}^{\frac{8}{5}} )|S|_{4}^{\frac{12}{5}}\notag\\
&\hspace{5mm}+ C\|Z_{2}\|_{3}^{\frac{8}{5}}\|Z_{1}\|_{3}^{\frac{8}{5}}|S|_{4}^{\frac{12}{5}}.
\end{align}
		
		Now we evaluate the last term in Eq. (\ref{eqn:L4ofS}) separately by
		\begin{align}
			\label{eqn:L4ofS2}
			&\bigg|\int_{\mathbf D}\frac{br_s}{r}\bigg(\int_{\xi}^{1}\text{div }(Z_1+\bu)d\xi'\bigg)S^3d\mathbf D\bigg|\notag\\
			&\leq \bigg|\frac{br_s}{r}\int_{\xi}^{1}\text{div }(Z_1+\bu)d\xi'\bigg|_2|S^3|_2\notag\\
			&\leq C(\|Z_{1}\|_1+|\nabla_{e_\theta}\bu|_2+ |\nabla_{e_\varphi}\bu|_2)|S^{2}|_{2}^{3/4}(|\nabla S^{2}|_{2}^{3/4}+|\partial_{\xi} S^{2}|_{2}^{3/4} +|S^{2}|_{2}^{3/4} )\notag\\
&\leq \varepsilon(|\nabla S^{2}|_{2}^{2} +|\partial_{\xi}S^{2}|_{2}^{2} )+C(\|Z_{1}\|_{1}^{8/5}+|\nabla_{e_{\theta}} \bu|_{2}^{8/5}+ |\nabla_{e_{\varphi}} \bu|_{2}^{8/5}  )|S|_{4}^{12/5}\notag\\
&\hspace{5mm}+C( \|Z_{1}\|_{1}+|\nabla_{e_{\theta}} \bu|_{2}+ |\nabla_{e_{\varphi}} \bu|_{2} )|S|_{4}^{3}.
		\end{align}

By Lemma 2.3 and Young's inequality,  we get
\begin{align}
\label{eqn:5.21}
\bigg|\int_{\D}\frac{abr_{s}}{r}(p+Z_{3})\Big{(}\int_{\xi}^{1} \mathrm{div} \bu d\xi'   \Big{)}S^{3}  \bigg|
\leq \varepsilon |\nabla S^{2}|_{2}^{2}+C(|\nabla_{e_{\theta}}\bu |_{2}^{2}+ |\nabla_{e_{\varphi}}\bu |_{2}^{2}    )|p+Z_{3}|_{4}^{2}|S|_{4}^{2}.
\end{align}

		Now back to the equation (\ref{eqn:L4ofS}), together with estimates in \eqref{eqn:5.18}-\eqref{eqn:5.21}, we have
		\begin{align}
		\label{eqn:5.22}
			\frac{ |S|_4^2}{2}\frac{d|S|_4^2}{dt}&+(\frac{3}{4}-4\varepsilon)|\nabla  S^2|_2^2+(\frac{3}{4}-4\varepsilon)|\partial_\xi S^2|_2^2+(\alpha-2\varepsilon)\int_{ \mathbf S^2}|S(\xi=1)|_4d\mathbf S^2\notag\\
&\leq C(|Q_{T}|_{2}^{8/5}+|Z_{2}|_{2}^{8/5}+\|Z_{1}\|_{1}^{8/5}+ \|\bu\|_{1}^{8/5}\notag\\
&\ \ \ +\|Z_{2}\|_{3}^{8/5}\|\bu\|_{1}^{8/5}+\|Z_{1}\|_{3}^{8/5}\|Z_{2}\|_{3}^{8/5} )|S|_{4}^{12/5}\notag\\
&\ \ \ + C(\|Z_{1}\|_{1}+\|\bu\|_{1} )|S|_{4}^{3}+C\|\bu\|_{1}^{2}(|p|_{4}^{2}+|Z_{3}|_{4}^{2} )|S|_{4}^{2}.
		\end{align}
Since $\|S^{2}\|_{1}^{2}$ is equivalent to $|\nabla  S^2|_2^2+|\partial_\xi S^2|_2^2+ \alpha |S^{2}(\xi=1)|_2^{2}$, there exists a constant $C$ such that
$$|S|_{4}^{4}=|S^{2}|_{2}^{2}\leq C |\nabla  S^2|_2^2+|\partial_\xi S^2|_2^2+ \alpha |S^{2}(\xi=1)|_2^{2}.$$
Then by \eqref{eqn:5.22} we have
\begin{align}
\label{eqn:5.23}
\frac{d|S|_4^2}{dt}+C|S|_4^2
&\leq C(|Q_{T}|_{2}^{8/5}+|Z_{2}|_{2}^{8/5}+\|Z_{1}\|_{1}^{8/5}+ \|\bu\|_{1}^{8/5}\notag\\
&\ \ \ +\|Z_{2}\|_{3}^{8/5}\|\bu\|_{1}^{8/5}+\|Z_{1}\|_{3}^{8/5}\|Z_{2}\|_{3}^{8/5} )|S|_{4}^{2/5}\notag\\
&\ \ \ + C(\|Z_{1}\|_{1}+\|\bu\|_{1} )|S|_{4}+C\|\bu\|_{1}^{2}(|p|_{4}^{2}+|Z_{3}|_{4}^{2} )\notag\\
&\leq C(|Q_{T}|_{2}^{2}+|Z_{2}|_{2}^{2}+\|Z_{1}\|_{1}^{2}+\|Z_{1}\|_{3}^{2}\|Z_{2}\|_{3}^{2}+ \|\bu\|_{1}^{2}\notag\\
&\ \ \ +\|Z_{3}\|_{1}^{2}\|\bu\|_{1}^{2}+\|Z_{2}\|_{3}^{2}\|\bu\|_{1}^{2} +\|\bu\|_{1}^{2}|p|_{4}^{2}).
\end{align}

		Applying the Gronwall inequality, $\text{ for }t<\tau$,
		\begin{align}
		\label{eqn:5.24}
			|S(t)|_4^2 & \leq |S_0|_4^2e^{-Ct}+C\int_{t_0}^te^{-C(t-s)}(|Q_{T}|_{2}^{2}+|Z_{2}|_{2}^{2}+\|Z_{1}\|_{1}^{2}+\|Z_{1}\|_{3}^{2}\|Z_{2}\|_{3}^{2})ds\notag\\
& \hspace{5mm}+C\int_{t_0}^te^{-C(t-s)}(1+\|Z_{3}\|_{1}^{2}+ \|Z_{2}\|_{3}^{2}+|p|_{4}^{2})\|\bu\|_{1}^{2}ds.
		\end{align}
		\subsection{$L^4$ estimates of $\bu$.}
		To estimate $L^4$ norm of $\bu$, we denote by
		$$\bar{\bu}(\theta,\varphi)=\int_0^1\bu(\theta,\varphi,\xi)d\xi,\text{ for }(\theta,\varphi)\in \mathbf S^2,$$
		and
		$\wt{\bu}=\bu-\bar{\bu}$. Note that
		$$\bar{\wt{\bu}}=0\text{ and div }\bar{\bu}=0.$$
		Now we take the average value of Eq. (\ref{eqn:u}) with respect to $\xi$
		\begin{align*}
			&\partial_t\bar{\bu}+\overline{\nabla_{Z_1+\bu}(Z_1+\bu)+w(Z_1+\bu)\partial_\xi(Z_1+\bu)}+\frac{f}{R_0}\vec k\times(\bar Z_1+\bar{\bu})+\nabla\Phi_s\\
			&=\Delta\bar{\bu}+\gamma\bar Z_1-\int_0^1\int_{\xi}^{1}\frac{br_s}{r}\nabla[1+a(Z_3+p)(Z_2+S)]d\xi'd\xi.
		\end{align*}
		Now applying integration by parts and boundary conditions,
		$$\int_0^1\nabla_{Z_1+\bu}(Z_1+\bu)d\xi=\int_0^1\nabla_{\wt Z_1+\wt \bu}(\wt Z_1+\wt \bu)d\xi+\nabla_{\bar Z_1+\bar{\bu}}(\bar Z_1+\bar{\bu}),$$
		and
		$$\int_0^1w(Z_1+\bu)\partial_\xi(Z_1+\bu)d\xi=\int_0^1(Z_1+\bu)\text{div }(Z_1+\bu)d\xi=\int_0^1(\wt Z_1+\wt\bu)\text{div }(\wt Z_1+\wt\bu)d\xi.$$
		Thus, we have $\bar{\bu}$ satisfy the following equation and boundary conditions
		\begin{subequations}
			\label{eqn:baru}
			\begin{align}
				&\partial_t\bar{\bu}+\overline{\nabla_{\wt Z_1+\wt \bu}(\wt Z_1+\wt \bu)+(\wt Z_1+\wt\bu)\text{div }(\wt Z_1+\wt\bu)}+\nabla_{\bar Z_1+\bar{\bu}}(\bar Z_1+\bar{\bu})\notag\\
				&+\frac{f}{R_0}\vec k\times(\bar Z_1+\bar{\bu})+\nabla\Phi_s=\Delta\bar{\bu}+\gamma\bar Z_1-\int_0^1\int_{\xi}^{1}\frac{br_s}{r}\nabla[1+a(Z_3+p)(Z_2+S)]d\xi'd\xi,\\
				&\text{div }\bar{\bu}=0,\text{ on }\mathbf S^2.
			\end{align}
		\end{subequations}
		By subtracting Eq. (\ref{eqn:baru}) from Eq. (\ref{eqn:u}), we can also conclude that $\wt\bu$ satisfies the following equation and boundary conditions
		\begin{subequations}
			\label{eqn:wtu}
			\begin{align}
				&\partial_t\wt\bu+\nabla_{\wt Z_1+\wt\bu}(\wt Z_1+\wt\bu)+w(\wt Z_1+\wt\bu)\partial_\xi(\wt Z_1+\wt\bu)+\nabla_{\wt Z_1+\wt\bu}(\bar Z_1+\bar{\bu})+\nabla_{\bar Z_1+\bar{\bu}}(\wt Z_1+\wt\bu)-\notag\\
				&\overline{\nabla_{\wt Z_1+\wt \bu}(\wt Z_1+\wt \bu)+(\wt Z_1+\wt\bu)\text{div }(\wt Z_1+\wt\bu)}+\frac{f}{R_0}\vec k\times(\wt Z_1+\wt\bu)-\Delta\wt\bu-\partial_{\xi\xi}\wt\bu\notag\\
				&=\int_\xi^1\frac{br_s}{r}\nabla[1+a(Z_3+p)(Z_2+S)]d\xi'-\int_0^1\int_{\xi}^{1}\frac{br_s}{r}\nabla[1+a(Z_3+p)(Z_2+S)]d\xi'd\xi,\\
				&\partial_\xi\wt\bu=0,\text{ when }\xi=0\text{ and }\xi=1.
			\end{align}
		\end{subequations}

By the definition of covariant derivative, for $h \in C^{\infty}(\mathbf{D}) $ and $u=(u_{\theta}, u_{\varphi} )\in C^{\infty}(T \mathbf{D}|T \mathbf{S}^{2} ),$ we have
\begin{align}	
 &\nabla_{e_{\theta}}(h\wt \bu  )=h\nabla_{e_{\theta}}\wt \bu+ \wt \bu\nabla_{e_{\theta}}h,\\	
 &\nabla_{e_{\varphi}}(h\wt \bu  )=h\nabla_{e_{\varphi}}\wt \bu+ \wt \bu\nabla_{e_{\varphi}}h,\\	
 &\nabla_{e_{\theta}}(u \cdot\wt \bu  )=u \cdot\nabla_{e_{\theta}}\wt \bu+ \wt \bu \cdot\nabla_{e_{\theta}}u,\\	
 &\nabla_{e_{\varphi}}(u \cdot\wt \bu  )=u \cdot\nabla_{e_{\varphi}}\wt \bu+ \wt \bu \cdot\nabla_{e_{\varphi}}u,\\	
 &\nabla_{u}\wt \bu =u_{\theta}\nabla_{e_{\theta}}\wt \bu+ u_{\varphi}\nabla_{e_{\varphi}}\wt \bu.
\end{align}
		Applying integration by parts, we have
		$$\int_{\mathbf D}\bigg[\nabla_{\wt\bu}\wt\bu+\bigg(\int_{\xi}^{1}\text{div }\wt\bu d\xi'\bigg)\partial_\xi\wt\bu\bigg]\cdot|\wt\bu|^2 \wt\bu d\mathbf D=0.$$
		 Using integration by parts together with $\text{div }\bar{\bu}=0$,
		$$\int_{\mathbf D}\nabla_{\bar Z_1+\bar{\bu}}\wt\bu\cdot(|\wt\bu|^2\wt\bu) d\mathbf D=-\frac{1}{4}\int_{\mathbf D}|\wt\bu|^4\text{div }(\bar Z_1+\bar{\bu})d\mathbf D=0.$$
		Similarly, we will also have
	\begin{align*}		
&\int_{\mathbf D}\overline{\nabla_{\wt Z_1+\wt \bu}(\wt Z_1+\wt \bu)+(\wt Z_1+\wt\bu)\text{div }(\wt Z_1+\wt\bu)}\cdot|\wt\bu|^2 \wt\bu d\mathbf D\\
			&=-\int_{\mathbf D}\overline{(\wt Z_{1,\theta}+\wt u_\theta)(\wt Z_1+\wt\bu)}\cdot\nabla_{e_\theta}(|\wt\bu|^2 \wt\bu) d\mathbf D-\int_{\mathbf D}\overline{(\wt Z_{1,\varphi}+\wt u_\varphi)(\wt Z_1+\wt\bu)}\cdot\nabla_{e_\varphi}(|\wt\bu|^2 \wt\bu )d\mathbf D.
		\end{align*}

		Now taking inner product of the Eq. (\ref{eqn:wtu}) with $|\wt\bu|^2 \wt\bu$ in $(L^2(\mathbf D))^2$, and using the above equalities, we get
		\begin{align}
			&\frac{1}{4}\frac{d|\wt\bu|_4^4}{dt}+\frac{1}{2}\int_{\mathbf D}\Big(|\nabla_{e_\theta}|\wt\bu^2||^2+|\nabla_{e_\varphi}|\wt\bu^2||^2+|\partial_\xi|\wt\bu^2||^2\Big)d\mathbf D\notag\\
			&\hspace{5mm}+\int_{\mathbf D}|\wt\bu|^2\Big(|\nabla_{e_\theta}\wt\bu|^2+|\nabla_{e_\varphi}\wt\bu|^2+|\partial_\xi\wt\bu|^2+|\wt\bu|^2\Big)d\mathbf D\notag\\
			&=-\int_{\mathbf D}[\nabla_{\wt Z_1}\wt\bu+\nabla_{\wt\bu}\wt Z_1+\nabla_{\wt Z_1}\wt Z_1]\cdot|\wt\bu|^2 \wt\bu d\mathbf D\notag\\
			&\hspace{5mm}-\int_{\mathbf D}[w(\wt Z_1)\partial_\xi\wt\bu+w(\wt\bu)\partial_\xi\wt Z_1+w(\wt Z_1)\partial_\xi\wt Z_1]\cdot|\wt\bu|^2 \wt\bu d\mathbf D\notag\\
			&\hspace{5mm}+\int_{\mathbf D}\nabla_{\wt Z_1+\wt\bu}  (\bar Z_1+\bar\bu)  \cdot |\wt\bu|^2 \wt\bu  d\mathbf D\notag\\
			&\hspace{5mm}-\int_{\mathbf D}\nabla_{\bar Z_1+\bar{\bu}}\wt Z_1\wt\bu^3 d\mathbf D-\int_{\mathbf D}\bigg(\frac{f}{R_0}\vec k\times\wt Z_1\bigg)\cdot|\wt\bu|^2 \wt\bu d\mathbf D\notag\\
			&\hspace{5mm}-\int_{\mathbf D}\overline{(\wt Z_{1,\theta}+\wt u_\theta)(\wt Z_1+\wt\bu)}\cdot\nabla_{e_\theta}(|\wt\bu|^2\wt\bu)d\mathbf D-\int_{\mathbf D}\overline{(\wt Z_{1,\varphi}+\wt u_\varphi)(\wt Z_1+\wt\bu)}\cdot\nabla_{e_\varphi}(|\wt\bu|^2\wt\bu) d\mathbf D\notag\\
			&\hspace{5mm}+\int_{\mathbf D}\bigg(\int_\xi^1\frac{br_s}{r}\nabla[1+a(Z_3+p)](Z_2+S)d\xi'\bigg)\cdot|\wt\bu|^2 \wt\bu d\mathbf D\notag\\
			&\hspace{5mm}-\int_{\mathbf D}\bigg(\int_0^1\int_{\xi}^{1}\frac{br_s}{r}\nabla[1+a(Z_3+p)](Z_2+S)d\xi'd\xi\bigg)\cdot|\wt\bu|^2 \wt\bu d\mathbf D:=\sum_{i=1}^{9}I_i.
		\end{align}
		We will estimate $I_i$ respectively for $i=1,\cdots,9$. Now applying integration by parts to first terms of $I_1$ and $I_2$,
		\begin{align*}
			&\int_{\mathbf D}\nabla_{\wt Z_1}\wt\bu\cdot|\wt\bu|^2 \wt\bu d\mathbf D+\int_{\mathbf D}w(\wt Z_1)\partial_\xi\wt\bu\cdot|\wt\bu|^2 \wt\bu d\mathbf D\\
			&=-\frac{1}{4}\int_{\mathbf D}\nabla_{\wt Z_1}|\wt\bu|^4d\mathbf D-\frac{1}{4}\int_{\mathbf D}w(\wt Z_1)\partial_\xi|\wt\bu|^4=0.
		\end{align*}
		Then applying the H\"older inequality and the interpolation inequality to the other terms in $I_1$, we get
		\begin{align*}
			|I_1|\leq&\bigg|\int_{\mathbf D}\nabla_{\wt\bu}\wt Z_1\cdot(|\wt\bu|^2 \wt\bu)d\mathbf D\bigg|+\bigg|\int_{\mathbf D}\nabla_{\wt Z_1}\wt Z_1\cdot(|\wt\bu|^2 \wt\bu )d\mathbf D\bigg|\\
			\leq&C\bigg(\int_{\mathbf D}\Big(|\nabla_{e_\theta}\wt Z_1|^2+|\nabla_{e_\varphi}\wt Z_1|^2\Big)^{3/2}d\mathbf D\bigg)^{1/3}\bigg\{\bigg(\int_{\mathbf D}|\wt\bu|^{2\times 3}d\mathbf D\bigg)^{2/3}+|\wt Z_1|_6\bigg(\int_{\mathbf D}|\wt\bu|^{2\times 3}d\mathbf D\bigg)^{1/2}\bigg\}\\
			\leq&C\bigg(\int_{\mathbf D}\Big(|\nabla_{e_\theta}\wt Z_1|^2+|\nabla_{e_\varphi}\wt Z_1|^2\Big)^{3/2}d\mathbf D\bigg)^{1/3}\bigg\{||\wt\bu|^2|_2\norm{|\wt\bu|^2}_1+||\wt\bu|^2|_2^{3/4}\norm{|\wt\bu|^2}_1^{3/4}|\wt Z_1|_6\bigg\}\\
			\leq&\eps\int_{\mathbf D}\Big(|\nabla_{e_\theta}|\wt\bu|^2|^2+|\nabla_{e_\varphi}|\wt\bu|^2|^2+|\partial_\xi|\wt\bu|^2|^2+|\wt\bu|^{4}\Big)d\mathbf D+C\norm{ Z_1}_2^2|\wt\bu|_4^4+C\norm{ Z_1}_2^{16/5}|\wt\bu|_4^{12/5}.
		\end{align*}
		Applying integration by parts and the interpolation inequality on $\mathbf S^2$ to the other terms of $I_2$,
\begin{align*}
-\int_{\mathbf D}w(\wt\bu)\partial_\xi\wt Z_1\cdot(|\wt\bu|^{2} \wt\bu)d\mathbf D&\leq |\partial_\xi Z_1 |_{\infty}\int_{\mathbf{S}^{2}}(\int_{0}^{1}|\mathrm{div} \wt\bu|d\xi
 \int_{0}^{1} |\wt\bu|^{3}d\xi)d\mathbf{S}^{2}\\
 &\leq C\|Z_{1}\|_{3}(|\nabla_{e_{\theta}}\wt\bu |_{2}+ |\nabla_{e_{\varphi}}\wt\bu |_{2} )\int_{0}^{1}\Big{(}\int_{\mathbf{S}^{2}}|\wt\bu|^{3\times 2}  d \mathbf{S}^{2} \Big{)}^{1/2}d\xi\\
  &\leq C\|Z_{1}\|_{3}\|\bu\|_{1}\int_{0}^{1}||\wt\bu|^{2} |_{L^{3}(\mathbf{S}^{2} )}^{3/2}d\xi\\
  &\leq C\|Z_{1}\|_{3}\|\bu\|_{1}\int_{0}^{1}||\wt\bu|^{2}|_{L^{2}(\mathbf{S}^{2} )}(|\nabla_{e_{\theta}}|\wt\bu|^{2}|_{L^{2}(\mathbf{S}^{2} )}^{1/2}+
 |\nabla_{e_{\varphi}}|\wt\bu|^{2}|_{L^{2}(\mathbf{S}^{2} )}^{1/2}+||\wt\bu|^{2}|_{L^{2}(\mathbf{S}^{2} )}^{1/2}  )d\xi\\
 &\leq \varepsilon ( |\nabla_{e_{\theta}}|\wt\bu|^{2}|_{2}^{2}+|\nabla_{e_{\varphi}}|\wt\bu|^{2}|_{2}^{2}  )\\
 &\hspace{5mm}+C\|Z_{1}\|_{3}^{4/3}\|\bu\|_{1}^{4/3}|\wt\bu|_{4}^{8/3}\\
 &\hspace{5mm}+C\|Z_{1}\|_{3}\|\bu\|_{1}|\wt\bu|_{4}^{3}.
\end{align*}

Using the H\"older inequality and the Sobolev embedding theorem yields that
\begin{align*}
-\int_{\mathbf{D}}(\nabla_{\wt Z_{1}}\wt Z_{1})|\wt \bu|^{2}\wt \bu d \mathbf{D}\leq |Z_{1}|_{\infty}( |\nabla_{e_{\theta}}Z_{1}|_{\infty}+ |\nabla_{e_{\varphi}}Z_{1}|_{\infty} )|\wt \bu|_{4}^{3}\leq C\|Z_{1}\|_{3}^{2} |\wt \bu|_{4}^{3}.
\end{align*}
By the argument above, we have estimates for $I_2$ as
		$$|I_2|\leq \eps\Big(|\nabla_{e_\theta}|\wt\bu|^2|_2^{2}+|\nabla_{e_\varphi}|\wt\bu|^2|_2^{2}\Big)+C\norm{ Z_1}_3^2|\wt\bu|_4^3+C\norm{ Z_1}_3^{4/3}\|\bu\|_1^{4/3}|\wt\bu|_{4}^{8/3}+C\|Z_{1}\|_{3}\|\bu\|_{1}|\wt\bu|_{4}^{3}.$$
		By the interpolation inequality, the H\"older inequality and the Minkowski inequality,
\begin{align*}
&\int_{\mathbf D}(\nabla_{\wt\bu}   \bar{\bu}    )  \cdot |\wt \bu|^{2}\wt \bu d \mathbf D\\
&\leq \int_{\mathbf{S}^{2}}\Big{(}( |\nabla_{e_{\theta}} \bar{\bu}|+ |\nabla_{e_{\varphi}} \bar{\bu}|  )\int_{0}^{1}|\wt \bu|^{4}d\xi\Big{)} d \mathbf{S}^{2}\\
&\leq (|\nabla_{e_{\theta}} \bar{\bu}|_{2}+|\nabla_{e_{\varphi}} \bar{\bu}|_{2} )\int_{0}^{1}\Big{(} \int_{\mathbf{S}^{2}}|\wt \bu|^{2\times 4}d \mathbf{S}^{2}   \Big{)}^{1/2}d\xi\\
&\leq (|\nabla_{e_{\theta}} \bar{\bu}|_{2}+|\nabla_{e_{\varphi}} \bar{\bu}|_{2} )\int_{0}^{1}\Big{(} ||\wt \bu|^{2}|_{L^{2}(\mathbf{S}^{2})}
(|\nabla_{e_{\theta}}   |\wt \bu|^{2}|_{L^{2}(\mathbf{S}^{2})} +|\nabla_{e_{\varphi}}  | \wt \bu|^{2}|_{L^{2}(\mathbf{S}^{2})} +||\wt \bu|^{2}|_{L^{2}(\mathbf{S}^{2})} )\Big{)} d\xi\\
&\leq \varepsilon (|\nabla_{e_{\theta}}|\wt \bu|^{2} |_{2}^{2}+ |\nabla_{e_{\varphi}}|\wt \bu|^{2} |_{2}^{2}  )+C(1+|\nabla_{e_{\theta}}\bar{\bu} |_{2}^{2}+|\nabla_{e_{\varphi}}\bar{\bu} |_{2}^{2} )|\wt \bu|_{4}^{4}.
\end{align*}
Similarly, we have
\begin{align*}
&\int_{\mathbf D}(\nabla_{\wt Z_{1}} \bar Z_{1})\cdot (|\wt\bu|^{2}\wt \bu) d \mathbf{D}+\int_{\mathbf D}(\nabla_{\wt Z_{1}} \bar \bu)\cdot (|\wt\bu|^{2}\wt \bu) d \mathbf{D}
+\int_{\mathbf D}(\nabla_{\wt \bu} \bar Z_{1})\cdot (|\wt\bu|^{2}\wt \bu) d \mathbf{D}\\
&\leq |Z_{1}|_{\infty}(|\nabla_{e_{\theta}} Z_{1} |_{4}+|\nabla_{e_{\varphi}} Z_{1} |_{4} )|\wt \bu|_{4}^{3}\\
&\hspace{5mm}+C|Z_{1}|_{\infty}( |\nabla_{e_{\theta}} \bar \bu |_{2}+|\nabla_{e_{\varphi}} \bar \bu |_{2} )||\wt \bu|^{3}|_{2}+|\nabla_{e_{\theta}} Z_{1}+\nabla_{e_{\varphi}} Z_{1} |_{\infty}|\wt \bu|_{4}^{4}\\
&\leq C\|Z_{1}\|_{2}^{2}|\wt \bu|_{4}^{3}+C\|Z_{1}\|_{2}(|\nabla_{e_{\theta}} \bar \bu |_{2}+|\nabla_{e_{\varphi}} \bar \bu |_{2} )||\wt \bu|^{2}|_{3}^{3/2}
+C\|Z_{1}\|_{3}|\wt \bu|_{4}^{4}\\
&\leq C\|Z_{1}\|_{2}^{2}|\wt \bu|_{4}^{3}+C\|Z_{1}\|_{3}|\wt \bu|_{4}^{4}\\
&\hspace{5mm}+C\|Z_{1}\|_{2}\|\bu\|_{1}||\wt \bu|^{2}|_{2}^{3/4}(| \nabla_{e_{\theta}}|\wt\bu|^{2}|_{2}^{3/4} + | \nabla_{e_{\varphi}}|\wt\bu|^{2}|_{2}^{3/4}
+|\partial_{\xi} |\wt\bu|^{2}|_{2}^{3/4}+||\wt\bu|^{2}|_{2}^{3/4} )\\
&\leq \varepsilon (| \nabla_{e_{\theta}}|\wt\bu|^{2}|_{2}^{2} + | \nabla_{e_{\varphi}}|\wt\bu|^{2}|_{2}^{2}
+|\partial_{\xi} |\wt\bu|^{2}|_{2}^{2} )+C\|Z_{1}\|_{2}^{8/5}\|\bu\|_{1}^{8/5}|\wt\bu|_{4}^{12/5}\\
&\hspace{5mm}+C\|Z_{1}\|_{2}\|\bu\|_{1}|\wt\bu|_{4}^{3}+C\|Z_{1}\|_{2}^{2}|\wt\bu|_{4}^{3}+C\|Z_{1}\|_{3}|\wt \bu|_{4}^{4}.
\end{align*}
Therefore, combining the above estimates we obtain
\begin{align*}
I_{3}&\leq \varepsilon (| \nabla_{e_{\theta}}|\wt\bu|^{2}|_{2}^{2} + | \nabla_{e_{\varphi}}|\wt\bu|^{2}|_{2}^{2}+|\partial_{\xi} |\wt\bu|^{2}|_{2}^{2} )+C\|Z_{1}\|_{2}^{8/5}\|\bu\|_{1}^{8/5}|\wt\bu|_{4}^{12/5}\\
&\hspace{5mm}+C\|Z_{1}\|_{2}\|\bu\|_{1}|\wt\bu|_{4}^{3}+C\|Z_{1}\|_{2}^{2}|\wt\bu|_{4}^{3}+C(1+\|Z_{1}\|_{2}+\|\bu\|_{1}^{2})|\wt \bu|_{4}^{4}.
\end{align*}

		Analogously, we have
		\begin{align*}
			|I_4|=&\bigg|\int_{\mathbf D}\nabla_{\bar Z_1+\bar{\bu}}\wt Z_1\cdot(|\wt\bu|^2 \wt\bu )d\mathbf D\bigg|\\
			\leq&||\wt\bu|^{3}|_2\bigg(\int_{\mathbf D}\Big(|\nabla_{e_\theta}\wt Z_1|^4+|\nabla_{e_\varphi}\wt Z_1|^4\Big)d \mathbf{D}\bigg)^{1/4}(|\bar Z_1|_{L^4(\mathbf S^2)}+|\bar{\bu}|_{L^4(\mathbf S^2)})\\
\leq&||\wt\bu|^{2}|_3^{3/2}\bigg(\int_{\mathbf D}\Big(|\nabla_{e_\theta}\wt Z_1|^4+|\nabla_{e_\varphi}\wt Z_1|^4\Big)d \mathbf{D}\bigg)^{1/4}(|\bar Z_1|_{L^4(\mathbf S^2)}+|\bar{\bu}|_{L^4(\mathbf S^2)})\\
			\leq&C||\wt\bu|^2|_2^{3/4}\norm{|\wt\bu|^2}_1^{3/4}\norm{\wt Z_1}_2(|\bar Z_1|_{L^4(\mathbf S^2)}+|\bar{\bu}|_{L^4(\mathbf S^2)})\\
			\leq&\eps\int_\mathbf D \Big(|\nabla_{e_\theta}|\wt\bu|^2|^2+|\nabla_{e_\varphi}|\wt\bu|^2|^2+|\partial_\xi |\wt\bu|^2|^2\Big)d\mathbf D\\
			&+C\norm{\wt Z_1}_2(|\bar{\bu}|_4+|\bar Z_1|_4)|\wt\bu|_4^3+C\norm{\wt Z_1}_2^{8/5}(|\bar{\bu}|_{L^4(\mathbf S^2)}^{8/5}+|\bar Z_1|_{L^4(\mathbf S^2)}^{8/5})|\wt\bu|_4^{12/5}.
		\end{align*}
		By the H\"older inequality and the Sobolev embedding theorem, we have
		$$|I_5|\leq C\norm{ Z_1}_1|\wt\bu|_4^3.$$
To estimate $I_{5},$ we first consider
\begin{align}
\label{eqn:5.33}
&\hspace{5mm}\int_{\mathbf{D} }\overline{\wt u_{\theta} \wt \bu}\cdot(\wt \bu \nabla_{e_{\theta}}|\wt \bu|^{2}  )\nonumber\\
&\leq \int_{\mathbf{S}^{2}}\overline{|\wt\bu|^{2}}\Big{(} \int_{0}^{1}|\nabla_{e_{\theta}}|\wt \bu|^{2} |^{2}d\xi  \Big{)}^{1/2}\Big{(} \int_{0}^{1}|\wt \bu|^{2}d\xi\Big{)}^{1/2}d\mathbf{S}^{2}\nonumber\\
&\leq \int_{\mathbf{S}^{2}}\Big{(} \int_{0}^{1}|\wt \bu|^{2} d\xi   \Big{)}^{3/2}\Big{(} \int_{0}^{1}|\nabla_{e_{\theta}}|\wt \bu|^{2} |^{2}d\xi  \Big{)}^{1/2}d\mathbf{S}^{2}\nonumber\\
&\leq |\nabla_{e_{\theta}}|\wt \bu|^{2} |_{2}\Big{(}\int_{0}^{1}  \Big{(}\int_{\mathbf{S}^{2}}|\wt \bu|^{6} d\mathbf{S}^{2}\Big{)}^{1/3}d\xi\Big{)}^{3/2}\nonumber\\
&\leq \varepsilon|\nabla_{e_{\theta}}|\wt \bu|^{2} |_{2}^{2}+C\Big{(}\int_{0}^{1}\Big{(}\int_{\mathbf{S}^{2}}|\wt \bu|^{6} d\mathbf{S}^{2}     \Big{)}^{1/3}d\xi    \Big{)}^{3}\nonumber\\
&\leq \varepsilon|\nabla_{e_{\theta}}|\wt \bu|^{2} |_{2}^{2}+C\Big{(}\int_{0}^{1}|\wt \bu|_{L^{4}(\mathbf{S}^{2})}^{4/3}( |\nabla_{e_{\theta}}\wt \bu|_{L^{2}(\mathbf{S}^{2})}^{2/3}+ |\nabla_{e_{\varphi}}\wt \bu|_{L^{2}(\mathbf{S}^{2})}^{2/3}+|\wt \bu|_{2}^{2/3}  )d\xi\Big{)}^{3}\nonumber\\
&\leq \varepsilon|\nabla_{e_{\theta}}|\wt \bu|^{2} |_{2}^{2}+C|\wt \bu|_{4}^{4}\|\bu \|_{1}^{2}.
\end{align}
Then
\begin{align}
&\hspace{5mm}-\int_{\mathbf{D}}	\overline{\wt Z_{1,\theta} \wt Z_{1}}\nabla_{e_{\theta}}(|\wt \bu|^{2}\wt \bu )d\mathbf{D}\nonumber\\
&=-\int_{\mathbf{D}}	\overline{\wt Z_{1,\theta} \wt Z_{1}}(|\wt \bu|^{2}\nabla_{e_{\theta}}\wt \bu+ \wt \bu \nabla_{e_{\theta}}|\wt \bu|^{2})d\mathbf{D} \nonumber\\
&\leq |Z_{1}|_{\infty}^{2}|\bu|_{2}   \Big{(}\int_{\mathbf{D}}(|\wt \bu|^{2}| \nabla_{e_{\theta}}\wt \bu|^{2} +  |\nabla_{e_{\theta}}|\wt \bu|^{2} |^{2})d \mathbf{D}\Big{)}^{1/2}\nonumber\\
&\leq \varepsilon \Big{(}\int_{\mathbf{D}}(|\wt \bu|^{2}| \nabla_{e_{\theta}}\wt \bu|^{2} +  |\nabla_{e_{\theta}}|\wt \bu|^{2} |^{2})d \mathbf{D}\Big{)}+C\|Z_{1}\|_{2}^{4}|\bu|_{2}^{2}.
\end{align}
Using the H\"older inequality and the Minkowski inequality, the interpolation inequality and Young's inequality we obtain
\begin{align}
&\hspace{5mm}-\int_{\mathbf{D}}	\overline{ \wt Z_{1,\theta} \wt \bu}\nabla_{e_{\theta}}(|\wt \bu|^{2}\wt \bu )d\mathbf{D}\nonumber\\
&\leq\int_{\mathbf{S}^{2}}\int_{0}^{1}|\wt Z_{1,\theta} \wt \bu |d\xi\int_{0}^{1}(|\wt\bu|^{2}|\nabla_{e\theta} \wt\bu|+ |\wt\bu||\nabla_{e_{\theta}}|\wt\bu|^{2} |  )d\mathbf{D}\nonumber\\
&\leq |Z_{1}|_{\infty}\int_{\mathbf{S}^{2}}\Big{[}\int_{0}^{1}|\wt \bu |d\xi\Big{(}\int_{0}^{1}|\wt\bu|^{2}d\xi\Big{)}^{1/2}\Big{(}\int_{0}^{1}|\wt\bu|^{2} |\nabla_{e_{\theta}}\wt\bu|^{2} d\xi\Big{)}^{1/2}\Big{]}d\mathbf{S}^{2}\nonumber\\
&\hspace{5mm}+ |Z_{1}|_{\infty}\int_{\mathbf{S}^{2}}\Big{[}\int_{0}^{1}|\wt \bu |d\xi\Big{(}\int_{0}^{1}|\wt\bu|^{2}d\xi\Big{)}^{1/2}\Big{(}\int_{0}^{1} |\nabla_{e_{\theta}}|\wt\bu|^{2}|^{2} d\xi\Big{)}^{1/2}\Big{]}d\mathbf{S}^{2}\nonumber\\
&\leq C\|Z_{1}\|_{2}|\wt\bu|_{3}^{3/2}\Big{(}\int_{\mathbf{D}}(|\wt \bu|^{2}||\nabla_{e_{\theta}}\wt\bu|^{2}+ |\nabla_{e_{\theta}}|\wt\bu|^{2}|^{2}  )d\mathbf{D}\Big{)}^{1/2}\nonumber\\
&\leq \varepsilon( \int_{\mathbf{D}}(|\wt \bu|^{2}||\nabla_{e_{\theta}}\wt\bu|^{2}+ |\nabla_{e_{\theta}}|\wt\bu|^{2}|^{2}  )d\mathbf{D}  )+C\|Z_{1}\|_{2}^{2}|\wt\bu|_{4}^{3}.
\end{align}
Analogously, we have
\begin{equation}
\label{eqn:5.36}
-\int_{\mathbf{D}}	\overline{\wt u_{\theta} \wt Z_{1} }\nabla_{e_{\theta}}(|\wt \bu|^{2}\wt \bu )d\mathbf{D}\leq \varepsilon( \int_{\mathbf{D}}(|\wt \bu|^{2}||\nabla_{e_{\theta}}\wt\bu|^{2}+ |\nabla_{e_{\theta}}|\wt\bu|^{2}|^{2}  )d\mathbf{D}  )+C\|Z_{1}\|_{2}^{2}|\wt\bu|_{4}^{3}.
\end{equation}
Combining the estimates \eqref{eqn:5.33}-\eqref{eqn:5.36}  yields that
\begin{align}
\label{eqn:5.37}
I_{6}=&-\int_{\mathbf D}\overline{(\wt Z_{1,\theta}+\wt u_\theta)(\wt Z_1+\wt\bu)}\cdot\nabla_{e_\theta}(|\wt\bu|^2\wt\bu)d\mathbf D\nonumber\\
\leq& \varepsilon(|\nabla_{e_{\theta}}|\wt\bu|^{2} |_{2}^{2} +||\wt\bu| | \nabla_{e_{\theta}}\wt\bu||_{2}^{2} )+C\|Z_{1}\|_{2}^{4}|\bu|_{2}^{2}+
C\|Z_{1}\|_{2}^{2}|\wt\bu|_{4}^{3}+C\|\bu\|_{1}^{2}|\wt\bu|_{4}^{4}.
\end{align}
Repeating the argument of \eqref{eqn:5.37} we get
\begin{align}
I_{7}=&-\int_{\mathbf D}\overline{(\wt Z_{1,\varphi}+\wt u_\varphi)(\wt Z_1+\wt\bu)}\cdot\nabla_{e_\varphi}(|\wt\bu|^2\wt\bu) d\mathbf D\notag\\
\leq& \varepsilon(|\nabla_{e_{\varphi}}|\wt\bu|^{2} |_{2}^{2} +||\wt\bu| | \nabla_{e_{\varphi}}\wt\bu||_{2}^{2} )+C\|Z_{1}\|_{2}^{4}|\bu|_{2}^{2}+
C\|Z_{1}\|_{2}^{2}|\wt\bu|_{4}^{3}+C\|\bu\|_{1}^{2}|\wt\bu|_{4}^{4}.
\end{align}
By the definition of the horizontal covariant derivative and  the  horizontal divergence  (see \eqref{eqn:2.2a} and \eqref{eqn:2.3a}), we have
\begin{align*}
\mathrm{div }(|\wt\bu|^{2}\wt\bu  )=(\nabla_{e_{\theta}} |\wt \bu|^{2} +\nabla_{e_{\varphi}} |\wt \bu|^{2} )\cdot \wt \bu+  |\wt \bu|^{2} \mathrm{div} \wt \bu,
 \end{align*}
and
\begin{align*}
|\mathrm{div} \wt\bu |\leq |\nabla_{e_{\theta}}\wt\bu|+| \nabla_{e_{\varphi}}\wt\bu|.
 \end{align*}
Therefore applying integration by parts to $I_8$ and $I_{9}$,we obtain
		\begin{align*}
			|I_8+I_{9}|\leq&C\bigg|\int_{\mathbf D}\bigg(\int_{\xi}^{1}[1+a(Z_3+p)](Z_2+S)d\xi'\cdot\text{div }(|\wt\bu|^2 \wt\bu) \bigg)d\mathbf D\\
			&+\int_{\mathbf D}\bigg(\int_0^1\int_{\xi}^{1}[1+a(Z_3+p)](Z_2+S)d\xi'd\xi\cdot\text{div }(|\wt\bu|^2 \wt\bu) \bigg)d\mathbf D\bigg|\\
			\leq&C\int_{\mathbf S^2}(|\overline{Z_2+S}|+|\overline{(Z_3+p)(Z_2+S)}|)\int_0^1|\wt\bu|^2\Big(|\nabla_{e_\theta}\wt\bu|^2+|\nabla_{e_\varphi}\wt\bu|^2\Big)^{1/2}d\xi d\mathbf S^2\\
&+C\int_{\mathbf S^2}(|\overline{Z_2+S}|+|\overline{(Z_3+p)(Z_2+S)}|)\int_0^1|\wt\bu|\Big(|\nabla_{e_\theta}|\wt\bu|^2|+|\nabla_{e_\varphi}|\wt\bu|^2|\Big) d\xi d\mathbf S^2\\
		=&J_{1}+J_{2}.
		\end{align*}
We first estimate $J_{1},$ then the estimate of $J_{2}$ follows similarly.
By the H\"older inequality, the interpolation inequality and Young's inequality we have
\begin{align*}
J_{1}&\leq |\wt\bu |\nabla_{e_\theta}\wt\bu | |_{2}|\wt\bu|_{4}\Big(|\overline{Z_2+S}|_{L^{4}(\mathbf{S}^{2})} +|\overline{(Z_3+p)(Z_2+S)}|_{L^{4}(\mathbf{S}^{2})} \Big)\\
&\leq |\wt\bu |\nabla_{e_\theta}\wt\bu | |_{2}|\wt\bu|_{4}\Big(|\bar{Z_2}|_{L^{4}(\mathbf{S}^{2})}+|\bar{S}|_{L^{4}(\mathbf{S}^{2})} +|\bar{Z_3} |_{L^{8}(\mathbf{S}^{2})}\\
&\hspace{3.5cm} +|\bar p|_{L^{8}(\mathbf{S}^{2})} + |\bar Z_2|_{L^{8}(\mathbf{S}^{2})}   +|\bar S|_{L^{8}(\mathbf{S}^{2})} \Big)\\
&\leq |\wt\bu |\nabla_{e_\theta}\wt\bu | |_{2}|\wt\bu|_{4}\Big(\|Z_2\|_{2}+\|Z_{3} \|_{2} +|\bar p|_{L^{4}(\mathbf{S}^{2})}^{1/2}(|\bar p|_{L^{4}(\mathbf{S}^{2})}^{1/2}+ |\bar\nabla p|_{L^{2}(\mathbf{S}^{2})}^{1/2} )\\
&\hspace{3.5cm}  +|\bar S|_{L^{4}(\mathbf{S}^{2})}^{1/2}(|\bar S|_{L^{4}(\mathbf{S}^{2})}^{1/2}+ |\bar\nabla S|_{L^{2}(\mathbf{S}^{2})}^{1/2} )\Big)\\
&\leq |\wt\bu |\nabla_{e_\theta}\wt\bu | |_{2}|\wt\bu|_{4}\Big(\|Z_2\|_{2}+\|Z_{3} \|_{2} +|p|_{4}^{1/2}(|p|_{4}^{1/2}+ |\nabla p|_{2}^{1/2} )\\ &\hspace{3.5cm}+|S|_{4}^{1/2}(|S|_{4}^{1/2}+ |\nabla S|_{2}^{1/2} )\Big)\\
&\leq \varepsilon \int_{\mathbf{D}}|\wt\bu|^{2} |\nabla_{e_\theta}\wt\bu |^{2} d\mathbf{D}+C|\wt\bu|_{4}^{2}( \|Z_2\|_{2}^{2}+\|Z_{3} \|_{2}^{2}+|p|_{4}\|p\|_{1} +|S|_{4}\|S\|_{1}).
\end{align*}
Similarly, we have the estimate of $J_{2},$
\begin{align*}
J_{2}&\leq\varepsilon \int_{\mathbf{D}}|\nabla_{e_\theta}|\wt\bu |^{2}|^{2} d\mathbf{D}+C|\wt\bu|_{4}^{2}( \|Z_2\|_{2}^{2}+\|Z_{3} \|_{2}^{2}+|p|_{4}\|p\|_{1} +|S|_{4}\|S\|_{1}).
\end{align*}
 By virtue of Estimates of $J_1$ and $J_{2}$ we have
 \begin{align*}
 |I_8+I_{9}|\leq& \varepsilon \int_{\mathbf{D}}(|\nabla_{e_\theta}|\wt\bu |^{2}|^{2} +   |\wt\bu|^{2} |\nabla_{e_\theta}\wt\bu |^{2})d\mathbf{D}+C|\wt\bu|_{4}^{2}( \|Z_2\|_{2}^{2}+\|Z_{3} \|_{2}^{2}+|p|_{4}\|p\|_{1} +|S|_{4}\|S\|_{1}).
 \end{align*}

		Throughout the estimates $I_1$-$I_{10}$, we have
		\begin{align}
		\label{eqn:5.39}
			&\frac{d|\wt\bu|_4^4}{dt}+\int_{\mathbf D}\Big(|\nabla_{e_\theta}\wt\bu^2|^2+|\nabla_{e_\varphi}\wt\bu^2|^2+|\partial_\xi|\wt\bu|^2|^2\Big)d\mathbf D+\int_{\mathbf D}|\wt\bu|^2\Big(|\nabla_{e_\theta}\wt\bu|^2+|\nabla_{e_\varphi}\wt\bu|^2+|\partial_\xi\wt\bu|^2+|\wt\bu|^2\Big)d\mathbf D\notag\\
			&\leq C(1+\norm{Z_1}_2^4+\norm{\bu}_1^2)|\wt\bu|_4^4\notag\\
            &\hspace{5mm}+C(\norm{Z_1}_3^2+\norm{Z_1}_3\norm{\bu}_1+\norm{Z_1}_1)|\wt\bu|_4^3\notag\\
			&\hspace{5mm}+C \|Z_{1}\|_{3}^{4/3}\|\bu\|_{1}^{4/3}|\wt\bu|_{4}^{8/3}+C(\norm{Z_1}_{2}^{16/5}+\norm{Z_1}_2^{8/5}\norm{\bu}_1^{8/5})|\wt\bu|_4^{12/5}\notag\\
			&\hspace{5mm}+C(\norm{Z_1}_2^4+ \|Z_2\|_{2}^{2}+\|Z_3\|_{2}^{2} +\norm{Z_1}_2^2\norm{\bu}_1^2+|S|_4\|S\|_{1}+|p|_4\|p\|_{1})|\wt\bu|_4^2,
		\end{align}
		and we also have
		\begin{align}
		\label{eqn:5.40}
&			\frac{d|\wt\bu|_4^2}{dt}+ |\wt \bu|_{4}^{2}\notag\\
&\leq C(1+\|Z_{1}\|_{2}^{4}+\|\bu\|_{1}^{2})|\wt\bu|_{4}^{2}+C\|Z_{1}\|_{3}^{2}\|\bu\|_{1}^{2}\notag\\
&\hspace{5mm}+C |S|_{4}\|S\|_{1}+C|p|_{4}\|p\|_{1}+C(\|Z_{1}\|_{1}^{2}+\|Z_{2}\|_{2}^{2}+\|Z_{3}\|_{2}^{2}+\|Z_{1}\|_{3}^{4}).
		\end{align}
		Now applying the Gronwall's inequality gives that
		\begin{align}
			\label{eqn:L4ofwtu}
			&\sup_{t\in[0,\tau)}|\wt\bu(t)|_4^4+\int_0^\tau\int_{\mathbf D}\Big(|\nabla_{e_\theta}\wt\bu^2|^2+|\nabla_{e_\varphi}\wt\bu^2|^2+|\partial_\xi|\wt\bu|^2|^2\Big)d\mathbf D\notag\\
			&\hspace{5mm}+\int_0^\tau \int_{\mathbf D}|\wt\bu|^2\Big(|\nabla_{e_\theta}\wt\bu|^2+|\nabla_{e_\varphi}\wt\bu|^2+|\partial_\xi\wt\bu|^2+|\wt\bu|^2\Big)d\mathbf D \notag\\
			&\leq C(\tau,Z_1,Z_2,Z_3, U_0).
		\end{align}
		Taking the inner product of Eq. (\ref{eqn:baru}) with $-\Delta\bar{\bu}$ in $L^2(\mathbf S^2)$, we get
		\begin{align}
		\label{eqn:5.42}	 &\hspace{5mm}\frac{1}{2}\partial_{t}(|\nabla_{e_\theta}\bar{\bu}|_{2}^2+|\nabla_{e_\varphi}\bar{\bu}|_{2}^2+|\bar{\bu}|_{2}^2)+|\Delta\bar{\bu}|_2^2\notag\\
			&=\int_{\mathbf S^2}\overline{\nabla_{\wt Z_1+\wt \bu}(\wt Z_1+\wt \bu)+(\wt Z_1+\wt\bu)\text{div }(\wt Z_1+\wt\bu)}\cdot\Delta\bar{\bu}d\mathbf S^2+\int_{\mathbf S^2}\nabla_{\bar Z_1+\bar{\bu}}(\bar Z_1+\bar{\bu})\cdot\Delta\bar{\bu}d\mathbf S^2\notag\\
			&\hspace{5mm}+\int_{\mathbf S^2}\frac{f}{R_0}\vec k\times\bar Z_1\cdot\Delta\bar{\bu}d\mathbf S^2,
		\end{align}
		because by integration by parts, we have
		$$\int_{\mathbf S^2}\frac{f}{R_{0}}\vec{k}\times \bar\bu \cdot \Delta \bar\bu d\mathbf S^2=0,\ \   \int_{\mathbf S^2}\nabla\Phi_s\cdot\Delta\bar{\bu}d\mathbf S^2=0,$$
		and
		$$\int_{\mathbf S^2}\nabla\int_0^1\int_{\xi}^{1}\frac{br_s}{r}[1+a(Z_3+p)(Z_2+S)]d\xi'd\xi\cdot\Delta\bar{\bu}d\mathbf S^2=0.$$
		Applying the H\"older inequality, the Minkowski inequality and the Sobolev embedding theorem, we first have
		\begin{align*}
			&\bigg|\int_{\mathbf S^2}\overline{\nabla_{\wt Z_1+\wt \bu}(\wt Z_1+\wt \bu)+(\wt Z_1+\wt\bu)\text{div }(\wt Z_1+\wt\bu)}\cdot\Delta\bar{\bu}d\mathbf S^2\bigg|\\
			&\leq 2|\Delta\bar{\bu}|_2\bigg(\int_{\mathbf D}(|\wt Z_1|^2+|\wt\bu|^2)(|\nabla_{e_\theta}(\wt Z_1+\wt\bu)|^2+|\nabla_{e_\varphi}(\wt Z_1+\wt\bu)|^2)d\mathbf D\bigg)^{1/2}\\
			&\leq \eps|\Delta\bar{\bu}|_2^2+C\int_{\mathbf{D}}|\wt\bu|^{2}(|\nabla_{e_{\theta}}\wt\bu |^{2}+ |\nabla_{e_{\varphi}}\wt\bu |^{2}) d\mathbf{D}\\
 &\hspace{5mm}+C(|\wt\bu|_4^2\norm{Z_1}_2^2+|Z_1|_\infty\norm{\bu}_1^2+|Z_1|_\infty\norm{Z_1}_1^2).
		\end{align*}
By the H\"older inequality and Young's inequality we have
\begin{align*}
&\hspace{5mm}\int_{\mathbf{S}^{2}}\bar \bu_{\theta}\nabla_{e_{\theta}}\bar \bu\cdot \Delta \bar \bu d\mathbf{S}^{2}\\
&\leq |\Delta \bar \bu|_{L^{2}(\mathbf{S}^{2})}|\nabla_{e_{\theta}} \bar \bu|_{L^{4}(\mathbf{S}^{2})}|\bar \bu_{\theta}|_{L^{4}(\mathbf{S}^{2})}\\
&\leq C|\Delta \bar \bu|_{L^{2}(\mathbf{S}^{2})}|\nabla_{e_{\theta}} \bar \bu|_{L^{2}(\mathbf{S}^{2})}^{1/2}|\Delta \bar\bu|_{L^{2}(\mathbf{S}^{2})}^{1/2}
|\bar \bu_{\theta}|_{L^{2}(\mathbf{S}^{2})}^{1/2}( |\nabla_{e_{\theta}} \bar \bu|_{L^{2}(\mathbf{S}^{2})}^{1/2}+ |\nabla_{e_{\varphi}} \bar \bu|_{L^{2}(\mathbf{S}^{2})}^{1/2} )\\
&\leq \varepsilon |\Delta \bar \bu|_{2}^{2}+C|\bu|_{2}^{2}\|\bu\|_{1}^{2}|\nabla_{e_{\theta}} \bar \bu|_{2}^{2}.
\end{align*}
		Similarly, we can get
		\begin{align*}
			 &\hspace{5mm}\bigg|\int_{\mathbf S^2}\nabla_{\bar Z_1+\bar{\bu}}(\bar Z_1+\bar{\bu})\cdot\Delta\bar{\bu}d\mathbf S^2\bigg|\\
			&\leq\bigg|\int_{\mathbf S^2} \bar \bu_{\theta} \nabla_{e_{\theta}} \bar \bu \cdot\Delta\bar{\bu}d\mathbf S^2
+\int_{\mathbf S^2} \bar \bu_{\theta} \nabla_{e_{\theta}} \bar Z_{1} \cdot\Delta\bar{\bu}d\mathbf S^2\\
&\hspace{5mm}+\int_{\mathbf S^2} \bar Z_{1,\theta} \nabla_{e_{\theta}} \bar Z_{1} \cdot\Delta\bar{\bu}d\mathbf S^2
+\int_{\mathbf S^2} \bar Z_{1,\theta} \nabla_{e_{\theta}} \bar \bu \cdot\Delta\bar{\bu}d\mathbf S^2\bigg|\\
&\hspace{5mm}+\bigg|\int_{\mathbf S^2} \bar \bu_{\varphi} \nabla_{e_{\varphi}} \bar \bu \cdot\Delta\bar{\bu}d\mathbf S^2
+\int_{\mathbf S^2} \bar \bu_{\varphi} \nabla_{e_{\varphi}} \bar Z_{1} \cdot\Delta\bar{\bu}d\mathbf S^2\\
&\hspace{5mm}+\int_{\mathbf S^2} \bar Z_{1,\varphi} \nabla_{e_{\varphi}} \bar Z_{1} \cdot\Delta\bar{\bu}d\mathbf S^2
+\int_{\mathbf S^2} \bar Z_{1,\varphi} \nabla_{e_{\varphi}} \bar \bu \cdot\Delta\bar{\bu}d\mathbf S^2\bigg|\\
&\leq \varepsilon |\Delta \bar \bu|_{2}^{2}+C|\bu|_{2}^{2}\|\bu\|_{1}^{2}(|\nabla_{e_{\theta}} \bar \bu|_{2}^{2}+|\nabla_{e_{\varphi}} \bar \bu|_{2}^{2})\\
&\hspace{5mm}+C(|\nabla_{e_{\theta}}Z_{1}|_{\infty}^{2}+  |\nabla_{e_{\varphi}}Z_{1}|_{\infty}^{2} )|\bu|_{2}^{2}+C|Z_{1}|_{\infty}\|\bu\|_{1}^{2}+C\|Z_{1}\|_{2}^{4}\\
&\leq \varepsilon |\Delta \bar \bu|_{2}^{2}+C|\bu|_{2}^{2}\|\bu\|_{1}^{2}(|\nabla_{e_{\theta}} \bar \bu|_{2}^{2}+|\nabla_{e_{\varphi}} \bar \bu|_{2}^{2})\\
&\hspace{5mm}+C\|Z_{1}\|_{3}^{2}|\bu|_{2}^{2}+C\|Z_{1}\|_{2}\|\bu\|_{1}^{2}+C\|Z_{1}\|_{2}^{4}.
		\end{align*}

Therefore, combining the above arguments and \eqref{eqn:5.42} yields
\begin{align}	
\label{eqn:5.43}
 &\hspace{5mm}\frac{1}{2}\partial_{t}(|\nabla_{e_\theta}\bar{\bu}|_{2}^2+|\nabla_{e_\varphi}\bar{\bu}|_{2}^2+|\bar{\bu}|_{2}^2)+|\Delta\bar{\bu}|_2^2\notag\\
&\leq C\int_{\mathbf{D}}|\wt\bu|^{2}(|\nabla_{e_{\theta}} \wt\bu|^{2}+|\nabla_{e_{\varphi}} \wt\bu|^{2} )d\mathbf{D}\notag\\
&\hspace{5mm}+C|\bu|_{2}^{2}\|\bu\|_{1}^{2}( |\nabla_{e_{\theta}} \bar\bu|_{2}^{2}+|\nabla_{e_{\varphi}} \bar\bu|_{2}^{2}  )\notag\\
&\hspace{5mm}+C\|Z_{1}\|_{2}\|\bu\|_{1}^{2}+C|\wt\bu|_{4}^{2}\|Z_{1}\|_{2}^{2}+C|\bu|_{2}^{2}\|Z_{1}\|_{3}^{2}\notag\\
&\hspace{5mm}+C\|Z_{1}\|_{2}^{4}+C\|Z_{1}\|_{2}\|Z_{1}\|_{1}^{2}+C|Z_{1}|_{2}^{2}.
	\end{align}
Applying the Gronwall's inequality and \eqref{eqn:L4ofwtu} to \eqref{eqn:5.43} we obtain
		\begin{equation}
		\label{eqn:L2ofgradientbaru}
			\sup_{t\in[t_0,\tau)}\Big(|\nabla_{e_\theta}\bar{\bu}(t)|_{2}^2+|\nabla_{e_\varphi}\bar{\bu}(t)|_{2}^2+|\bar\bu(t)|_{2}^{2}\Big)\leq C(\tau,Z_1,Z_2,Z_3, U_0).
		\end{equation}
		\subsection{$H^1$ estimates of $\v, T, q$.}
		Taking the derivative of (\ref{eqn:u}) with respect to $\xi$, we get
		\begin{align}
			&\partial_t\bu_\xi-\Delta\bu_\xi-\partial_{\xi\xi}\bu_\xi+\nabla_{Z_1+\bu}(\partial_\xi Z_1+\bu_\xi)+\nabla_{\partial_\xi Z_1+\bu_\xi}(Z_1+\bu)\notag\\
			&-(\text{div }(Z_1+\bu))(\partial_\xi Z_1+\bu_\xi)+w(Z_1+\bu)(\partial_{\xi\xi}Z_1+\bu_{\xi\xi})\notag\\
			&=-\frac{f}{R_0}\vec k\times(\partial_\xi Z_1+\bu_\xi)+\frac{br_s}{r}\nabla[(1+aq)T]+\gamma\partial_\xi Z_1.
		\end{align}
		Applying integration by parts, we have
		$$\int \Big(\nabla_\bu\bu_\xi+w(\bu)\bu_{\xi\xi}\Big)\cdot\bu_\xi d\mathbf D=0.$$
		Since $\Phi_s$ is independent of $\xi$, also with $[\frac{f}{R_0}\vec k\times\bu_\xi]\cdot\bu_\xi=0$, then taking the inner product with $\bu_\xi$ in $(L^2(\D))^2$, we have
		\begin{align}
			&\frac{1}{2}\frac{d|\bu_\xi|_2^2}{dt}+\int_{\mathbf D}[|\nabla_{e_\theta}\bu_\xi|^2+|\nabla_{e_\varphi}\bu_\xi|^2+|\bu_\xi|^2]d\mathbf D+\int_{\mathbf D}|\bu_{\xi\xi}|^2\notag\\
			&=-\int_{\mathbf D}\nabla_{\bu}\partial_\xi Z_1\cdot\bu_\xi d\mathbf D-\int_{\mathbf D}\nabla_{Z_1}(\partial_\xi Z_1+\bu_\xi)\cdot\bu_\xi d\mathbf D-\int_{\mathbf D}\nabla_{\partial_\xi Z_1+\bu_\xi}(Z_1+\bu)\cdot\bu_\xi d\mathbf D\notag\\
			&\hspace{5mm}+\int_{\mathbf D}\text{div }(Z_1+\bu)(\partial_\xi Z_1+\bu_\xi)\cdot\bu_\xi d\mathbf D+\int_{\mathbf D}w(Z_{1})\bu_{\xi\xi}\cdot\bu_\xi d\mathbf D+\int_{\mathbf D}w(Z_1+\bu)\partial_{\xi\xi}Z_1\cdot\bu_\xi d\mathbf D\notag\\
			&\hspace{5mm}+\int_{\mathbf D}\frac{br_s}{r}\nabla[(1+aq)T]\cdot\bu_\xi d\mathbf D+\gamma\int_{\mathbf D}\partial_\xi Z_1\cdot\bu_\xi d\mathbf D  - \int_{\mathbf D} \frac{f}{R_0}\vec k\times \partial_\xi Z_1  \cdot\bu_\xi d\mathbf D.
		\end{align}
		By the H\"older inequality, the interpolation inequality, the Sobolev embedding theorem and Young's inequality, we have
		\begin{align}
		\label{eqn:5.47}
			&\int_{\mathbf D}\nabla_{\bu}\partial_\xi Z_1\cdot\bu_\xi d\mathbf D\leq\norm{Z_1}_2|\bu|_4|\bu_\xi|_4\notag\\
			 &\leq\norm{Z_1}_2|\bu|_4|\bu_\xi|^{1/4}_2\Big(|\nabla_{e_\theta}\bu_\xi|_2^{3/4}+|\nabla_{e_\varphi}\bu_\xi|_2^{3/4}+|\bu_{\xi\xi}|_2^{3/4}+|\bu_\xi|_2^{3/4}\Big)\notag\\
			 &\leq\eps\Big(|\nabla_{e_\theta}\bu_\xi|_2^2+|\nabla_{e_\varphi}\bu_\xi|_2^2+|\bu_{\xi\xi}|_2^2\Big)\notag\\
&\hspace{5mm}+C\norm{Z_1}_2^{8/5}|\bu|_{4}^{8/5}|\bu_{\xi}|_{2}^{2/5}+C\|Z_{1}\|_{2}|\bu|_{4}|\bu_{\xi}|_{2}\notag\\
&\leq\eps\Big(|\nabla_{e_\theta}\bu_\xi|_2^2+|\nabla_{e_\varphi}\bu_\xi|_2^2+|\bu_{\xi\xi}|_2^2\Big)\notag\\
&\hspace{5mm}+C\norm{Z_1}_2^{8/5}(|\wt\bu|_{4}^{8/5}+\|\bar\bu\|_{1}^{8/5})|\bu_{\xi}|_{2}^{2/5}+C\|Z_{1}\|_{2}(|\wt\bu|_{4}+
\|\bar\bu\|_{1})|\bu_{\xi}|_{2}\notag\\
&\leq\eps\Big(|\nabla_{e_\theta}\bu_\xi|_2^2+|\nabla_{e_\varphi}\bu_\xi|_2^2+|\bu_{\xi\xi}|_2^2+ |\bu_{\xi}|_{2}^{2}\Big)\notag\\
&\hspace{5mm}+C\|Z_{1}\|_{2}^{2}(|\wt\bu|_{4}^{2}+
\|\bar\bu\|_{1}^{2}).
\end{align}
		Similarly, we also obtain
		\begin{align}
			\int_{\mathbf D}\nabla_{Z_1}(\partial_\xi Z_1+\bu_\xi)\cdot\bu_\xi d\mathbf D\leq &|Z_1|_\infty\Big(|\nabla_{e_\theta}(\partial_\xi Z_1+\bu_\xi)|_2+|\nabla_{e_\varphi}(\partial_\xi Z_1+\bu_\xi)|_2\Big)|\bu_\xi|_2\notag\\
			\leq & \eps\Big(|\nabla_{e_\theta}\bu_\xi|_2^2+|\nabla_{e_\varphi}\bu_\xi|_2^2\Big)+C\norm{Z_1}_2^2(|\bu_\xi|_2^2+1).
		\end{align}
	By integration by parts, following the calculus in \eqref{eqn:5.47} we have
		\begin{align}
			&\int_{\mathbf D}\nabla_{\partial_\xi Z_1}(Z_1+\bu)\cdot\bu_\xi d\mathbf D\notag\\
&=\int_{\mathbf D}  \Big(\bu_\xi  \cdot(Z_1+\bu)\Big) \mathrm{div} \partial_\xi Z_1  d\mathbf D-\int_{\mathbf D}\nabla_{\partial_\xi Z_1}\bu_\xi  \cdot(Z_1+\bu) d\mathbf D\notag\\
			&\leq |\bu_\xi|_4|\bu+Z_1|_4\Big(|\nabla_{e_\theta}\partial_\xi Z_1|_2+|\nabla_{e_\varphi}\partial_\xi Z_1|_2\Big)+|\partial_\xi Z_1|_4|\bu+Z_1|_4\Big(|\nabla_{e_\theta}\bu_\xi|_2+|\nabla_{e_\varphi}\bu_\xi|_2\Big)\notag\\
			&\leq C\norm{Z_1}_2|\bu_\xi|_2^{1/4}|\bu+Z_1|_4\Big(|\nabla_{e_\theta}\bu_\xi|_2^{3/4}+|\nabla_{e_\varphi}\bu_\xi|_2^{3/4}+|\bu_\xi|_2^{3/4}+|\bu_{\xi\xi}|_2^{3/4}\Big)\notag\\
			&\hspace{5mm}+C\norm{Z_1}_2|\bu+Z_1|_4\Big(|\nabla_{e_\theta}\bu_\xi|_2+|\nabla_{e_\varphi}\bu_\xi|_2\Big)\notag\\
			&\leq \eps\Big(|\nabla_{e_\theta}\bu_\xi|_2^2+|\nabla_{e_\varphi}\bu_\xi|_2^2+|\bu_{\xi\xi}|_2^2\Big)+C\norm{Z_1}_2^2|\bu+Z_1|_4^2+C|\bu_\xi|_2^2.
		\end{align}
		Now repeating similar discussion with the H\"older inequality, the interpolation inequality and the Sobolev embedding theorem, we have
		\begin{align}
			&\int_{\mathbf D}\nabla_{\bu_\xi}(Z_1+\bu)\cdot\bu_\xi d\mathbf D\notag\\
&=-\int_{\mathbf D}\nabla_{\bu_\xi} \bu_\xi  \cdot (Z_1+\bu) d\mathbf D+\int_{\mathbf D} \bu_\xi  \cdot (Z_1+\bu) \mathrm{div} \bu_\xi  d\mathbf D\notag\\
			&\leq C|\bu|_4|\bu_\xi|_4\Big(|\nabla_{e_\theta}\bu_\xi|_2+|\nabla_{e_\varphi}\bu_\xi|_2\Big)+C|Z_1|_\infty|\bu_\xi|_2\Big(|\nabla_{e_\theta}\bu_\xi|_2+|\nabla_{e_\varphi}\bu_\xi|_2\Big)\notag\\
			&\leq C|\bu|_4|\bu_\xi|_2^{1/4}\Big(|\nabla_{e_\theta}\bu_\xi|_2^{3/4}+|\nabla_{e_\varphi}\bu_\xi|_2^{3/4}+|\bu_{\xi\xi}|_2^{3/4}+|\bu_\xi|_2^{3/4}\Big)\Big(|\nabla_{e_\theta}\bu_\xi|_2+|\nabla_{e_\varphi}\bu_\xi|_2\Big)\notag\\
			&\hspace{5mm}+C|Z_1|_\infty|\bu_\xi|_2\Big(|\nabla_{e_\theta}\bu_\xi|_2+|\nabla_{e_\varphi}\bu_\xi|_2\Big)\notag\\
			&\leq \eps\Big(|\nabla_{e_\theta}\bu_\xi|_2^2+|\nabla_{e_\varphi}\bu_\xi|_2^2+|\bu_{\xi\xi}|_2^2 +|\bu_{\xi}|_{2}^{2}\Big)+C(|\bu|_4^8+\norm{Z_1}_2^2)|\bu_\xi|_2^2.
		\end{align}
By integration by parts, using the H\"{o}lder inequality, interpolation inequality and Young's inequality we obtain
\begin{align*}
&\hspace{5mm}\int_{\mathbf{D}}(\mathrm{div} \bu)\bu_{\xi}\cdot \bu_{\xi}d\mathbf D\\
&=2\int_{\mathbf D}\Big(\nabla_{\bu}\bu_{\xi}\Big) \bu_{\xi} d \mathbf{D}\\
&\leq 2\int_{\mathbf D} (|\nabla_{e_{\theta}}\bu_{\xi} |+ |\nabla_{e_{\varphi}}\bu_{\xi} |)|\bu||\bu_{\xi}|d \mathbf{D}\\
&\leq C(|\nabla_{e_{\theta}}\bu_{\xi}|_{2}+ |\nabla_{e_{\varphi}}\bu_{\xi}|_{2})|\bu|_{4}|\bu_{\xi}|_{2}^{1/4}\|\bu_{\xi}\|_{1}^{3/4}\\
&\leq \varepsilon \|\bu_{\xi}\|_{1}^{2}+C|\bu_{\xi}|_{2}^{2}|\bu|_{4}^{8}.
\end{align*}
Then
		\begin{align}
			&\int_{\mathbf D}\text{div }(Z_1+\bu)(\partial_\xi Z_1+\bu_\xi)\cdot \bu_\xi d\mathbf D\notag\\
			&=\int_{\mathbf{D}}(\mathrm{div} \bu)\bu_{\xi}\cdot \bu_{\xi}d\mathbf D+\int_{\mathbf{D}}(\mathrm{div} Z_{1})(\partial_{\xi}Z_{1})\cdot \bu_{\xi}d\mathbf D\notag\\
&\hspace{5mm}+\int_{\mathbf{D}}(\mathrm{div} Z_{1})\bu_{\xi}\cdot \bu_{\xi}d\mathbf D+ \int_{\mathbf{D}}(\mathrm{div} \bu) (\partial_{\xi}Z_{1})\cdot \bu_{\xi}d\mathbf D\notag\\
&\leq \varepsilon \|\bu_{\xi}\|_{1}^{2}+C|\bu_{\xi}|_{2}^{2}|\bu|_{4}^{8}+|\mathrm{div} Z_{1}|_{4}|\partial_{\xi} Z_{1}|_{4}|\bu_{\xi}|_{2}\notag\\
&\hspace{5mm}+|\mathrm{div }Z_{1}|_{\infty}|\bu_{\xi}|_{2}^{2}+|\partial_{\xi}Z_{1}|_{\infty}\|\bu\|_{1}|\bu_{\xi}|_{2}\notag\\
&\leq \varepsilon \|\bu_{\xi}\|_{1}^{2}+C|\bu_{\xi}|_{2}^{2}|\bu|_{4}^{8}+C\|Z_{1}\|_{2}^{2}|\bu_{\xi}|_{2}\notag\\
&\hspace{5mm}+C\|Z_{1}\|_{3}|\bu_{\xi}|_{2}^{2}+C\|Z_{1}\|_{3}\|\bu\|_{1}|\bu_{\xi}|_{2}.
		\end{align}
		Applying integration by parts,
		\begin{align*}
		\int_{\mathbf D}w(\bu)\partial_{\xi\xi}Z_1\cdot\bu_\xi d\mathbf D=&\int_{\mathbf D}\text{div }\bu\partial_\xi Z_1\cdot\bu_\xi d\mathbf D-\int_{\mathbf D}w(\bu)\partial_\xi Z_1\cdot\bu_{\xi\xi} d\mathbf D\notag\\
		=&-\int_{\D}(\nabla_{\bu}\partial_{\xi}Z_{1})\bu_{\xi} d\mathbf D-\int_{\mathbf D}\partial_\xi Z_1\cdot(\nabla_{\bu}\bu_\xi )d\mathbf D	\notag\\
		&-\int_{\mathbf D}w(\bu)\partial_\xi Z_1\cdot\bu_{\xi\xi} d\mathbf D.
		\end{align*}
		Due to the H\"older inequality, the Sobolev embedding theorem, the Minkowski inequality and the interpolation inequality, one can obtain
		\begin{align}
			&\int_{\D}w(\bu)\partial_{\xi\xi}Z_{1}\cdot\bu_\xi d\mathbf D\notag\\
& \leq \|Z_{1}\|_{2}|\bu|_{4}|\bu_{\xi}|_{4}+\|\bu_{\xi}\|_{1}|\bu|_{4}|\partial_{\xi}Z_{1}|_{4}+|\partial_{\xi}Z_{1}|_{\infty}\|\bu\|_{1}|\bu_{\xi
\xi}|_{2} \notag\\
&\leq\eps\Big(|\nabla_{e_\theta}\bu_\xi|_2^2+|\nabla_{e_\varphi}\bu_\xi|_2^2+|\bu_{\xi\xi}|_2^2\Big)+C|\bu_\xi|_2^2\notag\\ &+C\norm{Z_1}_2^2\Big(|\nabla_{e_\theta}\bar\bu|_2^2+|\nabla_{e_\varphi}\bar\bu|_2^2+|\wt\bu|_4^2\Big)+C\norm{Z_1}_3^2\|\bu\|_{1}^{2}.
		\end{align}
		Similarly, we also have
		\begin{align}
			&\int_{\mathbf D}w(Z_1)(\partial_{\xi\xi}Z_1+\bu_{\xi\xi})\cdot\bu_\xi d\mathbf D\notag\\\leq &|\bu_{\xi\xi}|_2|\bu_\xi|_4|w(Z_1)|_4+|\partial_{\xi\xi}Z_1|_2\Big(|\nabla_{e_\theta}Z_1|_4+|\nabla_{e_\varphi}Z_1|_4\Big)|\bu_\xi|_4\notag\\
			 \leq&C(|\bu_{\xi\xi}|_2\norm{Z_1}_2+\norm{Z_1}_2^2)|\bu_\xi|_2^{1/4}\Big(|\nabla_{e_\theta}\bu_\xi|_2^{3/4}+|\nabla_{e_\varphi}\bu_\xi|_2^{3/4}+|\bu_{\xi\xi}|_2^{3/4}\Big)\notag\\
			 \leq&\eps\Big(|\nabla_{e_\theta}\bu_\xi|_2^2+|\nabla_{e_\varphi}\bu_\xi|_2^2+|\bu_{\xi\xi}|_2^2\Big)+C\norm{Z_1}_2^8|\bu_\xi|_2^2+C\norm{Z_1}_2^{16/5}|\bu_\xi|_2^{2/5}.
		\end{align}
	   Finally we obtain that
		\begin{align}
			\int_{\mathbf D}\frac{br_s}{r}\nabla[(1+aq)T]\cdot\bu_\xi d\mathbf D=&-\int\frac{br_s}{r}(1+aq)T\text{div }\bu_\xi d\mathbf D\notag\\
			\leq&\eps\Big(|\nabla_{e_\theta}\bu_\xi|_2^2+|\nabla_{e_\varphi}\bu_\xi|_2^2+|\bu_{\xi\xi}|_2^2\Big)+C|T|_2^2+C|q|_4^2|T|_4^2.
		\end{align}
		With the above estimates altogether, we arrive at
		\begin{align}
			\label{eqn:H1ofv}
			&\partial_t|\bu_\xi|_2^2+|\nabla_{e_\theta}\bu_\xi|_2^2+|\nabla_{e_\varphi}\bu_\xi|_2^2+|\bu_{\xi\xi}|_2^2\notag\\
			&\leq C(|\wt\bu|_4^8+\|\bar\bu\|_1^8+\norm{Z_1}_2^8+\|Z_{1}\|_{3}+1)|\bu_\xi|_2^2\notag\\
&\hspace{5mm}+C\|\bu\|_{1}^{2}\norm{Z_1}_3^2+C\|Z_{1}\|_{2}^{4}
			+C|T|_2^2+C|q|_4^2|T|_4^2.
		\end{align}
		Therefore, by the Gronwall inequality, and $L^2, L^4$ estimates of $T$ and $q$, one can have
		\begin{align}		
		\label{eqn:5.56}
		 &\sup_{t\in[t_0,\tau)}|\bu_\xi(t)|_2^2+\int_0^\tau\Big(|\nabla_{e_\theta}\bu_\xi(s)|_2^2+|\nabla_{e_\varphi}\bu_\xi(s)|_2^2+|\bu_{\xi\xi}(s)|_2^2\Big)ds\notag\\
			&\leq C(\tau,Q_T,Q_q,Z_1,Z_2,Z_3,U_0).
		\end{align}
Taking the derivative of \eqref{eqn:p} with respect to $\xi,$ we obtain
\begin{align}
&\partial_{t} p_{\xi} +\nabla_{\partial_{\xi}Z_{1}+\bu_{\xi}}(Z_{3}+p)+\nabla_{Z_{1}+\bu}(\partial_{\xi}Z_{3}+p_{\xi} )\notag\\
&\hspace{5mm}-(\mathrm{div} (Z_{1}+\bu) )\partial_{\xi}(Z_{3}+p)+w(Z_{1}+\bu)\partial_{\xi\xi}(Z_{3}+p)\notag\\
&\hspace{5mm}=\Delta p_{\xi}+\partial_{\xi\xi}p_{\xi}+\partial_{\xi}Q_{q}+\gamma \partial_{\xi}Z_{3}.
\end{align}
Applying integration by parts, we have
\begin{align*}
\int_{\D}[(\nabla_{\bu} p_{\xi})p_{\xi}+w(\bu)p_{\xi\xi}p_{\xi} ]d\D=0.
\end{align*}
Then taking inner product with $p_{\xi}$ in $L^{2}(\D)$ yields that
\begin{align}
\label{eqn:5.58}
&\frac{1}{2}\partial_{t}|p_{\xi}|_{2}^{2}+|\nabla p_{\xi}|_{2}^{2}+|p_{\xi\xi}|_{2}^{2}-\int_{\mathbf{S}^{2}}p_{\xi\xi}(\xi=1)p_{\xi}(\xi=1)d\mathbf{S}^{2}\notag\\
&=-\int_{\D}\nabla_{\partial_{\xi}Z_{1}+\bu_{\xi}}(Z_{3}+p)p_{\xi}d\D\notag\\
&\hspace{5mm}-\int_{\D}(\nabla_{Z_{1}}p_{\xi}+ \nabla_{Z_{1}+\bu}\partial_{\xi}Z_{3})p_{\xi}d\D\notag\\
&\hspace{5mm}+\int_{\D}(\mathrm{div} (Z_{1}+\bu) )\partial_{\xi}(Z_{3}+p)p_{\xi}d\D\notag\\
&\hspace{5mm}-\int_{\D}w(Z_{1})p_{\xi\xi} p_{\xi}d\D-\int_{\D}w(Z_{1}+\bu)\partial_{\xi\xi}Z_{3} p_{\xi}d\D\notag\\
&\hspace{5mm}+\int_{\D}(\partial_{\xi}Q_{q}+\gamma \partial_{\xi}Z_{3} )p_{\xi}d\D.
\end{align}
Taking the trace on $\xi=1$ of \eqref{eqn:p}, we have
\begin{align}
\label{eqn:5.59}
&\hspace{5mm}-\int_{\mathbf{S}^{2}}p_{\xi}|_{\xi=1}p_{\xi\xi}|_{\xi=1}d\mathbf{S}^{2}\notag\\
&=\frac{\beta}{2} \partial_{t}|p|_{\xi=1}|_{2}^{2}+\beta|\nabla p|_{\xi=1} |_{2}^{2}+\beta\int_{\mathbf{S}^{2}}(p\nabla_{Z_{1}+\bu}(Z_{3}+p))|_{\xi=1}d\mathbf{S}^{2}\notag\\
&\hspace{5mm}-\beta\int_{\mathbf{S}^{2}}(pQ_{q})|_{\xi=1}d\mathbf{S}^{2}-\gamma \beta\int_{\mathbf{S}^{2}}(pZ_{3})|_{\xi=1}d\mathbf{S}^{2}.
\end{align}
Since
\begin{align}
\label{eqn:5.60}
&\hspace{5mm}\int_{\mathbf{S}^{2}}(p\nabla_{\bu+Z_{1}}Z_{3})|_{\xi=1}d\mathbf{S}^{2}\notag\\
&\leq |(\bu+Z_{1})|_{\xi=1}|_{L^{2}(\mathbf{S}^{2})}|p|_{\xi=1}|_{L^{4}(\mathbf{S}^{2})}|\nabla Z_{3}|_{\xi=1}|_{L^{4}(\mathbf{S}^{2})}\notag\\
&=\int_{0}^{1}|p|_{\xi=1}|_{L^{4}(\mathbf{S}^{2})}d\xi\int_{0}^{1} |(\bu+Z_{1})|_{\xi=1}|_{L^{2}(\mathbf{S}^{2})}d\xi\notag\\
&\hspace{5mm}\times \int_{0}^{1}|\nabla Z_{3}|_{\xi=1}|_{L^{4}(\mathbf{S}^{2})}d\xi\notag\\
&\leq |p|_{\xi=1}|_{4}\Big(\int_{\D} |(\bu+Z_{1})|_{\xi=1}|^{2}d\D\Big)^{1/2}\notag\\
&\hspace{5mm}\times \Big(\int_{\D}|\nabla Z_{3}|_{\xi=1}|^{4}d\D\Big)^{1/4}.
\end{align}
To obtain the estimate of \eqref{eqn:5.60}, we first consider
\begin{align}
\label{eqn:5.61}
&\Big(\int_{\mathbf{S}^{2}} |(\bu+Z_{1})|_{\xi=1}|^{2}d\mathbf{S}^{2}\Big)^{1/2}=\Big(\int_{\D} |(\bu+Z_{1})|_{\xi=1}|^{2}d\D\Big)^{1/2}\notag\\
&=\Big(\int_{\D}(\int_{\xi}^{1}\partial_{\xi}|\bu +Z_{1} |^{2}d\xi+|\bu +Z_{1}|^{2}) d\D\Big)^{1/2}\notag\\
&\leq\Big(\int_{\D}2 |\bu_{\xi} + \partial_{\xi}Z_{1}|\cdot |\bu +Z_{1} |d\D\Big)^{1/2}+|\bu+Z_{1}|_{2}\notag\\
&\leq |\bu_{\xi}+\partial_{\xi}Z_{1}|_{2}+2|\bu+Z_{1}|_{2},
\end{align}
where the first equality follows by $\mathbf{D}=\mathbf{S}^{2}\times (0,1).$   Then
\begin{align}
\label{eqn:5.62}
&\hspace{5mm}\Big(\int_{\D}|\nabla Z_{3}|_{\xi=1}|^{4}d\D\Big)^{1/4}\notag\\
&= \Big(\int_{\D}(\int_{\xi}^{1}\partial_{\xi}|\nabla Z_{3}|^{4}d\xi +|\nabla Z_{3}|^{4})d\D \Big)^{1/4}\notag\\
&\leq |\nabla Z_{3}|_{4}+\Big(\int_{\D}4 (\partial_{\xi}\nabla Z_{3}\cdot \nabla Z_{3})|\nabla Z_{3}|^{2}d\D\Big)^{1/4}\notag\\
&\leq C\|Z_{3}\|_{2}+C|\partial_{\xi}\nabla Z_{3} |_{4}^{1/4}|\nabla Z_{3}|_{4}^{3/4}
\leq C\|Z_{3}\|_{3}.
\end{align}
Combining \eqref{eqn:5.60}-\eqref{eqn:5.62} yields that
\begin{align}
\label{eqn:5.63}
&\hspace{5mm}\int_{\mathbf{S}^{2}}(p\nabla_{\bu+Z_{1}}Z_{3})|_{\xi=1}d\mathbf{S}^{2}\notag\\
&\leq C|p|_{\xi=1}|_{4}(|\bu_{\xi}+\partial_{\xi}Z_{1}|_{2}+|\bu +Z_{1}|_{2})\|Z_{3}\|_{3}.
\end{align}
Since $\D= \mathbf{S}^{2}\times (0,1),$ by integration by parts,
\begin{align}
&\hspace{5mm}\int_{\mathbf{S}^{2}}(p\nabla_{Z_{1}+\bu}p)|_{\xi=1}d\mathbf{S}^{2}\notag\\
&=\frac{1}{2}\int_{\mathbf{S}^{2}}(p^{2}\mathrm{div} (Z_{1}+\bu))|_{\xi=1}d\mathbf{S}^{2}\notag\\
&=\frac{1}{2}\int_{\D}(p^{2}\mathrm{div} (Z_{1}+\bu))|_{\xi=1}d\D\notag\\
&=\frac{1}{2}\int_{\D}(p^{2}|_{\xi=1}\int_{\xi}^{1}\partial_{\xi}\mathrm{div}(Z_{1}+\bu)d\xi +p^{2}|_{\xi=1}\mathrm{div} (Z_{1}+\bu))d\mathbf{S}^{2}\notag\\
&\leq |p|_{\xi=1}|_{4}^{4}+\|\partial_{\xi}Z_{1}+\bu_{\xi}\|_{1}^{2}+\|Z_{1}+\bu\|_{1}^{2}.
\end{align}
Taking an analogous argument of \eqref{eqn:5.61} we get
\begin{equation}
\label{eqn:5.65}
|Q_{q}|_{\xi=1}|_{L^{2}(\mathbf{S}^{2})}\leq 2|Q_{q}|_{2}+|\partial_{\xi}Q_{q}|_{2},
\end{equation}
and
\begin{align}
\label{eqn:5.66}
|Z_{3}|_{\xi=1}|_{L^{2}(\mathbf{S}^{2})}\leq 2|Z_{3}|_{2}+|\partial_{\xi}Z_{3}|_{2}.
\end{align}
In view of \eqref{eqn:5.59}, \eqref{eqn:5.63}-\eqref{eqn:5.66}, the H\"older inequality  and Young's inequality we obtain
\begin{align}
\label{eqn:5.67}
&\hspace{5mm}-\int_{\mathbf{S}^{2}}p_{\xi}|_{\xi=1}p_{\xi\xi}|_{\xi=1}d\mathbf{S}^{2}\notag\\
&\leq\frac{\beta}{2} \partial_{t}|p|_{\xi=1}|_{2}^{2}+\beta|\nabla p|_{\xi=1} |_{2}^{2}\notag\\
&\hspace{5mm}+C|p|_{\xi=1}|_{4}(|\bu_{\xi}+\partial_{\xi}Z_{1}|_{2}+|\bu+Z_{1}|_{2} )\|Z_{3}\|_{3}\notag\\
&\hspace{5mm}+C|p|_{\xi=1}|_{4}^{4}+C\|\partial_{\xi}Z_{1}+\bu_{\xi}\|_{1}^{2}+C\|Z_{1}+\bu\|_{1}^{2}\notag\\
&\hspace{5mm}+C|p|_{\xi=1}|_{2}(|Q_{q}|_{2}+|\partial_{\xi}Q_{q}|_{2}+|Z_{3}|_{2}+|\partial_{\xi}Z_{3}|_{2})\notag\\
&\leq\frac{\beta}{2} \partial_{t}|p|_{\xi=1}|_{2}^{2}+\beta|\nabla p|_{\xi=1} |_{2}^{2}\notag\\
&\hspace{5mm}+|p|_{\xi=1}|_{4}^{2}+C(|\bu_{\xi}|_{2}^{2}\|Z_{3}\|_{3}^{2}+|\bu|_{2}^{2}\|Z_{3}\|_{3}^{2}+\|Z_{1}\|_{1}^{2}\|Z_{3}\|_{3}^{2} )\notag\\
&\hspace{5mm}+C|p|_{\xi=1}|_{4}^{4}+C\|Z_{1}\|_{2}^{2}+C\|\bu_{\xi}\|_{1}^{2}+C\|\bu\|_{1}^{2}\notag\\
&\hspace{5mm}+C(|Q_{q}|_{2}^{2}+ |\partial_{\xi}Q_{q}|_{2}^{2}+\|Z_{3}\|_{1}^{2}  ).
\end{align}
By the H\"older inequality, the interpolation inequality and Young's inequality, we have
\begin{align}
\label{eqn:5.68}
&\hspace{5mm}-\int_{\D}\nabla_{\partial_{\xi}Z_{1}+\bu_{\xi}}(Z_{3}+p)p_{\xi}d\D\notag\\
&\leq |\nabla Z_{3}+\nabla p|_{2}|p_{\xi}|_{4}|\partial_{\xi}Z_{3}+\bu_{\xi}|_{4}\notag\\
&\leq |\nabla Z_{3}+\nabla p|_{2}|p_{\xi}|_{2}^{1/4}(|\nabla p_{\xi}|_{2}^{3/4}+|p_{\xi\xi}|_{2}^{3/4}+|p_{\xi}|_{2}^{3/4} )|\partial_{\xi} Z_{3}+\bu_{\xi}|_{4}\notag\\
&\leq \varepsilon ( |\nabla p_{\xi}|_{2}^{2}+|p_{\xi\xi}|_{2}^{2}  )+C|\nabla Z_{3}+\nabla p |_{2}^{8/5}|\partial_{\xi} Z_{3}+\bu_{\xi} |_{4}^{8/5}|p_{\xi}|_{2}^{2/5}\notag\\
&\hspace{5mm}+C|\nabla Z_{3}+\nabla p|_{2}|p_{\xi}|_{2}|\partial_{\xi} Z_{3}+\bu_{\xi}|_{4}\notag\\
&\leq \varepsilon ( |\nabla p_{\xi}|_{2}^{2}+|p_{\xi\xi}|_{2}^{2}  )+ C(\| Z_{3}\|_{1}^{2}+\|p\|_{1}^{2})|p_{\xi}|_{2}^{2}\notag\\
&\hspace{5mm}+C\|Z_{3}\|_{2}^{2}+C\|\bu_{\xi}\|_{1}^{2}.
\end{align}
By integration by parts, the H\"older inequality, the interpolation inequality and Young's inequality, we get
\begin{align}
\label{eqn:5.69}
&\hspace{5mm}-\int_{\D}(\nabla_{Z_{1}}p_{\xi}+\nabla_{Z_{1}+\bu}\partial_{\xi}Z_{3}  )p_{\xi}d\D\notag\\
&\leq \frac{1}{2}\int_{\D}(\mathrm{div} Z_{1})p_{\xi}^{2}d\D+ |p_{\xi}|_{2}\|Z_{3}\|_{3}|Z_{1}+\bu|_{4}\notag\\
&\leq C\|Z_{1}\|_{1}|p_{\xi}|_{4}^{2}+C|p_{\xi}|_{2}\|Z_{3}\|_{3}(\|Z_{1}\|_{1}+|\wt\bu|_{4}+\|\bar\bu\|_{1})\notag\\
&\leq C\|Z_{1}\|_{1}|p_{\xi}|_{2}^{1/2}(|\nabla p_{\xi}|_{2}^{3/2}+|p_{\xi\xi}|_{2}^{3/2}+|p_{\xi}|_{2}^{3/2}  )\notag\\
&\hspace{5mm}+C|p_{\xi}|_{2}\|Z_{3}\|_{3}(\|Z_{1}\|_{1}+|\wt\bu|_{4}+\|\bar\bu\|_{1})\notag\\
&\leq \varepsilon (|\nabla p_{\xi}|_{2}^{2}+|p_{\xi\xi}|_{2}^{2}    )+C(\|Z_{1}\|_{1}^{4}+\|Z_{3}\|_{3}^{2}+\|Z_{1}\|_{1})|p_{\xi}|_{2}^{2}\notag\\
&\hspace{5mm}+C\|Z_{1}\|_{1}^{2}+C|\wt\bu|_{4}^{2}+C\|\bar\bu\|_{1}^{2}.
\end{align}
By integration by parts, the H\"older inequality and Young's inequality we have
\begin{align}
\label{eqn:5.70}
&\hspace{5mm}\int_{\D} (\mathrm{div} \bu) \partial_{\xi} p \partial_{\xi} p d\D\notag\\
&=-2\int_{\D}\bu \nabla p_{\xi} p_{\xi}d\D\leq 2|\nabla p_{\xi}||\bu|_{4}|p_{\xi}|_{4}\notag\\
&\leq C|\nabla p_{\xi}||\bu|_{4}|p_{\xi}|_{2}^{1/4}(|\nabla p_{\xi}|_{2}^{3/4}+|p_{\xi\xi}|_{2}^{3/4}+|p_{\xi}|_{2}^{3/4}    )\notag\\
&\leq \varepsilon(|\nabla p_{\xi}|_{2}^{2}+|p_{\xi\xi}|_{2}^{2}   )+C(|\bu|_{4}^{2}+|\bu|_{4}^{8})|p_{\xi}|_{2}^{2}.
\end{align}
To estimate the third term on the righthand side of \eqref{eqn:5.58},
\begin{align}
\label{eqn:5.71}
&\hspace{5mm}\int_{\D}( \mathrm{div} (Z_{1}+\bu))\partial_{\xi}(Z_{3}+p)p_{\xi}d\D\notag\\
&=\int_{\D}( \mathrm{div} Z_{1})\partial_{\xi}(Z_{3}+p)p_{\xi}d\D+  \int_{\D}( \mathrm{div} \bu )(\partial_{\xi}Z_{3}) p_{\xi}d\D\notag\\
&\hspace{5mm}+\int_{\D} (\mathrm{div} \bu) \partial_{\xi} p \partial_{\xi} p d\D.
\end{align}
For the first two terms on the righthand side of  \eqref{eqn:5.71} we have
\begin{align}
\label{eqn:5.72}
&\hspace{5mm}\int_{\D}( \mathrm{div} Z_{1})\partial_{\xi}(Z_{3}+p)p_{\xi}d\D+  \int_{\D}( \mathrm{div} \bu )(\partial_{\xi}Z_{3}) p_{\xi}d\D\notag\\
&\leq |\mathrm{div} Z_{1}|_{\infty} |\partial_{\xi} Z_{3}|_{2}^{2}+|\mathrm{div} Z_{1}|_{\infty}|p_{\xi}|_{2}^{2}
+|\partial_{\xi}Z_{3}|_{\infty} |\mathrm{div} \bu|_{2}|p_{\xi}|_{2}\notag\\
&\leq C\|Z_{1}\|_{3}\|Z_{3}\|_{1}^{2}+C(\|Z_{1}\|_{3}+\|Z_{3}\|_{3}^{2})|p_{\xi}|_{2}^{2}+C\|\bu\|_{1}^{2}.
\end{align}
From \eqref{eqn:5.70}-\eqref{eqn:5.72}, we conclude that
\begin{align}
\label{eqn:5.73}
&\hspace{5mm}\int_{\D}( \mathrm{div} (Z_{1}+\bu))\partial_{\xi}(Z_{3}+p)p_{\xi}d\D\notag\\
&\leq \varepsilon(|\nabla p_{\xi}|_{2}^{2}+|p_{\xi\xi}|_{2}^{2}   )+C\|\bu\|_{1}^{2}\notag\\
&\hspace{5mm}+C(|\bu|_{4}^{2}+|\bu|_{4}^{8}+\|Z_{1}\|_{3} +\|Z_{3}\|_{3}^{2})|p_{\xi}|_{2}^{2}.
\end{align}
By integration by parts, the Sobolev imbedding theorem we get
\begin{align}
&\hspace{5mm}-\int_{\D}w(Z_{1})p_{\xi\xi}p_{\xi}d\D\notag\\
&=-\frac{1}{2}\int_{\D}w(Z_{1})\partial_{\xi}p_{\xi}^{2}=\frac{1}{2}\int_{\D}(\mathrm{div }Z_{1})p_{\xi}^{2} d\D\notag\\
&\leq\frac{1}{2}|\mathrm{div} Z_{1}|_{\infty}|p_{\xi}|_{2}^{2}\leq C\|Z_{1}\|_{3}|p_{\xi}|_{2}^{2}.
\end{align}
Similarly,  we have
\begin{align}
&\hspace{5mm}-\int_{\D}w(Z_{1}+\bu)\partial_{\xi\xi}Z_{3} p_{\xi}d\D\notag\\
&=\int_{\D}w(Z_{1}+\bu)\partial_{\xi}Z_{3}p_{\xi\xi}d\D+\int_{\D}\mathrm{div} (Z_{1}+\bu)\partial_{\xi}Z_{3}p_{\xi}d\D\notag\\
&\leq |\partial_{\xi}Z_{3}|_{\infty}(|p_{\xi\xi}|_{2}+|p_{\xi}|_{2})(\|Z_{1}\|_{1}+\|\bu\|_{1})\notag\\
&\leq \varepsilon |p_{\xi\xi}|_{2}^{2}+C|p_{\xi}|_{2}^{2}+C\|Z_{3}\|_{3}^{2}(\|Z_{1}\|_{1}^{2}+\|\bu\|_{1}^{2}).
\end{align}
In view of the H\"older inequality and Young's inequality we otain
\begin{equation}
\label{eqn:5.76}
\int_{\D}(\partial_{\xi}Q_{q}+\gamma \partial_{\xi}Z_{3})p_{\xi}d\D\leq C|\partial_{\xi} Q_{q}|_{2}^{2}+C\|Z_{3}\|_{1}^{2}+C|p_{\xi}|_{2}^{2}.
\end{equation}
By virtue of \eqref{eqn:5.67}-\eqref{eqn:5.69} and \eqref{eqn:5.73}-\eqref{eqn:5.76}, we conclude
\begin{align}
\label{eqn:5.77}
&\partial_{t}(|p_{\xi}|_{2}^{2}+|p(\xi=1)|_{2}^{2} )+|\nabla p_{\xi}|_{2}^{2}+|p_{\xi\xi}|_{2}^{2}+|\nabla p|_{\xi=1}|_{2}^{2}\notag\\
&\leq C(|p|_{\xi=1}|_{4}^{2}+|p|_{\xi=1}|_{4}^{4} )+C(|Q_{q}|_{2}^{2}+|\partial_{\xi}Q_{q}|_{2}^{2})\notag\\
&\hspace{5mm}+C(\|Z_{3}\|_{3}^{2}+1 )\|\bu\|_{1}^{2}+C\|Z_{1}\|_{1}^{2}\|Z_{3}\|_{3}^{2}\notag\\
&\hspace{5mm}+C\|Z_{1}\|_{2}^{2}+C\|\bu_{\xi}\|_{1}^{2}+C\|Z_{3}\|_{2}^{2}\notag\\
&\hspace{5mm}+C(1+\|Z_{3}\|_{3}^{2}+\|Z_{1}\|_{3}^{4}+\|p\|_{1}^{2}+\|\bar \bu\|_{1}^{8}+|\wt\bu|_{4}^{8}  )|p_{\xi}|_{2}^{2}.
\end{align}
Applying the Gronwall inequality to \eqref{eqn:5.77} yields that
\begin{align}
\label{eqn:5.78}
&\sup\limits_{t\in [t_0, \tau)}(|p_{\xi}|_{2}^{2}+|p(\xi=1)|_{2}^{2} )+\int_{t_0}^{\tau}( |\nabla p_{\xi}(t)|_{2}^{2}+|p_{\xi\xi}(t)|_{2}^{2}+|\nabla p|_{\xi=1}(t)|_{2}^{2})dt\notag\\
&\hspace{10mm}<C(\tau, Z_{1}, Z_{2}, Z_{3}, Q_{q}).
\end{align}
Taking the derivative of \eqref{eqn:S} with respect to $\xi,$ we obtain
\begin{align}
\label{eqn:5.79}
&\partial_{t} S_{\xi} +\nabla_{\partial_{\xi}Z_{1}+\bu_{\xi}}(Z_{2}+S)+\nabla_{Z_{1}+\bu}(\partial_{\xi}Z_{2}+S_{\xi} )\notag\\
&\hspace{5mm}-(\mathrm{div} (Z_{1}+\bu) )\partial_{\xi}(Z_{2}+S)+w(Z_{1}+\bu)\partial_{\xi\xi}(Z_{2}+S)\notag\\
&=\Delta S_{\xi}+\partial_{\xi\xi}S_{\xi}+\partial_{\xi}Q_{T}+\gamma \partial_{\xi}Z_{2}\notag\\
&\hspace{5mm}-\frac{br_{s}^{2}}{r^{2}}(1+a(Z_{3}+p))w(Z_{1}+\bu)\notag\\
&\hspace{5mm}+\frac{br_{s}}{r}a(\partial_{\xi}Z_{3}+p_{\xi})w(Z_{1}+\bu)\notag\\
&\hspace{5mm}-\frac{br_{s}}{r}(1+a(Z_{3}+p))\mathrm{div}(Z_{1}+\bu).
\end{align}
According to \eqref{eqn:5.79}, in order to estimate $S_{\xi}$ we only need to estimate the last three terms on the right hand side of \eqref{eqn:5.79} because other terms can be estimated in an similar way used in \eqref{eqn:5.78}. Applying Lemma \ref{lemma:2.3}, Young's inequality and the interpolation inequality to the following term yields that
\begin{align}
\label{eqn:5.80}
&\hspace{5mm}\int_{\D}\frac{-br_{s}^{2}}{r^{2}}(1+a(Z_{3}+p))w(Z_{1}+\bu)S_{\xi}d\D\notag\\
&\leq C\|Z_{1}\|_{1}^{2}+C\|\bu\|_{1}^{2}+C|S_{\xi}|_{2}^{2}\notag\\
&\hspace{5mm}+C|Z_{3}|_{\infty}^{2}(\|Z_{1}\|_{1}^{2}+\|\bu\|_{1}^{2})\notag\\
&\hspace{5mm}+C\|Z_{1}+\bu\|_{1}|p|_{4}|S_{\xi}|_{2}^{1/2}(|\nabla S_{\xi}|_{2}^{1/2}+|S_{\xi}|_{2}^{1/2} )\notag\\
&\leq \varepsilon |\nabla S_{\xi}|_{2}^{2}+C|S_{\xi}|_{2}^{2}+C\|Z_{1}\|_{1}^{2}+C\|Z_{1}\|_{1}^{2}\|Z_{3}\|_{2}^{2}\notag\\
&\hspace{5mm}+C(1+\|Z_{3}\|_{2}^{2}+|p|_{4}^{2})\|\bu\|_{1}^{2}+C\|Z_{1}\|_{1}^{2}|p|_{4}^{2}.
\end{align}
Using Lemma \ref{lemma:2.3}, the Sobolev embedding theorem and Young's inequality, we have
\begin{align}
\label{eqn:5.81}
&\hspace{5mm}\int_{\D}\frac{abr_{s}}{r}(\partial_{\xi}Z_{3}+p_{\xi})w(Z_{1}+\bu)S_{\xi}\notag\\
&\leq C\|Z_{1}+\bu\|_{1}|p_{\xi}|_{2}^{1/2}(|\nabla p_{\xi}|_{2}^{1/2}+|p_{\xi}|_{2}^{1/2}  )|S_{\xi}|_{2}^{1/2}(|\nabla S_{\xi}|_{2}^{1/2}+|S_{\xi}|_{2}^{1/2})\notag\\
&\hspace{5mm}+C|\partial_{\xi} Z_{3}|_{\infty}\|Z_{1}+\bu\|_{1}|S_{\xi}|_{2}\notag\\
&\leq \varepsilon |\nabla S_{\xi}|_{2}^{2}+C\|Z_{1}+\bu\|_{1}^{2}+C(|p_{\xi}|_{2}^{4}+|p_{\xi}|_{2}^{2}|\nabla p_{\xi}|_{2}^{2}+\|Z_{3}\|_{3}^{2}+1 )|S_{\xi}|_{2}^{2}.
\end{align}
By the H\"older inequality, the interpolation inequality and Young's inequality we get
\begin{align}
\label{eqn:5.82}
&\hspace{5mm}-\int_{\D}\frac{br_{s}}{r}(1+a(Z_{3}+p)) \mathrm{div} (Z_{1}+\bu)S_{\xi}d\D\notag\\
&\leq C\|Z_{1}+\bu\|_{1}|Z_{3}+p|_{4}|S_{\xi}|_{4}+C\|Z_{1}+\bu\|_{1}|S_{\xi}|_{2}\notag\\
&\leq C\|Z_{1}+\bu\|_{1}|Z_{3}+p|_{4}|S_{\xi}|_{2}^{1/4}(|\nabla S_{\xi}|_{2}^{3/4}+|S_{\xi\xi}|_{2}^{3/4}+|S_{\xi}|_{2}^{3/4} )\notag\\
&\hspace{5mm}+C\|Z_{1}+\bu\|_{1}^{2}+C|S_{\xi}|_{2}^{2}\notag\\
&\leq \varepsilon ( |\nabla S_{\xi}|_{2}^{2}+|S_{\xi\xi}|_{2}^{2})+C|S_{\xi}|_{2}^{2}+C\|Z_{1}+\bu\|_{1}^{2}(1+\|Z_{3}\|_{1}^{2}+|p|_{4}^{2}).
\end{align}
By virtue of \eqref{eqn:5.80}-\eqref{eqn:5.82}, taking an analogous argument of \eqref{eqn:5.77} yields
\begin{align}
\label{eqn:5.83}
&\partial_{t}(|S_{\xi}|_{2}^{2}+|S|_{\xi=1}|_{2}^{2} )+|\nabla S_{\xi}|_{2}^{2}+|S_{\xi\xi}|_{2}^{2}+|\nabla S|_{\xi=1}|_{2}^{2}\notag\\
&\leq C(|S|_{\xi=1}|_{4}^{2}+|S|_{\xi=1}|_{4}^{4} )+C(|Q_{T}|_{2}^{2}+|\partial_{\xi}Q_{T}|_{2}^{2})\notag\\
&\hspace{5mm}+C(\|Z_{2}\|_{3}^{2}+\|Z_{3}\|_{2}^{2}+|p|_{4}^{2} +1 )\|\bu\|_{1}^{2}+C\|Z_{1}\|_{1}^{2}(\|Z_{2}\|_{3}^{2}+\|Z_{3}\|_{2}^{2})\notag\\
&\hspace{5mm}+C\|Z_{1}\|_{2}^{2}(1+|p|_{4}^{2})+C\|\bu_{\xi}\|_{1}^{2}+C\|Z_{2}\|_{2}^{2}\notag\\
&\hspace{5mm}+C(1+\|Z_{2}\|_{3}^{2}+\|Z_{1}\|_{3}^{4}+\|Z_{3}\|_{3}^{2}+|p_{\xi}|_{2}^{4}\notag\\
&\hspace{5mm}+|p_{\xi}|_{2}^{2}|\nabla p_{\xi}|_{2}^{2}+\|S\|_{1}^{2}+\|\bar \bu\|_{1}^{8}+|\wt\bu|_{4}^{8}  )|S_{\xi}|_{2}^{2}.
\end{align}
Applying the Gronwall inequality to \eqref{eqn:5.83} yields that
\begin{align}
\label{eqn:5.84}
&\sup\limits_{t\in [t_0, \tau)}(|S_{\xi}|_{2}^{2}+|S|_{\xi=1}|_{2}^{2} )
+\int_{0}^{\tau}( |\nabla S_{\xi}(t)|_{2}^{2}+|S_{\xi\xi}(t)|_{2}^{2}+|\nabla S|_{\xi=1}(t)|_{2}^{2})dt\nonumber\\
&\hspace{5mm}< C(\tau,Q_T,Z_1,Z_2,Z_3,U_0).
\end{align}

		Now we take the inner product of (\ref{eqn:u}) with $-\Delta\bu$, and because $\int_{\D}(\frac{f}{R_0}\vec k\times \bu)\cdot\Delta\bu d\D=0$, $\int_{\D}\nabla\Phi_s\cdot\Delta\bu d\D=0$, we get
		\begin{align}
		\label{eqn:5.85}
			&\partial_{t}(|\nabla_{e_\theta}\bu|^2_2+|\nabla_{e_\varphi}\bu|^2_2+|\bu|_2^2)+|\Delta\bu|_2^2+|\nabla_{e_{\theta}} \bu_{\xi}|^2_2+ |\nabla_{e_{\varphi}} \bu_{\xi}|^2_2+ |\bu_{\xi}|^2_2\notag\\
			&=-\int_{\mathbf D}\nabla_{Z_1+\bu}(Z_1+\bu)\cdot\Delta\bu d\mathbf D-\int_{\mathbf D}w(Z_1+\bu)\partial_\xi(Z_1+\bu)\cdot\Delta\bu d\mathbf D\notag\\
			&\hspace{5mm}-\int_{\mathbf D}\frac{f}{R_0}\vec k\times Z_1\cdot\Delta\bu d\mathbf D+\gamma \int_{\mathbf D}Z_1\cdot\Delta\bu d\mathbf D\notag\\
&\hspace{5mm}-\int_{\mathbf D}\int_{\xi}^{1}\frac{br_s}{r}\nabla[(1+a(Z_{3}+p))(Z_{2}+S)] d\xi'\cdot \Delta\bu d\mathbf D.
		\end{align}
	To obtain estimates of \eqref{eqn:5.85}, we first consider the first term on the right hand side of \eqref{eqn:5.85}. Using the H\"older inequality, the interpolation inequality and Young's inequality yields that
\begin{align}
&\hspace{5mm}-\int_{\D}(\nabla_{Z_{1}+\bu}(Z_{1}+\bu))\cdot \Delta \bu d\D\notag\\
&\leq |\Delta \bu|_{2}|Z_{1}+\bu|_{4}(|\nabla_{e_{\theta}}Z_{1}|_{4}+ |\nabla_{e_{\theta}}\bu|_{4}+ |\nabla_{e_{\varphi}}Z_{1}|_{4}+   |\nabla_{e_{\varphi}}\bu|_{4} )\notag\\
&\leq |\Delta \bu|_{2}|Z_{1}+\bu|_{4}\bigg(\|Z_{1}\|_{2}+ |\nabla_{e_{\theta}}\bu|_{2}^{1/2}(|\Delta \bu|_{2}^{3/2}+|\nabla_{e_{\theta}}\bu_{\xi}|_{2}^{3/2}+|\nabla_{e_{\theta}}\bu|_{2}^{3/2} )\notag\\
&\hspace{5mm}+ |\nabla_{e_{\varphi}}\bu|_{2}^{1/2}(|\Delta \bu|_{2}^{3/2}+|\nabla_{e_{\varphi}}\bu_{\xi}|_{2}^{3/2}+|\nabla_{e_{\varphi}}\bu|_{2}^{3/2} ) \bigg)\notag\\
&\leq \varepsilon (|\Delta \bu|_{2}^{2}+ |\nabla_{e_{\theta}} \bu_{\xi}|_{2}^{2}+|\nabla_{e_{\varphi}} \bu_{\xi}|_{2}^{2})\notag\\
&\hspace{5mm}+C(|Z_{1}+\bu|_{4}^{8}+1)(|\nabla_{e_{\theta}} \bu|_{2}^{2}+  |\nabla_{e_{\varphi}} \bu|_{2}^{2})+C|Z_{1}+\bu|_{4}^{2}\|Z_{1}\|_{2}^{2}.
\end{align}
Applying the Lemma \ref{lemma:2.2} to the second term  on the right hand side of \eqref{eqn:5.85}, we have
\begin{align}
\label{eqn:5.87}
&\hspace{5mm}-\int_{\D}w(Z_{1}+\bu)\partial_{\xi}(Z_{1}+\bu)\cdot\Delta \bu d\D\notag\\
&\leq C|\Delta \bu|_{2}|\partial_{\xi}(Z_{1}+\bu) |_{2}^{1/2}(|\nabla_{e_{\theta}}\partial_{\xi}(Z_{1}+\bu) |_{2}^{1/2}\notag\\
&\hspace{5mm}+ |\nabla_{e_{\varphi}}\partial_{\xi}(Z_{1}+\bu) |_{2}^{1/2}+|\partial_{\xi}(Z_{1}+\bu)|_{2}^{1/2})\notag\\
&\hspace{5mm}\times |\mathrm{div} (Z_{1}+\bu)|_{2}^{1/2}(|\Delta Z_{1}+\Delta \bu|_{2}^{1/2}+ |\mathrm{div} (Z_{1}+\bu)|_{2}^{1/2})\notag\\
&\leq  \varepsilon |\Delta \bu|_{2}^{2}+C(1+|\partial_{\xi}Z_{1}|_{2}^{2}+|\bu_{\xi}|_{2}^{2} )(\|Z_{1}\|_{2}^{2}+|\nabla_{e_{\theta}}\bu_{\xi} |_{2}^{2}\notag\\
&\hspace{5mm}+|\nabla_{e_{\varphi}}\bu_{\xi} |_{2}^{2}+|\bu_{\xi}|_{2}^{2})(\|Z_{1}\|_{2}^{2}+|\nabla_{e_{\theta}}\bu |_{2}^{2}+ |\nabla_{e_{\varphi}}\bu |_{2}^{2} ).
\end{align}
Similarly,
\begin{align}
\label{eqn:5.88}
&\hspace{5mm}-\int_{\mathbf D}\int_{\xi}^{1}\frac{br_s}{r}\nabla[(1+a(Z_{3}+p))(Z_{2}+S)] d\xi'\cdot \Delta\bu d\mathbf D\notag\\
&=-\int_{\mathbf D}\int_{\xi}^{1}\frac{abr_{s}}{r}(\nabla Z_{3}+\nabla p )(Z_{2}+S)d\xi'\cdot \Delta \bu d\D\notag\\
&\hspace{5mm}- \int_{\mathbf D}\int_{\xi}^{1}\frac{br_{s}}{r}(1+a(Z_{3}+p))(\nabla Z_{2}+\nabla S )d\xi'\cdot\Delta \bu d\D\notag\\
&\leq C|\Delta \bu|_{2}|Z_{2}+S|_{4}|\nabla Z_{3}+ \nabla p|_{2}^{1/2}|\Delta Z_{3}+\Delta p |_{2}^{1/2}\notag\\
&\hspace{5mm}+C|\Delta \bu|_{2}(1+|Z_{3}+p|_{4})|\nabla Z_{2}+\nabla S|_{2}^{1/2}|\Delta Z_{2}+\Delta S|_{2}^{1/2} \notag\\
&\leq \varepsilon( |\Delta \bu|_{2}^{2}+ |\Delta p|_{2}^{2}+ |\Delta S|_{2}^{2})+C(\|Z_{2}\|_{2}^{2}+ \|Z_{3}\|_{2}^{2})\notag\\
&\hspace{5mm}+C(|Z_{2}|_{4}^{4}+|S|_{4}^{4} )(|\nabla Z_{3}|_{2}^{2}+|\nabla p|_{2}^{2})\notag\\
&\hspace{5mm}+C(1+|Z_{3}|_{4}^{4}+|p|_{4}^{4})(|\nabla Z_{2}|_{2}^{2}+|\nabla S|_{2}^{2}).
		\end{align}
Combining \eqref{eqn:5.85}-\eqref{eqn:5.88} yields that
\begin{align}
\label{eqn:5.89}
&\partial_{t}(|\nabla_{e_\theta}\bu|^2_2+|\nabla_{e_\varphi}\bu|^2_2+|\bu|_2^2)+|\Delta\bu|_2^2+|\nabla_{e_{\theta}} \bu_{\xi}|^2_2+ |\nabla_{e_{\varphi}} \bu_{\xi}|^2_2+ |\bu_{\xi}|^2_2\notag\\
&\leq \varepsilon ( |\Delta \bu|_{2}^{2}+|\nabla_{e_{\theta}}\bu_{\xi}|_{2}^{2}+ |\nabla_{e_{\varphi}}\bu_{\xi}|_{2}^{2} + |\Delta p|_{2}^{2}+ |\Delta S|_{2}^{2})\notag\\
&\hspace{5mm}+C(|Z_{1}+\bu|_{4}^{8}+1)(|\nabla_{e_{\theta}} \bu|_{2}^{2}+  |\nabla_{e_{\varphi}} \bu|_{2}^{2})+C|Z_{1}+\bu|_{4}^{2}\|Z_{1}\|_{2}^{2}\notag\\
&\hspace{5mm}+C(1+|\partial_{\xi}Z_{1}|_{2}^{2}+|\bu_{\xi}|_{2}^{2} )(\|Z_{1}\|_{2}^{2}+|\nabla_{e_{\theta}}\bu_{\xi} |_{2}^{2}\notag\\
&\hspace{5mm}+|\nabla_{e_{\varphi}}\bu_{\xi} |_{2}^{2}+|\bu_{\xi}|_{2}^{2})(\|Z_{1}\|_{2}^{2}+|\nabla_{e_{\theta}}\bu |_{2}^{2}+ |\nabla_{e_{\varphi}}\bu |_{2}^{2} )\notag\\
&\hspace{5mm}+C(\|Z_{2}\|_{2}^{2}+ \|Z_{3}\|_{2}^{2})\notag\\
&\hspace{5mm}+C(|Z_{2}|_{4}^{4}+|S|_{4}^{4} )(|\nabla Z_{3}|_{2}^{2}+|\nabla p|_{2}^{2})\notag\\
&\hspace{5mm}+C(1+|Z_{3}|_{4}^{4}+|p|_{4}^{4})(|\nabla Z_{2}|_{2}^{2}+|\nabla S|_{2}^{2}).
\end{align}
Similarly, we have
\begin{align}
&\frac{1}{2}\partial_{t}|\nabla S|_{2}^{2}+|\Delta S|_{2}^{2}+|\nabla S_{\xi}|_{2}^{2}+\alpha|\nabla S|_{\xi=1}|_{2}^{2}\notag\\
&\leq\varepsilon(|\Delta S|_{2}^{2}+|\Delta \bu|_{2}^{2}  )+C(\|Z_{1}\|_{1}^{8}+|\bu|_{4}^{8})|\nabla S|_{2}^{2}\notag\\
&\hspace{5mm}+C\bigg( (\|Z_{2}\|_{1}^{2}+|S_{\xi}|_{2}^{2})(\|Z_{2}\|_{2}^{2}+|\nabla S_{\xi}|_{2}^{2})\notag\\
&\hspace{5mm}+(1+\|Z_{1}\|_{2}^{2}+\|Z_{3}\|_{1}^{4}+|p|_{4}^{4})\bigg)(|\nabla_{e_{\theta}}\bu|_{2}^{2}+ |\nabla_{e_{\varphi}}\bu|_{2}^{2})\notag\\
&\hspace{5mm}+C\|Z_{2}\|_{2}^{2}(\|Z_{1}\|_{1}^{2}+|\bu|_{4}^{2})+C\|Z_{1}\|_{1}^{2}(\|Z_{2}\|_{1}^{2}+|S_{\xi}|_{2}^{2})(\|Z_{2}\|_{2}^{2}+|\nabla S_{\xi}|_{2}^{2})\notag\\
&\hspace{5mm}+C(\|Z_{1}\|_{2}^{2}+1)(1+\|Z_{1}\|_{2}^{2}+\|Z_{3}\|_{1}^{4}+|p|_{4}^{4} )+C|Q_{T}|_{2}^{2}+C|Z_{2}|_{2}^{2},
\end{align}
and
\begin{align}
\label{eqn:5.91}
&\frac{1}{2}\partial_{t}|\nabla p|_{2}^{2}+|\Delta p|_{2}^{2}+|\nabla p_{\xi}|_{2}^{2}+\beta|\nabla S|_{\xi=1}|_{2}^{2}\notag\\
&\leq\varepsilon(|\Delta p|_{2}^{2}+|\Delta \bu|_{2}^{2}  )+C(\|Z_{1}\|_{1}^{8}+|\bu|_{4}^{8})|\nabla p|_{2}^{2}\notag\\
&\hspace{5mm}+C(\|Z_{3}\|_{1}^{2}+|p_{\xi}|_{2}^{2})(\|Z_{3}\|_{2}^{2}+|\nabla p_{\xi}|_{2}^{2})
(|\nabla_{e_{\theta}}\bu|_{2}^{2}+ |\nabla_{e_{\varphi}}\bu|_{2}^{2})\notag\\
&\hspace{5mm}+C\|Z_{3}\|_{2}^{2}(\|Z_{1}\|_{1}^{2}+|\bu|_{4}^{2})+C\|Z_{1}\|_{1}^{2}(\|Z_{3}\|_{1}^{2}+|p_{\xi}|_{2}^{2})(\|Z_{3}\|_{2}^{2}+|\nabla p_{\xi}|_{2}^{2})\notag\\
&\hspace{5mm}+C\|Z_{1}\|_{2}^{2}+C|Q_{q}|_{2}^{2}+C|Z_{3}|_{2}^{2}.
\end{align}
To estimate $(\bu, S, p)$ in $\mathcal{V}, $ we denote by
\begin{align*}
f:=|\nabla_{e_{\theta}}\bu|_{2}^{2}+|\nabla_{e_{\varphi}}\bu|_{2}^{2}+|\bu|_{2}^{2}+|\nabla S|_{2}^{2}+|\nabla p|_{2}^{2},
\end{align*}
\begin{align*}
g&:=|\Delta \bu|_{2}^{2}+|\nabla_{e_{\theta}}\bu_{\xi}|_{2}^{2}+|\nabla_{e_{\varphi}}\bu_{\xi}|_{2}^{2}\notag\\
&\hspace{5mm}+|\Delta S|_{2}^{2}+|\nabla S_{\xi}|_{2}^{2}+|\Delta p|_{2}^{2}+|\nabla p_{\xi}|_{2}^{2},
\end{align*}
\begin{align*}
h&:=1+|\bu|_{4}^{8}+|p|_{4}^{8}+|S|_{4}^{8}+\|Z_{1}\|_{2}^{8}+\|Z_{2}\|_{2}^{8}+\|Z_{3}\|_{2}^{8}\notag\\
&+(\|Z_{2}\|_{1}^{2}+|S_{\xi}|_{2}^{2})(\|Z_{2}\|_{2}^{2}+|\nabla S_{\xi}|_{2}^{2}) \notag\\
&+(\|Z_{3}\|_{1}^{2}+|p_{\xi}|_{2}^{2} )(\|Z_{3}\|_{2}^{2}+| \nabla p_{\xi}|_{2}^{2} )\notag\\
&+(1+\|Z_{1}\|_{1}^{2}+|\bu_{\xi}|_{2}^{2})(\|Z_{1}\|_{2}^{2}+\|\bu_{\xi}\|_{1}^{2}),
\end{align*}
and
\begin{align*}
k&:=(1+\|Z_{2}\|_{1}^{4}+|S|_{4}^{4}+\|Z_{3}\|_{1}^{4}+|p|_{4}^{4})(\|Z_{2}\|_{1}^{2}+\|Z_{3}\|_{1}^{2})\notag\\
&+(1+\|Z_{1}\|_{2}^{2})(1+\|Z_{1}\|_{2}^{2}+\|Z_{3}\|_{1}^{4}+|p|_{4}^{4})\notag\\
&+\|Z_{1}\|_{1}^{2}(\|Z_{3}\|_{1}^{2}+|p_{\xi}|_{2}^{2} )(\|Z_{3}\|_{2}^{2}+|\nabla p_{\xi}|_{2}^{2})\notag\\
&+\|Z_{1}\|_{1}^{2}(\|Z_{2}\|_{1}^{2}+|S_{\xi}|_{2}^{2} )(\|Z_{2}\|_{2}^{2}+|\nabla S_{\xi}|_{2}^{2})\notag\\
&+|Q_{T}|_{2}^{2}+|Q_{q}|_{2}^{2}+|Z_{2}|_{2}^{2}+|Z_{3}|_{2}^{2}\notag\\
&+(\|Z_{2}\|_{2}^{2}+\|Z_{3}\|_{2}^{2})(\|Z_{1}\|_{1}^{2}+|\bu|_{4}^{2}).
\end{align*}
In view of \eqref{eqn:gronwallofL4p}, \eqref{eqn:5.24}, \eqref{eqn:L4ofwtu}, \eqref{eqn:L2ofgradientbaru}, \eqref{eqn:5.56}, \eqref{eqn:5.78} and \eqref{eqn:5.84}, we  have $f, h, k\in L^{2}([t_0,\tau]; \mathbb{R}).$ Combining \eqref{eqn:5.89}-\eqref{eqn:5.91}, we obtain that
\begin{equation}
\label{eqn:5.92}
\partial_{t}f(t)+g(t)\leq h(t)f(t)+k(t),
\end{equation}
for $t\geq t_0.$
By the Gronwall inequality,  we  have
		\begin{align}
		\label{eqn:5.93}
			\sup_{t\in[t_0,\tau)}f(t)+\int_{t_0}^\tau gdt \leq C(\tau,Q_T,Q_q,Z_1,Z_2,Z_3,U_0).
		\end{align}

\qed

$\mathbf{Acknowledgments}.$  The authors are deeply grateful for  Professor Zhen-Qing Chen for the valuable discussions and suggestions.

	
\end{document}